\documentclass[10pt,twoside,dvipsnames]{siamart1116}
\usepackage[english]{babel}
\usepackage{graphicx,epstopdf,epsfig}
\usepackage{amsfonts,epsfig,fancyhdr,graphics, hyperref,amsmath,amssymb}

\usepackage{colortbl}

\usepackage{mathtools}
\usepackage{hyperref}
\usepackage{cleveref}
\usepackage{comment}
\usepackage{ragged2e}
\usepackage{xcolor}%
\usepackage{textcomp}%
\usepackage{multirow}%
\usepackage{diagbox}
\usepackage{algorithm}%
\usepackage{algorithmic}%
\numberwithin{algorithm}{section}
\usepackage{listings}%

\renewtheorem{theorem}{Theorem}[section]
\renewtheorem{definition}[theorem]{Definition}
\newtheorem{property}[theorem]{Property}
\renewtheorem{lemma}[theorem]{Lemma}
\renewtheorem{proposition}[theorem]{Proposition}
\newtheorem{remark}[theorem]{Remark}
\renewtheorem{corollary}[theorem]{Corollary} 
\renewtheorem{proposition}[theorem]{Proposition}%  proposition etc.
\numberwithin{equation}{section}

\crefname{thm}{Theorem}{Theorems}
\crefname{lem}{Lemma}{Lemmas}
\crefname{property}{Property}{Properties}
\crefname{defi}{Definition}{Definitions}
\crefname{algorithm}{Algorithm}{Algorithms}
\crefname{rem}{Remark}{Remarks}
\crefname{corollary}{Corollary}{Corollaries}
\crefname{appendix}{Appendix}{Appendices}
\Crefname{appendix}{Appendix}{Appendices}
\crefname{figure}{Figure}{Figures}

\newcommand{\wordlabel}[2]{\hypertarget{#1}{#2}}
\newcommand{\wordref}[2]{\hyperlink{#1}{#2}}
\newcommand{\fl}{\mathrm{fl}}

\newcommand{\simcom}{\textcolor{blue}{SM:} \textcolor{blue}}
\newcommand{\new}[1]{#1}

\newcommand{\Names}{Simon Mataigne and Kyle A.~Gallivan}
\newcommand{\Title}{The Eigenvalue Decomposition of Normal Matrices by the Skew-Symmetric Part}

\begin{document}

\bibliographystyle{plainurl}

%  Leave these commented lines here
%\input{ELAheader-template.tex}
% ELA insert correct page number
\setcounter{page}{1}

\thispagestyle{empty}

%Insert the title of the paper
 \title{\Title}%\thanks{Received by the editors on \DoS. Accepted for publication on \DoA.  Handling Editor: \HE. Corresponding Author: \CA}}

\author{
Simon Mataigne\thanks{ICTEAM Institute,
UCLouvain, Louvain-la-Neuve, Belgium.
(simon.mataigne@uclouvain.be). Simon Mataigne is a Research Fellow of the Fonds de la Recherche Scientifique - FNRS. This work was supported by the Fonds de la Recherche Scientifique - FNRS under Grant no T.0001.23.}
% Remember to put \and between any two authors
\and
Kyle A.~Gallivan\thanks{Department of Mathematics, Florida State University, Tallahasse, FL, USA. (kgallivan@fsu.edu). This work was initiated while Kyle~A.~Gallivan was a Visiting Professor at the ICTEAM Institute, UCLouvain.}}

\markboth{\Names}{\Title}

\maketitle

\begin{abstract}
We propose a new method for computing the eigenvalue decomposition of a dense real normal matrix $A$ through the decomposition of its skew-symmetric part. The method relies on algorithms that are known to be efficiently implemented, such as the bidiagonal singular value decomposition and the symmetric eigenvalue decomposition. The advantages of this method stand for normal matrices with few real eigenvalues, such as random orthogonal matrices. We provide a stability and a complexity analysis of the method. The numerical performance is compared with existing algorithms. In most cases, the method has the same operation count as the Hessenberg factorization of a dense matrix. Finally, we provide experiments for the application of computing a Riemannian barycenter on the special orthogonal group.
\end{abstract}

\begin{keywords}
Normal matrices, eigenvalue decomposition, orthogonal matrices, Schur decomposition, matrix logarithm, Riemannian barycenter, Karcher mean.
\end{keywords}
\begin{AMS}
65F15, 15B10, 15B57, 62R30.
\end{AMS}

% Sample article for the Electronic Journal of Linear Algebra

%%%%%%%%%%%%%%%%%%%%%%%%%%%%%%%%%%%%%%%

%%%%%%%%%%%%%%%%%%%%%%%%%%%%%%%%%%%%%%%%%%%%%%%%%%%%%%%%%%%%%

\section{Introduction}
\renewcommand\labelitemi{$\vcenter{\hbox{\tiny$\bullet$}}$}
Computing the eigenvalues and the eigenvectors of a matrix is one of the most fundamental problems in linear algebra. It is well-known that certain classes of matrices allow fast algorithms with convergence and accuracy guarantees. For example, Hermitian matrices ($A=\overline{A^\top }\eqcolon A^*$) enjoy a wide scope of applications and it is known that the QR algorithm \cite{Francis1961TheQT,Francis1961TheQT2} converges globally \cite{Wilkinson1965,Wilkinsonshifts,WILKINSON1968409}.
Additionally, for Hermitian matrices, a divide and conquer method and  a multisectioning method were designed~\cite{Cuppen80,DongarraSorensen87, Syshin87}. Robust symmetric eigensolvers are available in~\href{https://www.netlib.org/lapack/explore-html/index.html}{\texttt{LAPACK}}~\cite{lapack99,dhillon1997new}. All these methods consist of reducing the Hermitian matrix to a tridiagonal form under unitary ($UU^*=U^*U=I$) similarity transformations \cite{Paige72,Parlett79,SIMON1984,Wilkinson62}, and then solving the tridiagonal eigenvalue problem (EVP). For general matrices, the QR algorithm remains a state of the art approach to obtain the Schur form (upper triangular matrix with eigenvalues on the diagonal).

Although not nearly as popular as the symmetric EVP, the real skew-symmetric EVP ($A = -A^\top $) is known to be even more efficiently solvable \cite{Penke20,WardGray78} because, first, skew-symmetric matrices can also be reduced to a tridiagonal form \cite{HuangJia24,WardGray78}. Secondly, the zero diagonal makes the EVP equivalent to a bidiagonal singular value problem with half the size of the original EVP. It can be solved both rapidly and accurately~\cite{DemmelKahan88,DemmelKahan90, Kurukula94}. This procedure exploits that eigenvalues arise by pairs of conjugate purely imaginary numbers. Unfortunately, in the complex setting, skew-Hermitian matrices turn into Hermitian matrices when multiplied by the imaginary unit. This strategy cannot be extended; real-valued matrices are more structured.

The focus of this paper, normal matrices ($AA^*=A^*A$), includes the two previous classes. It was shown that the shifted QR algorithm converges for normal matrices~\cite{BATTERSON1994181}. However, normal matrices are diagonalizable under unitary similarity transformations, i.e., the Schur form is diagonal. Consequently, the costly QR iterations operate on a dense upper Hessenberg matrix while converging to a diagonal matrix without exploiting this property. The key observation for a better method is the following: the Hermitian part $H\coloneq\frac{A+A^*}{2}$ and the skew-Hermitian part $\Omega\coloneq\frac{A-A^*}{2}$ of a normal matrix $A$ commute: $H\Omega=\Omega H$. It is known that two normal matrices $H$ and $\Omega$ that commute are jointly diagonalizable, i.e., there exists a unitary matrix $U$ such that $U^*HU = \Lambda_H$ and $U^*\Omega U = \Lambda_\Omega$~\cite[Thm.~2.5.5]{horn13}, and thus $A=U(\Lambda_H+\Lambda_\Omega)U^*$. The latter decomposition can be obtained by the \emph{simultaneous diagonalization} of two Hermitian matrices (here $H$ and $i\Omega$) using the Jacobi method \cite{Gerstner98,Goldstine59}. These algorithms were notably developed for generalized EVPs \cite{CHARLIER1990,WATTERS1974}. Another approach, called the \emph{do-one-then-do-the-other} (DODO), considers the sequential diagonalization of $H$ and $\Omega$~\cite{Sutton23}. Moreover, a method leveraging random combinations of $H$ and $i\Omega$ has been proposed recently \cite{Haoze24}. Finally, the normal EVP can be transformed into a complex symmetric EVP, more details can be found in \cite{Aurentz2017ComputingTE,Ferranti2013,Vandebril2011}. However, all previously cited methods consider \emph{complex} matrices. For real matrices, a Jacobi method~\cite{ZhouBrent03} was designed. The aim of this paper is to show that the additional structure of real matrices enables a specific and computationally efficient approach.

Unitary matrices, called \emph{orthogonal matrices} when real-valued, are an important subset of normal matrices. A rich literature has emerged from the representation of $n\times n$ unitary Hessenberg matrices in their \emph{Schur parametric form}, i.e., as the product of $n$ elementary Householder matrices (or equivalently as Givens rotations): 
\begin{align*}
	 H = H_1H_2...H_n  \text{ where }&\\
	 H_k=\mathrm{diag}\left[I_{k-1},\begin{bmatrix}-\alpha_k&\beta_k\\\beta_k&\alpha_k^*
	\end{bmatrix}, I_{n-k-3} \right] \text{ with } |\alpha_k|^2 + &|\beta_k|^2=1 \text{ for } k=1,2,...,n-1,
\end{align*}
 and $H_n =\mathrm{diag}\left[I_{n-1},-\alpha_n\right]$ with $|\alpha_n|=1$. This compact storage format allows cheap and stable QR iterations, called UHQR, requiring only $\mathcal{O}(n)$ flops per iteration when only the eigenvalues are desired \cite{GRAGG19861, Stewart2006}. This method can be implemented as a bulge-chasing algorithm~\cite{Aurentz2015} and it was shown that appropriate shifts produce quadratic convergence~\cite{wanggragg02,wanggragg01}. The success of the Schur parametric form also led to the development of alternative algorithms such as divide-and-conquer \cite{graggreichel90,GuGuzzo2003} and SVD-based methods \cite{AmmarReichel86,CALVETTI2015}. The common step for all previously cited methods is the initial reduction to Hessenberg form. Our method avoids this computationally expensive step.% Finally, another approach leverages Hermitian matrices \cite{GEMIGNANI2005}.}%A rich literature has shown convergence of the upper Hessenberg form to the Schur form under QR iterations with shifts \cite{wanggragg02,wanggragg01}. Apart from QR, two other methods starting from the upper Hessenberg form have been proposed \cite{AmmarReichel86,graggreichel90}. The downside of the previously cited approaches is the need for initial reduction to the upper Hessenberg form. Consequently, the commutativity of the Hermitian and skew-Hermitian parts is not exploited.

%In view of the popularity of the QR algorithm, the standard way of obtaining the eigenvalue decomposition of a normal matrix, with orthogonal matrices as particular case, remains the classical path of first reduction to Hessenberg form, then perform QR iterations. For anyone performing this, it should be disturbing that the fact that the real Schur form is not block diagonal is not leveraged.
\paragraph{Contributions} In this paper, we show that real normal matrices, with orthogonal matrices as a special case,  admit a fast factorization to the eigenvalue decomposition (EVD) by exploiting the structure of their invariant spaces. This structure only exists in the real setting. In most cases,\footnote{The precise meaning of ``most cases'' is described in \cref{sec:best_worst_case}.} if a normal matrix $A\in\mathbb{R}^{n\times n}$ has $r$ real eigenvalues ($r$ does not need to be known in advance), its EVD can be obtained by solving the rank $n-r$ EVP of the skew-symmetric part of $A$, plus a $r\times r$ symmetric EVP. The method is designed for matrices where $r\ll n$, a frequent property of random orthogonal matrices (e.g., Haar-distributed). We provide a stability and a complexity analysis of the method and compare its numerical performance with existing algorithms. When all eigenvalues have distinct imaginary parts, i.e., well-separated in floating point arithmetic, the method has the same operation count as the Hessenberg factorization of a dense matrix. % if Householder reflectors are assembled, i.e., $\sim\frac{14}{3}n^3$ flops. 
When the imaginary parts of the eigenvalues are clustered, we propose refinement steps based on theoretical bounds that allow to maintain a user-specified accuracy of the decomposition. We corroborate our results in a suite of experiments. The method facilitates efficiently computing the matrix logarithm of orthogonal matrices, the main computational bottleneck for calculating Riemannian logarithms and barycenters on matrix manifolds. In particular in these applications, the EVP is an intermediate computation used in iterative algorithms where speed is more important than accuracy to machine precision. In addition, the normal EVP also finds applications in generalized EVPs \cite{BUNSEGERSTNER1991741,CHARLIER1990} and computation of Pisarenko frequency estimates \cite{Ammar1988determination, Cybenko85}.   
\paragraph{Overview} 
\Cref{sec:preliminaries} introduces the essential properties of the normal EVP and related concepts. \Cref{sec:normal_EVP} presents theoretical results that lead to the exact arithmetic formulation of our method for the normal EVP: \cref{alg:schur_decomposition}. In \cref{sec:accuracy}, we analyze the impact of floating-point arithmetic through a sensitivity analysis of the invariant spaces of the skew-symmetric part with respect to perturbations of the original matrix; this yields \cref{alg:schur_decomposition_floating}. \Cref{sec:exp_accuracy} corroborates the perturbation analysis through a series of experiments. In \cref{sec:best_worst_case}, we examine the computational complexity of the algorithm and validate the analysis with timing experiments. Finally, in \cref{sec:karcher_mean}, we apply \cref{alg:schur_decomposition_floating} to enhance the performance of Riemannian gradient descent for computing the Riemannian barycenter on $\mathrm{SO}(n)$.
\paragraph{Reproducibility} An implementation of the method, as well as the codes to perform all experiments presented in this paper, are available at \url{https://github.com/smataigne/NormalEVP.jl}.
\paragraph{Notation} We write $A=[a_{ij}]_{i,j=1}^n$ to denote matrix entries. We use $\mathrm{sym}(A)\coloneq\frac{A+A^\top }{2}$ and $\mathrm{skew}(A)\coloneq\frac{A-A^\top }{2}$ for, respectively, the symmetric part and the skew-symmetric part of the matrix $A$. $A_{a:b}$ denotes the selection of \emph{columns} indexed by $a$ to $b$ included and $A_{a:b,c:d}$ refers to the submatrix of \emph{rows} $a$ to $b$ and \emph{columns} $c$ to $d$ included. The Frobenius norm and the spectral norm are denoted $\|\cdot\|_\mathrm{F}$ and $\|\cdot\|_2$, respectively. The identity matrix of size $n\times n$ is denoted $I_n$. Since orthogonal transformations are central to normal EVPs, we recall the notation for the orthogonal group $\mathrm{O}(n)$ defined as 
\begin{equation}
	\nonumber
	\mathrm{O}(n) \coloneq\{Q\in\mathbb{R}^{n\times n}\ | Q^\top Q=QQ^\top =I_n\}.
\end{equation} 
$\mathrm{O}(n)$ admits two connected components: matrices with determinant $1$ and matrices with determinant $-1$, sometimes respectively referred to as the set of rotations and reflections. Rotations are best known as the special orthogonal group $\mathrm{SO}(n)\coloneq\{Q\in \mathrm{O}(n)\ | \det(Q)=1\}$. Finally, the set of $n\times p$ matrices, $n\geq p$, with orthonormal columns is called the Stiefel manifold $\mathrm{St}(n,p)$, namely
\begin{equation}
	\nonumber
	\mathrm{St}(n,p)\coloneq\{V\in\mathbb{R}^{n\times p}\ | \ V^\top V=I_p\}.
\end{equation}
In the special case $n=p$, we have $\mathrm{St}(n,n)= \mathrm{O}(n)$.

\section{Preliminaries}\label{sec:preliminaries}
\subsection{The normal eigenvalue problem}
In the context of \emph{real} normal matrices, eigenvalues and eigenvectors are, in general, complex-valued. Nonetheless, it is preferred for most eigenvalue algorithms to compute with real arithmetic. % since complex operations are more expensive.\footnote{Indeed, performing $(a+ib)(c+id)=(ac-bd)+(ad + cb)i$ requires $4$ real multiplications and $2$ additions, thus approximately $6$ times the cost of one real multiplication.} 
For a normal matrix~$A$, obtaining the real Schur decomposition (RSD), i.e., $A=QSQ^\top $ with $Q\in \mathrm{O}(n)$, is equivalent to obtaining the EVD. Indeed, the real Schur form $S$ is block-diagonal with $2\times 2$ and/or $1\times 1$ blocks \cite[Thm~2.5.8]{horn13}. Each $2\times 2$ block has the form $\lambda\left[\begin{smallmatrix}
	\cos(\theta)&-\sin(\theta)\\
	\sin(\theta)&\cos(\theta)
	\end{smallmatrix}\right]$, $\lambda>0,\ \theta\in(0,\pi)$, and is associated with the conjugate pair of eigenvalues $\{\lambda e^{\pm i\theta}\}$. The transition from the RSD to the EVD uses for each $2\times2$ block the relation
\begin{equation}\label{eq:RSDtoEVD}
	\lambda\begin{bmatrix}
	\cos(\theta)&-\sin(\theta)\\
	\sin(\theta)&\cos(\theta)
	\end{bmatrix} = \frac{1}{2}\begin{bmatrix}
	1&1\\
	i&-i
	\end{bmatrix}\begin{bmatrix}
	\lambda e^{-i\theta}&0\\
	0&\lambda e^{i\theta}
	\end{bmatrix}\begin{bmatrix}
	1&-i\\
	1&i
	\end{bmatrix}.
\end{equation}
In this paper, we consider a permuted version of the real Schur form~$S$, introduced in~\cref{property:schur}. The latter permuted version allows writing the matrix transformations easily and succinctly. Therefore, we will refer to \eqref{eq:schur_form} as the ``real Schur decomposition'', although it differs from the usual block-diagonal definition by a permutation.
\begin{property}[{\cite[Thm~2.5.8]{horn13}}]\label{property:schur}
Every normal matrix $A\in\mathbb{R}^{n\times n}$ is similar to a real Schur form $S\in\mathbb{R}^{n\times n}$ under orthogonal  transformation by a matrix $Q\in\mathrm{O}(n)$ such that
	\begin{equation}\label{eq:schur_form}
		A = Q S Q^\top  = Q \begin{bmatrix}
	\Lambda\cos(\Theta)&0&-\Lambda\sin(\Theta)\\
	0&\breve{\Lambda}&0\\
	 \Lambda\sin(\Theta)&0&\Lambda\cos(\Theta)
	\end{bmatrix}Q^\top ,
	\end{equation} where $\breve{\Lambda}=\mathrm{diag}(\breve{\lambda}_1,...,\breve{\lambda}_r)$ contains the $r$ real eigenvalues of $A$. Let $p\coloneq\frac{n-r}{2}$,\footnote{Note that $p$ is always an integer since $r$ is odd if and only if $n$ is odd.} then $\Lambda=\mathrm{diag}(\lambda_1,...,\lambda_p)$ and $\Theta=\mathrm{diag}(\theta_1,...,\theta_p)$ are real, can be assumed positive with $\theta_j\in(0,\pi)$ w.l.o.g.~and are associated with the pairs of eigenvalues~$\{\lambda_j e^{\pm i\theta_j}\}$ for $j=1,...,p$. The columns of the matrix $Q$ are called \emph{Schur vectors}.
	Special case~1: if $A$ is symmetric, $r=n$.
	Special case~2: if $A$ is skew-symmetric, $\cos(\Theta)=0$ and $\breve{\Lambda}=0_r$ \cite[Cor.~2.5.11]{horn13}.
\end{property}
 \begin{remark}
\emph{The Schur vectors of $A$ associated with real eigenvalues are also eigenvectors of $A$. Moreover, in view of \eqref{eq:RSDtoEVD}, if $q_1,q_2\in\mathrm{St}(n,1)$ is a pair of orthogonal Schur vectors associated with the same $(\lambda,\theta)$-block, then $\mathrm{span}([q_1\ q_2])$ is an invariant subspace of $A$ and $\frac{\sqrt{2}}{2}(q_1+iq_2)$ and $\frac{\sqrt{2}}{2}(q_1-iq_2)$ are eigenvectors of $A$.}
 \end{remark}
\subsection{The orthogonal groups}
The Lie groups $\mathrm{O}(n)$ and $\mathrm{SO}(n)$ are particularly important subsets of normal matrices for applications \cite{AbsMahSep2008,Chakraborty2019StatisticsOT,James1954}. Computing the RSD is a key primitive to obtain the matrix logarithm, which is required for computing Riemannian logarithms \cite{mataigne2024,RentmeestersQ,Szwagier23,Zimmermann17,ZimmermannHuper22} and Riemannian barycenters \cite{BINI2013,lirias165866,krakowski07,LIM2012,YING2016,ZimmermannHuper22}, notably on the Stiefel manifold. The principal matrix logarithm can be efficiently computed via the relation 
\begin{equation}\label{eq:log_normal}
	\log(A)=Q\log(S)Q^\top =Q\begin{bmatrix}
	\log(\Lambda)&0&-\Theta\\
	0&\log(\breve{\Lambda})&0\\
	 \Theta&0&\log(\Lambda)
	\end{bmatrix}Q^\top .
\end{equation}
 For both $\mathrm{O}(n)$ and $\mathrm{SO}(n)$, the property that eigenvalues lie on the unit circle yields $\Lambda = I_p$ and $\breve{\Lambda}=\mathrm{diag}(\pm 1,...,\pm 1)$. In particular, we have $\log(\Lambda)=0$. Note that $\log(\breve{\Lambda})$ admits a real element if and only if negative eigenvalues arise by pairs since $\left[\begin{smallmatrix}0&-\pi\\
 \pi&0\end{smallmatrix}\right]\in\log\left(\left[\begin{smallmatrix}-1&0\\
 0&-1\end{smallmatrix}\right]\right)$.

% If $n$ is odd, we know that $1$ is an eigenvalue and thus, $r\geq 1$. An interesting case for statistics on manifolds is to consider random orthogonal matrices. For such matrices if $n$ is even, $r=0$ with probability $1$ (no real eigenvalues). If $n$ is odd $r=1$ ($\breve{\Lambda}=1$) with probability one. Consequently, finding the decomposition \eqref{eq:schur_form} essentially reduces to finding $\Theta$.
\subsection{Invariance groups} The decomposition from~\eqref{eq:schur_form} is often not unique. It is key to this paper to describe the admissible transformations of the Schur vectors $Q$ such that \eqref{eq:schur_form} remains satisfied. To this end, we introduce the notion of the \emph{invariance group} of a matrix in~\cref{def:oic}.
\begin{definition}\label{def:oic}
Given a matrix $A\in\mathbb{R}^{n\times n}$, the \emph{invariance group} of $A$, written $\mathrm{ig}(A)$, is defined as the set
	\begin{equation}
	\label{eq:invariance_group}
	\mathrm{ig}(A)\coloneq\{Q\in\mathrm{O}(n)\ | \ QAQ^\top  = A\}.
	\end{equation}
	\end{definition}
	An equivalent definition of $\mathrm{ig}(A)$ is the set of orthogonal matrices that commute with~$A$. It is easily verified that $\mathrm{ig}(A)$ is a group under matrix multiplication since for all $Q_1,Q_2\in\mathrm{ig}(A)$, we have $Q_1Q_2\in\mathrm{ig}(A)$. The set $\mathrm{ig}\left(\left[\begin{smallmatrix}0&-I_m\\
	I_m&0\end{smallmatrix}\right]\right)$ is well-known as the \emph{ortho-symplectic group} $\mathrm{OSp}(2m)$, see, e.g., \cite{dopico09}.

Let $A$ be a normal matrix and $A=QSQ^\top$ be a RSD. For every $R\in\mathrm{ig}(S)$, we have \begin{equation*}
	A =(QR)(R^\top SR)(QR)^\top = (QR)S(QR)^\top.
\end{equation*}
Therefore, $QR$ is a also a matrix of Schur vectors. If $r<n$ ($A$ is not symmetric) or if $\breve{\Lambda}$ has repeated diagonal entries, then $\mathrm{ig}(S)$ is uncountable. Indeed, for each $2\times 2$ $(\lambda,\theta)$-block of $S$ and each real eigenvalue~$\breve{\lambda}$ of multiplicity $m>1$, we respectively have
\begin{equation}
\nonumber
\mathrm{ig}\left(\lambda\begin{bmatrix}
\cos(\theta)&-\sin(\theta)\\
\sin(\theta)&\cos(\theta)
\end{bmatrix}\right) = \mathrm{SO}(2)\quad\text{and,}\quad\mathrm{ig}(\breve{\lambda} I_m)=\mathrm{O}(m).
\end{equation}
Finally, we highlight a property of invariance groups that is important in \cref{sec:normal_EVP}.
\begin{lemma}\label{lem:oic_linearity}
For all pairs of matrices $A,B\in\mathbb{R}^{n\times n}$, we have $$\mathrm{ig}(A)\cap \mathrm{ig}(B)\subseteq \mathrm{ig}(A+B). $$
\end{lemma}
\begin{proof}
Let $Q\in \mathrm{ig}(A)\cap \mathrm{ig}(B)$. Therefore, we have $QAQ^\top =A$ and $QBQ^\top =B$, and by summation, $Q(A+B)Q^\top =A+B$. Hence, $Q\in\mathrm{ig}(A+B)$.
\end{proof}

\subsection{The even-odd permutation} To conclude the preliminaries, we introduce the notions of even-odd permutation and bidiagonalization of a tridiagonal skew-symmetric matrix. These notions are essential to the skew-symmetric EVP and are used in \cref{sec:normal_EVP}.
\begin{definition}
An even-odd permutation matrix $P_{\mathrm{eo}}\in\mathbb{R}^{n\times n}$ is a permutation such that
	for every vector $[a_1\ a_2\ ...\  a_n]\in\mathbb{R}^{1\times n}$, we have \begin{equation*}
	[a_1\ a_2\ ...\  a_n]P_{\mathrm{eo}} = [a_1\ a_3\ ...\ a_{\lceil\frac{n}{2}\rceil}\ a_2\ a_4\ ... \ a_{\lfloor\frac{n}{2}\rfloor}]. 
\end{equation*}
\end{definition}
The even-odd permutation is well-known in the field of discrete mathematics and it is related to the notion of Red-Black ordering \cite{TAO2023110008}. This permutation allows the reduction of the skew-symmetric EVP to a bidiagonal SVD problem by a similarity tranformation, described in \cref{prop:even_odd_T}.
\begin{property}[{\cite{WardGray78}}]\label{prop:even_odd_T}
Let $T\in\mathbb{R}^{n\times n}$ be a skew-symmetric tridiagonal matrix and $P_{\mathrm{eo}}\in\mathbb{R}^{n\times n}$ be an even-odd permutation matrix, then
\begin{equation}
	\label{eq:even_odd}
	P_{\mathrm{eo}}^\top  T P_{\mathrm{eo}} = \begin{bmatrix}
	0&-\widetilde{B}^\top \\
	\widetilde{B}&0
	\end{bmatrix},
\end{equation} 
where $\widetilde{B}\in \mathbb{R}^{\lfloor\frac{n}{2}\rfloor\times \lceil\frac{n}{2}\rceil}$ is bidiagonal.
\end{property}

\section{The eigenvalue decomposition of normal matrices}\label{sec:normal_EVP}
In this section, we describe our method in exact arithmetic for computing the RSD of a normal matrix~$A$. We show how the skew-symmetric part of a normal matrix reveals the structure of its Schur vectors and how these computations can be organized to design a fast RSD/EVD algorithm. In \cref{sec:accuracy}, we analyze the effect of floating point arithmetic. In particular, the method relies on the assumption that the three following routines are available:
\begin{itemize}
	\item \wordlabel{oracle1}{Routine~1}: returns the SVD of a real bidiagonal matrix.
	\item \wordlabel{oracle2}{Routine~2}: returns the RSD/EVD of a real symmetric matrix.
	\item \wordlabel{oracle3}{Routine~3}: returns the RSD of a (small) real normal matrix.
\end{itemize}
In the context of exact arithmetic, these routines can be viewed as oracles, whereas in practice, they are efficiently implemented, for example, in \href{https://www.netlib.org/lapack/explore-html/index.html}{\texttt{LAPACK}}. Indeed, \wordref{oracle1}{Routine~1} and \wordref{oracle2}{Routine~2} are respectively available as \href{https://www.netlib.org/lapack//explore-html/d6/d51/group__bdsqr_gade20fbf9c91aa7de0c3d565b39588dc5.html}{\texttt{bdsqr}}/\href{https://www.netlib.org/lapack/explore-html/d2/d37/group__bdsdc_ga343a26845ae796e1f95f0ff7e1b0af2d.html}{\texttt{bdsdc}} and \href{https://www.netlib.org/lapack/explore-html/d8/d1c/group__heev_ga8995c47a7578fef733189df3490258ff.html}{\texttt{syev}}/\href{https://www.netlib.org/lapack/explore-html/d8/d30/group__heevd_ga25b71a69f9921df0a6050aa5883f54f4.html}{\texttt{syevd}}/\href{https://www.netlib.org/lapack/explore-html/d1/d56/group__heevr_gaa334ac0c11113576db0fc37b7565e8b5.html}{\texttt{syevr}} routines. For the coherence of this discussion, the calls to \wordref{oracle3}{Routine~3} should be exceptional and only applied on small-size problems. \wordref{oracle3}{Routine~3} may be chosen to be \texttt{LAPACK}'s \href{https://www.netlib.org/lapack//explore-html/d5/d38/group__gees_gab48df0b5c60d7961190d868087f485bc.html}{\texttt{gees}} routine since the shifted QR algorithm is known to converge for normal matrices~\cite{BATTERSON1994181}. \wordref{oracle3}{Routine~3} may also be any of the previously cited methods for normal or orthogonal EVPs, e.g., \cite{AmmarReichel86,wanggragg01,ZhouBrent03}. All of these \texttt{LAPACK}'s routines are stable, see \cite{lapack99}.

\subsection{The skew-symmetric eigenvalue problem} The method starts by obtaining the RSD of $\Omega \coloneq \mathrm{skew}(A)$. It is known to be efficiently done in three steps \cite{WardGray78}. 

First, we compute $\Omega=\widetilde{Q}T\widetilde{Q}^\top $ with $\widetilde{Q}\in\mathrm{O}(n)$ and $T$ skew-symmetric and tridiagonal. This decomposition can be obtained in at least two ways, both adapted versions of symmetric tridiagonalization. Either by using Householder reflectors \cite{Wilkinson62}, see, e.g., \href{https://github.com/JuliaLinearAlgebra/SkewLinearAlgebra.jl/blob/main/src/hessenberg.jl}{SkewLinearAlgebra.jl}, or by using the Lanczos algorithm~\cite{HuangJia24}. Householder tridiagonalization is preferred for dense matrices for its stability and use of level 3 \texttt{BLAS} operations \cite{Bischof87,Schreiber89,Kazushige08}. Lanczos algorithm finds applications for large sparse matrices. 

Secondly, an even-odd permutation matrix $P_{\mathrm{eo}}$ yields 
\begin{equation*}%\label{eq:b_tilde}
	\Omega =\widetilde{Q}T\widetilde{Q}^\top = \widetilde{Q}P_{\mathrm{eo}}\begin{bmatrix}
	0&-\widetilde{B}^\top \\
	\widetilde{B}&0
	\end{bmatrix}P_{\mathrm{eo}}^\top \widetilde{Q}^\top,
\end{equation*}
where $\widetilde{B}\in\mathbb{R}^{\lfloor\frac{n}{2}\rfloor\times\lceil\frac{n}{2}\rceil}$ is bidiagonal. In the context of exact arithmetic, a transformation $G\in\mathrm{O}(n)$ can isolate the central square matrix $0_r$ in \eqref{eq:tridiagonalization}. $G$ is obtained by a finite process described in \url{https://github.com/smataigne/NormalEVP.jl}. In practice, we compute the SVD of $\widetilde{B}$ directly. For theoretical purposes, we consider the transformation $G$ in this section. This leads to
\begin{equation}\label{eq:tridiagonalization}
	\Omega = \widetilde{Q}P_{\mathrm{eo}}G\begin{bmatrix}
	0&0&-B^\top \\
	0&0_r&0\\
	B&0&0
	\end{bmatrix}G^\top P_{\mathrm{eo}}^\top \widetilde{Q}^\top ,
\end{equation}
where $B\in\mathbb{R}^{p\times p}$ is bidiagonal.   

The third and final step is a call to \wordref{oracle1}{Routine~1} to compute the SVD of the bidiagonal block $B =U\Sigma V^\top $. This yields
\begin{equation}\label{eq:skew_schur1}
	\Omega = \widetilde{Q}P_{\mathrm{eo}}G\begin{bmatrix}
	V&0&0\\
	0&I_r&0\\
	0&0&U
	\end{bmatrix}\begin{bmatrix}
	0&0&-\Sigma\\
	0&0_r&0\\
	 \Sigma&0&0
	\end{bmatrix}\begin{bmatrix}
	V^\top &0&0\\
	0&I_r&0\\
	0&0&U^\top 
	\end{bmatrix}G^\top P_{\mathrm{eo}}^\top  \widetilde{Q}^\top .
\end{equation}
By defining the matrix $\widehat{Q}\coloneq\widetilde{Q}P_{\mathrm{eo}}G\left[\begin{smallmatrix}
	V&0&0\\
	0&I_r&0\\
	0&0&U
	\end{smallmatrix}\right]$, the decomposition \eqref{eq:skew_schur1} becomes a RSD of $\Omega$. This decomposition can be compared with the skew-symmetric part of the (unknown) RSD of $A=QSQ^\top $ from~\eqref{eq:schur_form}:
\begin{equation}\label{eq:skew_schur2}
		\Omega = Q \mathrm{skew}(S)Q^\top  = Q \begin{bmatrix}
	0&0&-\Lambda\sin(\Theta)\\
	0&0_r&0\\
	 \Lambda\sin(\Theta)&0&0
	\end{bmatrix}Q^\top .
	\end{equation}
By comparing \eqref{eq:skew_schur1} and \eqref{eq:skew_schur2}, it must hold that $\Sigma=\Lambda\sin(\Theta)$ if we assume that both are sorted, e.g., in decreasing order. Unfortunately, it does not hold true for the Schur vectors: $Q\neq\widehat{Q}$ in general. However, by equating \eqref{eq:skew_schur1} and \eqref{eq:skew_schur2},  it holds that 
  \begin{equation*}%\label{eq:invariance_on_skew}
  	\widehat{Q}^\top Q\mathrm{skew}(S)Q^\top \widehat{Q} = \mathrm{skew}(S)\quad \Longrightarrow\quad \widehat{Q}^\top Q\in\mathrm{ig}(\mathrm{skew}(S)).
  \end{equation*}
We recall that $\widehat{Q}$ is a matrix of Schur vectors of $A$ if and only if $\widehat{Q}^\top Q\in\mathrm{ig}(S)$. We have thus not yet computed the Schur vectors of $A$ since, in general, $\mathrm{ig}(\mathrm{skew}(S))\nsubseteq \mathrm{ig}(S)$. Designing the next steps of the algorithm requires finding the relation between $\mathrm{ig}(\mathrm{skew}(S))$ and $\mathrm{ig}(S)$. We show in the next subsection that it depends on the number $r$ of real eigenvalues and the diagonal entries of $\Lambda\sin(\Theta)$.
\subsection{Invariant spaces of the skew-symmetric part: preparatory lemmas} The rest of \cref{sec:normal_EVP} is dedicated to the description of the relation between the invariance groups $\mathrm{ig}(\mathrm{skew}(S))$ and $\mathrm{ig}(S)$. By completing this analysis, we will be able to find a matrix $R\in\mathrm{O}(n)$ such that $ (\widehat{Q} R)^\top Q\in \mathrm{ig}(S)$, and thus $A = (\widehat{Q}R)S(\widehat{Q}R)^\top $. Let us start by dividing the eigenvalues of $A$ in three sets:
	\begin{itemize}
	\item \wordlabel{set1}{Set~1}: The eigenvalues with distinct nonzero imaginary parts.
	\item \wordlabel{set2}{Set~2}: The eigenvalues with repeated nonzero imaginary parts.
	\item \wordlabel{set3}{Set~3}: The real eigenvalues.
	\end{itemize}
	A preliminary result is showing that $R$ features a block structure related to the aforementioned sets of eigenvalues of $A$. This is done in \cref{lem:nullspace,lem:eigenspace_identified}. To this end, we partition $Q$ and $\widehat{Q}$ column-wise to match the blocks of the decomposition \eqref{eq:schur_form}:
	\begin{equation*}
		Q\coloneq\left[Q_{\mathrm{c},1}\ Q_\mathrm{r} \ Q_{\mathrm{c},2}\right]\quad\text{and,}\quad\widehat{Q}\coloneq[\widehat{Q}_{\mathrm{c},1}\ \widehat{Q}_\mathrm{r} \ \widehat{Q}_{\mathrm{c},2}],
	\end{equation*}
	where the notation follows because, by definition, $Q_\mathrm{r}$ contains the Schur vectors associated with the \emph{real} eigenvalues and $[Q_{\mathrm{c},1} \ Q_{\mathrm{c},2}]$ contains the Schur vectors associated with \emph{complex} eigenvalues. \cref{lem:nullspace} shows the equivalence between the real eigenspace of $A$ and the nullspace of $\mathrm{skew}(A)$.
\begin{lemma}\label{lem:nullspace}
Let $A\in\mathbb{R}^{n\times n}$ be a normal matrix with exactly $r>0$ real eigenvalues and $\breve{\Lambda}\in\mathbb{R}^{r\times r}$ be a diagonal matrix with these eigenvalues as diagonal entries. If a matrix $V\in \mathrm{St}(n,r)$ satisfies $AV = V\breve{\Lambda}$, then $\mathrm{skew}(A)V = 0$. Conversely, if $\mathrm{skew}(A)V = 0$, then there exists a matrix $\breve{R}\in\mathrm{O}(r)$ such that $A(V\breve{R}) = (V\breve{R})\breve{\Lambda}$.
	\end{lemma}
	\begin{proof}
	($\Longrightarrow$) Assume $AV = V\breve{\Lambda}$. Then, by definition, there is a RSD $A = QSQ^\top $ such that $V^\top  Q = [0_{r \times p} \ I_{r\times r}\ 0_{r\times p}]$. It follows that $\mathrm{skew}(A)V = Q\mathrm{skew}(S)Q^\top V = Q0=0$. 
	
	($\Longleftarrow$) Assume $\mathrm{skew}(A)V = 0$ and let $A=QSQ^\top $ be a RSD. Then it follows that $\ker(\mathrm{skew}(A)) = \mathrm{span}(Q_\mathrm{r})$. Moreover, since $\mathrm{skew}(A)V = 0$, $\mathrm{span}(V)\subseteq \ker(\mathrm{skew}(A))$.
%Then, it follows that $\mathrm{skew}(S)Q^\top V=0$. This yields
	%\begin{equation*}
	%\Lambda\sin(\Theta)[Q_{\mathrm{c},1}\ %Q_{\mathrm{c},2}]^\top V=0.
%\end{equation*}	 
%Since all diagonal entries of $\Lambda\sin(\Theta)$ are non-zero, it must hold that $[Q_{\mathrm{c},1}\ Q_{\mathrm{c},2}]^\top V=0$, i.e., the columns of $[Q_{\mathrm{c},1}\ Q_{\mathrm{c},2}]$ and $V$ span disjoint subspaces. 
Since both $Q_{\mathrm{r}}$ and $V$ have column rank~$r$, their columns span the same subspace of $\mathbb{R}^n$ and there is $\breve{R} \in \mathbb{R}^{r\times r}$ such that $Q_{\mathrm{r}}=V\breve{R}$. Finally, since $Q_{\mathrm{r}}\in\mathrm{St}(n,r)$, $I_r=Q_{\mathrm{r}}^\top Q_{\mathrm{r}}=\breve{R}^\top V^\top V\breve{R}=\breve{R}^\top \breve{R}$. Hence, $\breve{R}\in\mathrm{O}(r)$. 
	\end{proof}

 Since $\widehat{Q}$ satisfies $\mathrm{skew}(A)\widehat{Q}=\widehat{Q}\mathrm{skew}(S)$ by \eqref{eq:skew_schur1}, it holds in particular that $\widehat{Q}_\mathrm{r}$ spans the kernel of $\mathrm{skew}(A)$: $\mathrm{skew}(A)\widehat{Q}_\mathrm{r}=0$. Therefore, \cref{lem:nullspace} ensures the existence of $\breve{R}\in\mathrm{O}(r)$ such that $Q_\mathrm{r} = \widehat{Q}_\mathrm{r}\breve{R}$. As a corollary of \cref{lem:nullspace}, it also holds that
\begin{equation*}
Q_\mathrm{r}^\top [\widehat{Q}_{\mathrm{c},1} \ \widehat{Q}_{\mathrm{c},2} ] =0\quad \text{and,}\quad \widehat{Q}_\mathrm{r}^\top [Q_{\mathrm{c},1} \ Q_{\mathrm{c},2} ]=0.
\end{equation*}
Therefore, we can conclude that the matrix $R$ such that $(\widehat{Q}R)^\top Q \in \mathrm{ig}(S)$ has the block structure
\begin{equation}\label{eq:structure_M}
		R=\begin{bmatrix}
		R_{\mathrm{a}}&0&R_{\mathrm{b}}\\
		0&\breve{R}&0\\
		R_{\mathrm{c}}&0&R_{\mathrm{b}}\\
		\end{bmatrix}\in\mathrm{O}(n),
		\end{equation}
		 where $\left[\begin{smallmatrix}R_{\mathrm{a}}&R_{\mathrm{b}}\\R_{\mathrm{c}}&R_{\mathrm{d}}\end{smallmatrix}\right]\in\mathrm{O}(2p)$, $R_\mathrm{a}\in\mathbb{R}^{p\times p}$ and $\breve{R}\in\mathrm{O}(r)$. 
		 
We can show that $R$ is even more structured by considering \wordref{set1}{Set~1} and \wordref{set2}{Set~2} separately. \cref{lem:eigenspace_identified} shows that, for each set, $R$ acts on independent sets of columns of $\widehat{Q}$. 
\begin{lemma}\label{lem:eigenspace_identified}
For every normal matrix $A\in\mathbb{R}^{n\times n}$ that has $m\geq 1$ eigenvalues with imaginary part $\sigma>0$, it holds for every $V\in\mathrm{St}(n,2m)$ satisfying $\mathrm{skew}(A)V=V\left[\begin{smallmatrix}0&-\sigma I_m\\ \sigma I_m &0\end{smallmatrix}\right]$ that there exists a matrix $R\in\mathrm{O}(2m)$ such that \begin{equation*}
	A(VR) = (VR)\begin{bmatrix}
	D_m&-\sigma I_m\\
	\sigma I_m& D_m
	\end{bmatrix},
	\end{equation*}
	where $D_m=\mathrm{diag}(d_1,...,d_m)$ and $d_j\pm i\sigma$ is an eigenvalue of $A$ for $j=1,...,m$.
	\end{lemma}
	\begin{proof}
	 Let $k\coloneq 2(p-m)$ and $A=\breve{Q}\breve{S}\breve{Q}^\top $ be a RSD with
	 \begin{equation}\label{eq:permuted_Schur}
	 \breve{S} = {\scriptsize \begin{bmatrix}
	 D_m&-\sigma I_m&0&0&0\\
	 \sigma I_m&D_m& 0&0&0\\
	 0&0&\breve{\Lambda}&0&0\\
	 0&0&0&\Lambda_k\cos(\Theta_k)&-\Lambda_k\sin(\Theta_k)\\
	 0&0&0&\Lambda_k\sin(\Theta_k)&\Lambda_k\cos(\Theta_k)\\
	 \end{bmatrix}},
\end{equation}
where $\Lambda_k,\Theta_k$ are $k\times k$ diagonal matrices.\footnote{There is a permutation matrix $P$ such that if $A=QSQ^\top$ is a RSD from \eqref{eq:schur_form}, then $\breve{Q} = QP$ and $\breve{S} = P^\top S P$.} It follows that 
	 \begin{equation}
	 \label{eq:breveQ}
	 \mathrm{skew}(A)V = V\begin{bmatrix}
	0&-\sigma I_m\\
	\sigma I_m&0
	\end{bmatrix}\iff \mathrm{skew}(\breve{S})\breve{Q}^\top V = \breve{Q}^\top V \begin{bmatrix}
	0&-\sigma I_m\\
	\sigma I_m&0
	\end{bmatrix}.
	 \end{equation}
	 Let us divide $\breve{Q}^\top V$ row-wise in three blocks such that \begin{equation*}
	 \breve{Q}^\top V \coloneq \begin{bmatrix}
	 X\\
	 Y\\
	Z
	 \end{bmatrix}\text{ with } X\in\mathbb{R}^{2m\times 2m}, Y\in\mathbb{R}^{r\times 2m}\text{ and }Z \in\mathbb{R}^{2(p-m)\times 2m}.
	 \end{equation*}
	 Equation \eqref{eq:breveQ} implies that $0=Y\left[\begin{smallmatrix}0&-\sigma I_m\\ \sigma I_m &0\end{smallmatrix}\right]$ and therefore $Y=0$. Moreover, letting $Z = \left[\begin{smallmatrix}Z_{11}&Z_{12}\\ Z_{21}&Z_{22}\end{smallmatrix}\right]$ with appropriate blocks, we have
	 \begin{align}
	 \nonumber
	 &\begin{bmatrix}
	 0&-\Lambda_k\sin(\Theta_k)\\
	 \Lambda_k\sin(\Theta_k)&0
	 \end{bmatrix} Z = Z \begin{bmatrix}
	0&-\sigma I_m\\
	\sigma I_m&0
	\end{bmatrix}\\
	\label{eq:satisfaction}
	\iff& \begin{cases}
	\Lambda_k\sin(\Theta_k)Z_{21} = -\sigma Z_{12},\\ \Lambda_k\sin(\Theta_k)Z_{12} = -\sigma Z_{21},
	\end{cases}\quad\text{and,}\quad
	\begin{cases}
	\Lambda_k\sin(\Theta_k)Z_{11} = \sigma Z_{22},\\ \Lambda_k\sin(\Theta_k)Z_{22} = \sigma Z_{11}.
	\end{cases}
	 \end{align}
	 Since all diagonal entries of $\Lambda_k\sin(\Theta_k)$ are nonzero and distinct from $\sigma$, \eqref{eq:satisfaction} is possible if and only if $Z=0$. Finally, since $\breve{Q}\in\mathrm{O}(n)$ and $V\in\mathrm{St}(n,2m)$, we have $\breve{Q}^\top V\in\mathrm{St}(n,2m)$ and thus $X\in\mathrm{O}(2m)$. Letting $R = X^\top $, we have $\breve{Q}_{1:2m} = VR$ and this concludes the proof.
	\end{proof}
\cref{lem:nullspace,lem:eigenspace_identified} show that the matrix $R$ such that $\widehat{Q}R$ is a matrix of Schur vectors acts independently on three sets of columns of $\widehat{Q}$. This fact is the basis for the main structure of the algorithm: the separate computation of the three sets of eigenvalues and eigenvectors, as detailed in the next subsections.
	
	\subsection{The eigenvalues with distinct nonzero imaginary parts} \label{subsec:case1}
	In this subsection, we show in~\cref{thm:distinct_singular_values} that the Schur vectors of $A$ for eigenvalues with distinct nonzero imaginary parts can be obtained by computing the Schur vectors of the skew-symmetric part $\Omega$. Intermediate lemmas are given in \cref{sec:similarities}.
\begin{theorem}\label{thm:distinct_singular_values}
Let $A\in\mathbb{R}^{n\times n}$ be a normal matrix. If $\Lambda_k\sin(\Theta_k)$ is diagonal and contains $k$ distinct non-zero imaginary parts of the non-repeated eigenvalues of $\mathrm{skew}(A)$, then, for all $V\in\mathrm{St}(n,2k)$ such that
	\begin{equation}
	\nonumber
		\mathrm{skew}(A)V =V \begin{bmatrix}
	0&-\Lambda_k\sin(\Theta_k)\\
	 \Lambda_k\sin(\Theta_k)&0
	\end{bmatrix},
	\end{equation}
	 we have
	\begin{equation}
	\nonumber
	 AV = V \begin{bmatrix}
	\Lambda_k\cos(\Theta_k)&-\Lambda_k\sin(\Theta_k)\\
	 \Lambda_k\sin(\Theta_k)&\Lambda_k\cos(\Theta_k)
	\end{bmatrix}\eqcolon V S_{2k}.
	\end{equation}
	If, additionally, $k=\frac{n}{2}$, then $A = VS_{n}V^\top $ is a real Schur decomposition \eqref{eq:schur_form}.
	\end{theorem}
	\begin{proof}
	 Let $A=QSQ^\top $ be a RSD and take a permutation matrix $P$ such that $(P^\top SP)_{1:2k,1:2k}=S_{2k}$. By~\cref{lem:eigenspace_identified}, it follows that $QP = [VR\ \ (QP)_{(k+1):n}]$ for some $R\in \mathrm{O}(2k)$ and this yields 
\begin{align*}
	 \begin{bmatrix}
	0&-\Lambda_k\sin(\Theta_k)\\
	 \Lambda_k\sin(\Theta_k)&0
	\end{bmatrix}&= V^\top \mathrm{skew}(A)V\hspace{3cm}\text{By assumption.}\\
	&= V^\top Q\mathrm{skew}(S)Q^\top  V\\
	&= V^\top (QP)(P^\top \mathrm{skew}(S)P)(QP)^\top  V\\
	&= [R\ 0_{2k\times(n-2k)}](P^\top \mathrm{skew}(S)P)[R\ 0_{2k\times(n-2k)}]^\top \\
	  &= R\begin{bmatrix}
	0&-\Lambda_k\sin(\Theta_k)\\
	 \Lambda_k\sin(\Theta_k)&0
	\end{bmatrix}R^\top. 
\end{align*}	  
Therefore, we have $R\in\mathrm{ig}\left(\left[\begin{smallmatrix}
	0&-\Lambda_k\sin(\Theta_k)\\
	 \Lambda_k\sin(\Theta_k)&0
	\end{smallmatrix}\right]\right)$. Leveraging~\cref{lem:Invariance_lambda}, there is a diagonal matrix $\Phi\in\mathbb{R}^{k\times k}$ such that
	$
	 R=\left[\begin{smallmatrix}
		\cos(\Phi)&-\sin(\Phi)\\
		\sin(\Phi)&\cos(\Phi)
	\end{smallmatrix}\right]$.
	Moreover, it is easily verified that every such matrix $R$ also belongs to $\mathrm{ig}\left(\left[\begin{smallmatrix}
		\Lambda_k\cos(\Theta_k)&0\\
		0&\Lambda_k\cos(\Theta_k)
	\end{smallmatrix}\right]\right)$. This yields \begin{equation*}
	R\in\begin{bmatrix}
	0&-\Lambda_k\sin(\Theta_k)\\
	 \Lambda_k\sin(\Theta_k)&0
	\end{bmatrix}\subset \begin{bmatrix}
	\Lambda_k\cos(\Theta_k)&0\\
		0&\Lambda_k\cos(\Theta_k)
	\end{bmatrix}.
\end{equation*}%\footnote{Note that $\mathrm{ig}\left(\left[\begin{smallmatrix}\Lambda_k\cos(\Theta_k)&0\\0&\Lambda_k\cos(\Theta_k)\end{smallmatrix}\right]\right)$ also admits reflection matrices of the form $\left[\begin{smallmatrix}\cos(\Phi)&\sin(\Phi)\\\sin(\Phi)&-\cos(\Phi)\end{smallmatrix}\right] $, and potentially other families of matrices in case of repeated entries in $\Lambda_k\cos(\Theta_k)$.} 
	Finally, it follows by \cref{lem:oic_linearity} that \begin{equation*}
	R\in \mathrm{ig}\left(\begin{bmatrix}
	\Lambda_k\cos(\Theta_k)&-\Lambda_k\sin(\Theta_k)\\
	 \Lambda_k\sin(\Theta_k)&\Lambda_k\cos(\Theta_k)
	\end{bmatrix}\right)=: \mathrm{ig}(S_{2k}).
	\end{equation*} 
	Consequently, since $QP = [VR\ \ (QP)_{(k+1):n}]$, we have
	\begin{equation*}
		AV = (QP)(P^\top SP) (QP)^\top V =\begin{bmatrix}VR& (QP)_{(k+1):n}\end{bmatrix}\begin{bmatrix}
		S_{2k}&0\\
		0&*
		\end{bmatrix}\begin{bmatrix}
		R^\top \\
		0_{(n-2k)\times 2k}
\end{bmatrix}		 = V RS_{2k}R^\top  = V S_{2k}.
	\end{equation*}
	When $2k=n$, $V\in\mathrm{O}(n)$ and the real Schur decomposition is readily obtained since $S_n=S$ and $A=VSV^\top $. %Hence, $A = \widehat{Q}S\widehat{Q}^\top $ and the proof is complete.
	\end{proof}
	The conclusion of \cref{thm:distinct_singular_values} is that for all eigenvalues that have distinct, non-zero imaginary parts, the invariant subspaces of $\mathrm{skew}(A)$ are always invariant subspaces of~$A$. However, this does not hold true for \wordref{set2}{Set~2} and \wordref{set3}{Set~3} and this implies additional computations, as shown next.
	\subsection{The eigenvalues with repeated nonzero imaginary parts} \label{subsec:case2}
\begin{figure}[ht]
		\centering
		\includegraphics[scale = 1]{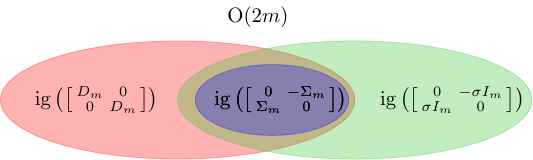}
		\caption{Venn diagram of the sets $\mathrm{ig}\left( \left[\begin{smallmatrix}D_m&0\\
		0&D_m
	\end{smallmatrix}\right]\right)$, $\mathrm{ig}\left(\left[\begin{smallmatrix}
	0&-\sigma I_m\\
	\sigma I_m&0
	\end{smallmatrix}\right]\right)$ and $\mathrm{ig}\left(\left[\begin{smallmatrix}
	0&-\Sigma_m\\
	\Sigma_m&0
	\end{smallmatrix}\right]\right)$ where $D_m$ and $\Sigma_m$ are diagonal with non-zero diagonal entries and $\sigma> 0$.}
	\label{fig:venn_diagram}
	\end{figure}

	Assume there is a diagonal entry $\sigma$ in $\Lambda\sin(\Theta)$ that has multiplicity $m>1$. Then, the associated $2m\times 2m$ block in $S$ can be isolated using a permutation~$P$ such that
	\begin{equation}\label{eq:highlight_multiplicity}
	 (P^\top SP)_{1:2m,1:2m} = 
		\begin{bmatrix}
		D_m&-\sigma I_m\\
		\sigma I_m&D_m
		\end{bmatrix},
	\end{equation}
where $D_m$ is a diagonal submatrix of $\Lambda\cos(\Theta)$. The proof of \cref{thm:distinct_singular_values} relies on the following inclusion of invariant groups: \begin{equation*}
\mathrm{ig}\left(\left[\begin{smallmatrix}
	0&-\Lambda_k\sin(\Theta_k)\\
	 \Lambda_k\sin(\Theta_k)&0
	\end{smallmatrix}\right]\right)\subset \mathrm{ig}\left(\left[\begin{smallmatrix}
		\Lambda_k\cos(\Theta_k)&0\\
		0&\Lambda_k\cos(\Theta_k)
	\end{smallmatrix}\right]\right).
	\end{equation*} Unfortunately,~\cref{lem:invariance_Im} shows that the block $\left[\begin{smallmatrix}
	0&-\sigma I_m\\
	\sigma I_m&0
	\end{smallmatrix}\right]$ admits a larger invariance group than when the eigenvalues have distinct imaginary parts, i.e., if $\Sigma_m$ is an $m\times m$ diagonal matrix with distinct non-zero entries, then $\mathrm{ig}\left(\left[\begin{smallmatrix}
	0&-\Sigma_m\\
	\Sigma_m&0
	\end{smallmatrix}\right]\right)\subset \mathrm{ig}\left(\left[\begin{smallmatrix}
	0&-\sigma I_m\\
	\sigma I_m&0
	\end{smallmatrix}\right]\right) $. However, if $D_m\neq d I_m$, then $\mathrm{ig}\left(\left[\begin{smallmatrix}
	0&-\sigma I_m\\
	\sigma I_m&0
	\end{smallmatrix}\right]\right) \nsubseteq \mathrm{ig}\left( \left[\begin{smallmatrix}D_m&0\\
		0&D_m
	\end{smallmatrix}\right]\right)$. This relation is illustrated in \cref{fig:venn_diagram}.
	
As shown in~\cref{thm:equal_singular_values}, it follows from~\cref{fig:venn_diagram} that if $V\in\mathrm{St}(n,2m)$ spans the invariant subspace of $\mathrm{skew}(A)$ associated to $\sigma$, then $V^\top AV$ is not a real Schur form. Instead, $V^\top AV$  reveals a very structured \emph{symmetric skew-Hamiltonian} matrix.
\begin{theorem}\label{thm:equal_singular_values}
For every normal matrix $A\in\mathbb{R}^{n\times n}$ that has $m> 1$ eigenvalues with imaginary part $\sigma>0$, for all $V\in\mathrm{St}(n, 2m)$ such that 
	\begin{equation}
	\nonumber
		\mathrm{skew}(A)V = V\begin{bmatrix}
	0&-\sigma I_m\\
	\sigma I_m&0
	\end{bmatrix},
	\end{equation}
 we have
	\begin{equation}
	\nonumber
	 AV = V\begin{bmatrix}
	\widetilde{H}&-\widetilde{\Omega}-\sigma I_m\\
	\widetilde{\Omega}+\sigma I_m&\widetilde{H}
	\end{bmatrix},
	\end{equation}
	where $\widetilde{H}\in\mathbb{R}^{m\times m}$ is symmetric and $\widetilde{\Omega}\in\mathbb{R}^{m\times m}$ is skew-symmetric.
	\end{theorem}
	\begin{proof}
	  Let $A=\breve{Q}\breve{S}\breve{Q}^\top $ be a RSD with the structure of \eqref{eq:permuted_Schur}.  By \cref{lem:eigenspace_identified}, there is $R\in\mathrm{O}(2m)$ such that $\breve{Q}_{1:2m}=VR$. Since $\breve{Q}\mathrm{skew}(\breve{S})\breve{Q}^\top V = V\left[\begin{smallmatrix}
	0&-\sigma I_m\\
	\sigma I_m&0
	\end{smallmatrix}\right]$, it follows that \begin{equation}
	R\begin{bmatrix}
	0&-\sigma I_m\\
	\sigma I_m&0
	\end{bmatrix}R^\top =\begin{bmatrix}
	0&-\sigma I_m\\
	\sigma I_m&0
	\end{bmatrix}.
\end{equation} Hence, $R\in \mathrm{ig}\left(\left[\begin{smallmatrix}
	0&-\sigma I_m\\
	\sigma I_m&0
	\end{smallmatrix}\right]\right)$. In~\cref{lem:invariance_Im}, we recall the following equivalence:
	\begin{equation}\label{eq:orthogonal_symplectic}
	R\in \mathrm{ig}\left(\begin{bmatrix}
		0&-\sigma I_m\\
		\sigma I_m&0
		\end{bmatrix}\right)\iff R\in\mathrm{O}(2m) \text{ and } R=\begin{bmatrix}
		E&-F\\
		F&E
		\end{bmatrix},
	\end{equation}
	where $E,F\in\mathbb{R}^{m\times m}$. Recall that such matrices are called \emph{ortho-symplectic}, see, e.g., \cite[Eq.~5.1]{dopico09}. In view of~\eqref{eq:orthogonal_symplectic}, it holds that
	\begin{align}
	\nonumber
	AV&=\mathrm{sym}(A) V+\mathrm{skew}(A)V\\
	&=\breve{Q}\mathrm{sym}(\breve{S})\breve{Q}^\top V+\mathrm{skew}(A)V\\
	\nonumber
	&=VR \begin{bmatrix}
	D_m&0\\
	0&D_m
	\end{bmatrix}R^\top  + V\begin{bmatrix}
	0&-\sigma I_m\\
	\sigma I_m& 0
	\end{bmatrix}\\
	\label{eq:symmetric_subproblem}
	&=V\begin{bmatrix}
	ED_mE^\top  + FD_mF^\top  & (FD_mE^\top  - ED_mF^\top )^\top  - \sigma I_m\\
	 FD_mE^\top  - ED_mF^\top +\sigma I_m&ED_mE^\top   + FD_mF^\top 
	\end{bmatrix}.
	\end{align}
	Defining $\widetilde{H} \coloneq ED_mE^\top  + FD_mF^\top $ and $\widetilde{\Omega}\coloneq FD_mE^\top  - ED_mF^\top $ concludes the proof.
	\end{proof}
	 Let us consider again $Q$, the Schur vectors of $A$, $\widehat{Q}$, the Schur vectors of $\mathrm{skew}(A)$, %, rather than $\mathrm{St}(n,2m)$ from~\cref{thm:equal_singular_values}, 
	 and the permutation matrix $P$ from \eqref{eq:highlight_multiplicity}. In view of~\cref{lem:eigenspace_identified}, $(QP)_{1:2m}  = (\widehat{Q}P)_{1:2m} R$ with $R\in\mathrm{O}(2m)$ and \cref{thm:equal_singular_values} holds with $V=(\widehat {Q}P)_{1:2m}$. 
 Therefore, for each diagonal entry $\sigma$ in $\Lambda\sin(\Theta)$ that has multiplicity $m>1$, it follows from~\cref{thm:equal_singular_values} that we must compute $V^\top AV$ and solve an additional $2m\times 2m$ EVP given by
	\begin{equation}\label{eq:subproblem_of_multiplicity}
			V^\top AV=\begin{bmatrix}
	\widetilde{H}&-\widetilde{\Omega}-\sigma I_m\\
	\widetilde{\Omega}+\sigma I_m&\widetilde{H}
	\end{bmatrix} = R\begin{bmatrix}
	 D_m&-\sigma I_m\\
	 \sigma I_m&D_m
	 \end{bmatrix}R^\top ,
		\end{equation}
		where $R\in\mathrm{O}(2m)$ and $D_m\in\mathbb{R}^{m\times m}$ are to be computed.
		Solving \eqref{eq:subproblem_of_multiplicity} is enough to obtain the correct Schur vectors by $(QP)_{1:2m}  \coloneq V R$.% Different structures of $D_m$ yield different strategies to solve \eqref{eq:subproblem_of_multiplicity}.
	\begin{enumerate}
		\item If $D_m=dI_m$, $d\in\mathbb{R}$, then it is the unique case where $\mathrm{ig}\left(\left[\begin{smallmatrix}
	0&-\sigma I_m\\
	\sigma I_m&0
	\end{smallmatrix}\right]\right)\subset\mathrm{ig}\left( \left[\begin{smallmatrix}d I_m&0\\
		0&d I_m
	\end{smallmatrix}\right]\right)$. Consequently, $\widetilde{H} = d I_m$ and $\widetilde{\Omega} = 0$. No additional work is needed, one can assume $R=I_{2m}$ without loss of generality.
		\item If $D_m$ is not a multiple of the identity, we can recover an ortho-symplectic eigenvector matrix $R$ from the \emph{symmetric skew-Hamiltonian} EVP of $\left[\begin{smallmatrix}\widetilde{H}&-\widetilde{\Omega}\\
	\widetilde{\Omega}&\widetilde{H}\end{smallmatrix}\right]$. Several methods have been designed for Hamiltonian and skew-Hamiltonian EVPs, see, e.g.,  \cite{BENNER199775,WATKINS200423,BennerHamiltonian,VANLOAN1984233,FABENDER1999125}, notably a QR algorithm relying exclusively on ortho-symplectic similarity transformations~\cite{VANLOAN1984233}. For completeness, we briefly discuss in \cref{app:sym_skew_hamiltonian} how the $m$ first steps of Lanczos tridiagonalization, with restarting and reorthogonalization for stability, are enough to construct an ortho-symplectic matrix $Y$ such that%\footnote{For Householder tridiagonalization, we verified numerically that $m-1$ steps were enough, see \href{https://github.com/smataigne/NormalEVP.jl/blob/main/src/wxeigen.jl}{wxeigen.jl} from the repository. Unfortunately, these enhanced Householder and Lanczos methods have been observed to have numerical difficulties so that \wordref{oracle2}{Routine~2} should be preferred, see \cref{fig:WXbenchmark}. The development of a numerically robust version that still benefits from a reduced complexity is an open problem for this particular case of repeated imaginary parts.}
		\begin{equation}\label{eq:skew_hamiltonian}
			\begin{bmatrix}
				\widetilde{H}&-\widetilde{\Omega}\\
	\widetilde{\Omega}&\widetilde{H}
			\end{bmatrix}=Y\begin{bmatrix}
			T_m&0\\
			0&T_m
			\end{bmatrix}Y^\top  = R\begin{bmatrix}
	 D_m&0\\
	 0&D_m
	 \end{bmatrix}R^\top ,
		\end{equation}
		where $T_m$ is an $m\times m$ symmetric tridiagonal matrix. Consequently, once the EVP of $T_m=Z D_m Z^\top $ is solved with a call to \wordref{oracle2}{Routine~2}, we obtain $R\coloneq Y\left[\begin{smallmatrix}Z&0\\ 0 &Z\end{smallmatrix}\right]$.
		
		As studied in \cref{sec:accuracy}, in practice, imaginary parts of the eigenvalues will not be exactly repeated but rather clustered and therefore, the symmetric skew-Hamiltonian structure cannot be guaranteed. This jeopardizes the attractive computational route based on~\eqref{eq:skew_hamiltonian}. In order to maintain the accuracy of the global RSD, it is more robust to call \wordref{oracle3}{Routine~3} to solve \eqref{eq:subproblem_of_multiplicity} directly. This choice has no influence on the computational complexity if $m$ remains small, which is expected. %on If $D_m$ has repeated diagonal entries, only a zero-measure set of the admissible eigenvector matrices of $\left[\begin{smallmatrix}\widetilde{H}&-\widetilde{\Omega}\\ \widetilde{\Omega}&\widetilde{H}\end{smallmatrix}\right]$ exhibit the symplectic structure required to satisfy $R\in\mathrm{ig}\left(\left[\begin{smallmatrix}0&-\sigma I_m\\\sigma I_m&0\end{smallmatrix}\right]\right)$. Working only with the symmetric part seems thus out of reach: one needs to obtain the non-symmetric (normal) $2m\times 2m$ RSD from \eqref{eq:subproblem_of_multiplicity} with a call to \wordref{oracle3}{Routine~3}.		
	\end{enumerate}
	
	\subsection{The real eigenvalues}\label{subsec:r_nonzero} We conclude this section by considering the computation of the $r$ real eigenvalues of $A$ and the associated eigenvectors. In \cref{thm:real_eigenvalues}, we show that an additional $r\times r$ symmetric EVP needs to be solved to obtain the eigenvectors of $A$.
\begin{theorem}\label{thm:real_eigenvalues}
For every normal matrix $A\in\mathbb{R}^{n\times n}$ that has $r>0$ real eigenvalues and for all $V\in\mathrm{St}(n,r)$ such that $\mathrm{skew}(A)V=0$, we have 
	\begin{equation}
	\nonumber
		AV = VH,
	\end{equation}
	where the matrix $H\in\mathbb{R}^{r\times r}$ is symmetric.
	\end{theorem}
	\begin{proof} For every RSD $A=QSQ^\top $ where $S$ follows the structure of~\eqref{eq:schur_form}, we know by \cref{lem:nullspace} that $Q=[Q_{1:p}\ | \ V\breve{R}\ | \ Q_{(p+r+1):n}]$ for some $\breve{R}\in\mathrm{O}(r)$. Therefore, it holds that 
	\begin{align*}
	AV &= QSQ^\top V \\
	&=Q\begin{bmatrix}
	\Lambda\cos(\Theta)&0&-\Lambda\sin(\Theta)\\
	0&\breve{\Lambda}&0\\
	 \Lambda\sin(\Theta)&0&\Lambda\cos(\Theta)
	\end{bmatrix} \begin{bmatrix}
	0_{p\times r}\\
	\breve{R}^\top\\
	0_{p\times r}
	\end{bmatrix}\\
	&=Q\begin{bmatrix}
	0_{p\times r}\\
	\breve{\Lambda}\breve{R}^\top\\
	0_{p\times r}
	\end{bmatrix}\\
	&= V\breve{R} \breve{\Lambda} \breve{R}^\top .
	\end{align*}
	Defining the symmetric matrix $H\coloneq\breve{R}\breve{\Lambda}\breve{R}^\top $ completes the proof.
	\end{proof}	
	\subsection{The algorithm in exact arithmetic} \new{We have gathered all pieces to design an algorithm in exact arithmetic. First, we obtain $\widehat{Q}$ from RSD of $\Omega=\mathrm{skew}(A)$~\eqref{eq:skew_schur1}. By \cref{thm:accuracy}, the invariant spaces associated to \wordref{set1}{Set~1} need no further action. For each invariant space of $\Omega$ associated to \wordref{set2}{Set~2}, i.e, a complex eigenvalue with a repeated imaginary part, an additional symmetric skew-Hamiltonian EVP needs to be solved. This gives the transformation $[Q_{\mathrm{c},1}\ Q_{\mathrm{c},2}] =[\widehat{Q}_{\mathrm{c},1}\ \widehat{Q}_{\mathrm{c},2}]\left[\begin{smallmatrix}R_{\mathrm{a}}&R_{\mathrm{b}}\\R_{\mathrm{c}}&R_{\mathrm{d}}\end{smallmatrix}\right]$ from~\eqref{eq:structure_M}.} Finally, by \cref{thm:accuracy_real}, the decomposition $H\coloneq\widehat{Q}_{\mathrm{r}}^\top A\widehat{Q}_{\mathrm{r}} =\breve{R}\breve{\Lambda}\breve{R}^\top $ using a call to \wordref{oracle2}{Routine~2} completes the determination of the real Schur decomposition 
	\begin{equation*}
	A = \widehat{Q}\begin{bmatrix}
	R_{\mathrm{a}}&0&R_{\mathrm{b}}\\
	0&\breve{R}&0\\
	R_{\mathrm{c}}&0&R_{\mathrm{d}}
\end{bmatrix} \begin{bmatrix}
	\Lambda\cos(\Theta)&0&-\Lambda\sin(\Theta)\\
	0&\breve{\Lambda}&0\\
	 \Lambda\sin(\Theta)&0&\Lambda\cos(\Theta)
	\end{bmatrix}\begin{bmatrix}
	R_{\mathrm{a}}&0&R_{\mathrm{b}}\\
	0&\breve{R}&0\\
	R_{\mathrm{c}}&0&R_{\mathrm{d}}
\end{bmatrix}^\top  \widehat{Q}^\top .
	\end{equation*}

The complete method in exact arithmetic is summarized in \cref{alg:schur_decomposition}. In a broad set of problems, steps~5 and~7 of \cref{alg:schur_decomposition} are either not necessary or very simple. Indeed, for Haar-distributed random orthogonal matrices~\cite{Anderson87, Stewart80}, it holds with probability~$1$ that $r\in\{0,1\}$ if $\det(A)=1$ and $r\in\{1,2\}$ if $\det(A)=-1$~\cite{Fasi21}. In particular, if $r=1$, $H=\pm 1$ and $\breve{R}=1$. Step 7 of \cref{alg:schur_decomposition} is thus trivial. Moreover, all complex eigenvalues have distinct imaginary parts with probability $1$ and thus step~5 is skipped. For general normal matrices, the distribution of the eigenvalues, and in particular the value of $r$, are specific to the context of the problem. Nonetheless, the case where all or a large majority of the eigenvalues belong to \wordref{set1}{Set~1} is common.
	\begin{algorithm}
	\caption{Exact arithmetic algorithm for the real Schur decomposition of a normal matrix}
	\label{alg:schur_decomposition}
	\begin{algorithmic} 
	\STATE \textbf{Input:} A normal matrix $A\in\mathbb{R}^{n\times n}$.
	\STATE \textbf{Output:} $Q$ and $S$ where $Q\in\mathrm{O}(n)$ and $A=QSQ^\top$ is a RSD.
		\STATE \textbf{step 1:} Compute the skew-symmetric part $\Omega\coloneq\frac12 (A-A^\top)$.
		\STATE \textbf{step 2.1:} Compute a tridiagonal reduction $\Omega = \widetilde{Q} T \widetilde{Q}^\top$.
		\STATE \textbf{step 2.2:} Apply the permutation $\left[\begin{smallmatrix}
		0&-\widetilde{B}^\top \\
		\widetilde{B}&0
		\end{smallmatrix}\right]\coloneq P_{\mathrm{eo}}^\top T P_{\mathrm{eo}}$ where $\widetilde{B}\in\mathbb{R}^{\lfloor\frac{n}{2}\rfloor\times \lceil\frac{n}{2}\rceil}$ is bidiagonal.
		\STATE \textbf{step 3:} Compute the (permuted) SVD of  $\widetilde{B}= U\left[\begin{smallmatrix} 0&0_{\lfloor\frac{r}{2}\rfloor\times \lceil\frac{r}{2}\rceil}\\ \Sigma  &0 \end{smallmatrix}\right]V^\top $ where $\mathrm{diag}(\Sigma)>0$.
		\STATE \textbf{step 4:} Define $\widehat{Q}\coloneq \widetilde{Q} P_{\mathrm{eo}} \left[\begin{smallmatrix}
		V&0\\
		0&U\end{smallmatrix}\right]$, $\Lambda \sin(\Theta) \coloneq \Sigma$ and $p \coloneq \frac{n-r}{2}$.
		\FOR{each non-zero diagonal entry $\sigma$ of $\Sigma$ with multiplicity $m > 1$}
			\STATE \textbf{step 5.1:} Define $V$ as the $2m$ columns of $\widehat{Q}$ such that $\Omega V = V\left[\begin{smallmatrix}
		0&-\sigma I_m\\
		\sigma I_m&0
		\end{smallmatrix}\right]$.
			\STATE \textbf{step 5.2:} Compute the matrix $ V^\top AV\eqcolon\left[\begin{smallmatrix}
		\widetilde{H}&-\widetilde{\Omega}-\sigma I_m\\
	\widetilde{\Omega}+\sigma I_m&\widetilde{H}
		\end{smallmatrix}\right]$ where $\widetilde{H}=\widetilde{H}^\top$ and $\widetilde{\Omega} = -\widetilde{\Omega}^\top$. 
			\STATE \textbf{step 5.3:} Compute the RSD $ \left[\begin{smallmatrix}
		\widetilde{H}&-\widetilde{\Omega}-\sigma I_m\\
	\widetilde{\Omega}+\sigma I_m&\widetilde{H}
		\end{smallmatrix}\right]=R\left[\begin{smallmatrix}
		D_m&-\sigma I_m\\
		\sigma I_m&D_m
		\end{smallmatrix}\right] R^\top $.
			\STATE \textbf{step 5.4:} Replace the columns of $\widehat{Q}$ selected at step 5.1 by $VR$.
		\ENDFOR
		\STATE \textbf{step 6:} Define the diagonal matrix $\Lambda\cos(\Theta)\coloneq\widehat{Q}_{1:p}^\top  A \widehat{Q}_{1:p}$.
		\IF{$r>0$}
		\STATE \textbf{step 7.1:} Define $V \coloneq \widehat{Q}_{(p+1):(p+r)}$  such that $\Omega V=0$.
		\STATE \textbf{step 7.2:} Compute the symmetric matrix $H\coloneq V^\top A V$. 
		\STATE \textbf{step 7.3:} Compute the EVD of $H=\breve{R}\breve{\Lambda} \breve{R}^\top $.
		\STATE \textbf{step 7.4:} Replace $\widehat{Q}_{(p+1):(p+r)}$ by $V\breve{R}$.
		\ENDIF
		\RETURN $Q\coloneq\widehat{Q}$ and $S\coloneq\left[\begin{smallmatrix}
		\Lambda \cos(\Theta)&0&-\Lambda\sin(\Theta)\\
		0&\breve{\Lambda}&0\\
		\Lambda\sin(\Theta)&0 &\Lambda \cos(\Theta)
		\end{smallmatrix}\right]$.
	\end{algorithmic}%
	\end{algorithm}
	\new{
	\subsection{The eigenvalue decomposition by the symmetric part} \cref{alg:schur_decomposition} shows that the real Schur decomposition of real normal matrix can be obtained from its skew-symmetric part. It is natural to ask if this approach can be adapted to start from the \emph{symmetric part}. The answer is positive but several drawbacks arise compared to the skew-symmetric approach.
	\begin{itemize}
		\item Drawback~1: The symmetric part does not allow separating real from complex eigenvalues easily. Indeed, given an EVD of $\mathrm{sym}(A)= V\Lambda V^\top $, one can not distinguish a real eigenvalue of $A$ from the real part of the complex eigenvalue, except if its multiplicity is $1$. Making this separation requires computing the product $V^\top A V$. In \cref{alg:schur_decomposition}, there is a direct correspondence between the real eigenspace of $A$ and the kernel of $\mathrm{skew}(A)$.
		\item Drawback~2: The aforementioned separation between real and complex eigenvalues can not be obtained if some real parts of complex eigenvalues are equal to real eigenvalues of $A$. Then, the associated sub-block of  $V^\top A V$ is a normal matrix with no additional structure. 
		\item Drawback~3: The symmetric eigensolver for $\mathrm{sym}(A)$ does not take advantage of the even multiplicity of the real parts of complex eigenvalues, while the skew-symmetric part uses this advantage by computing a half-the-size bidiagonal SVD.
		\item Drawback~4: For every two-dimensional invariant subspace of $A$, the problem arises that
		\begin{equation*}
			\mathrm{ig}\left(\begin{bmatrix}
			\lambda \cos(\theta)&-\lambda\sin(\theta)\\
			\lambda \sin(\theta)&\lambda\cos(\theta)\\
			\end{bmatrix}\right) =\mathrm{SO}(2)\subset \mathrm{ig}\left(\begin{bmatrix}
			\lambda \cos(\theta)&0\\
			0&\lambda\cos(\theta)\\
			\end{bmatrix}\right) =\mathrm{O}(2).
		\end{equation*}
		Since an element of $\mathrm{O}(2)$ belongs to $\mathrm{SO}(2)$ with probability $\frac12$, there is a probability $\frac12$ for each base to be incorrectly oriented.
	\end{itemize}
	The aforementioned drawbacks matter if a significant proportion of the eigenvalues are complex. If most eigenvalues are real, then, the approach by the symmetric part can be efficient. However, applications, e.g., with orthogonal matrices, are usually dominated by complex eigenvalues. This makes the approach by the skew-symmetric-part more relevant for applications to the best of the authors' knowledge.}
	 
	 \subsection{Exact versus floating point arithmetic: an example}\label{sec:sensitivity} We have not yet adressed an essential algorithmic question: how accurate is \cref{alg:schur_decomposition} in floating point arithmetic and how should it be implemented in order to always obtain a desired accuracy? Indeed, floating point arithmetic brings an additional challenge. In some cases, the invariant subspaces of $\Omega$ are less sensitive to perturbations of $A$ than the invariant subspaces of $A$ themselves, and, if not analyzed and addressed in the design of the algorithm, this can be a source of inaccuracies. We show next an illustrative example.
	 
	Consider a variable $\tau\geq 0$, a $p\times p$ random symmetric matrix $E$, $\|E\|_2=1$, and the normal matrix $A(\tau)$ defined by
	\begin{equation*}
	A(\tau) = \begin{bmatrix}
	\sin(\tau E)&-\cos(\tau E)\\
	\cos(\tau E)&\sin(\tau E)
	\end{bmatrix}\quad\text{such that,}\quad  A(0)=\begin{bmatrix}0&-I_p\\
	I_p&0 \end{bmatrix}.
	\end{equation*}
Notice that $A(0)=\Omega(0)$. For $\tau$ small enough, Taylor series yield $\cos(\tau E)\approx I_p -\frac{\tau^2}{2}E^2$ and $\sin(\tau E)\approx \tau E$ such that, using the spectral norm $\|\cdot\|_2$, we have
\begin{equation*}
	 \|A(\tau)-A(0)\|_2\approx \tau \quad\text{and,}\quad  \|\Omega(\tau)-\Omega(0)\|_2\approx \frac{\tau^2}{2}.
\end{equation*} 
Around $\tau=0$, it follows that $\Omega(\tau)$ is less sensitive to small variations of $\tau$ than $A(\tau)$. Let $\varepsilon_\mathrm{m}$ denote the machine precision. Then, applying \cref{alg:schur_decomposition} on $A(\sqrt{\varepsilon_\mathrm{m}})$ or $A(0)$ will yield similar Schur vectors in output since, in floating point arithmetic, $\Omega(\sqrt{\varepsilon_\mathrm{m}})\approx \Omega(0)$. Therefore, $\mathcal{O}(\sqrt{\varepsilon_\mathrm{m}})$-errors will be made in the computation of the Schur vectors of $A(\sqrt{\varepsilon_\mathrm{m}})$. This is verified in practice, here for $p=100$ and $\varepsilon_\mathrm{m}=2.22\cdot 10^{-16}$, we have
\begin{equation*}
	\text{Naive implementation of \cref{alg:schur_decomposition}:}\quad \|A(\sqrt{\varepsilon_\mathrm{m}}) \widehat{Q} - \widehat{Q} \widehat{S}\|_\mathrm{F} = 6.52\cdot 10^{-9} \|A(\sqrt{\varepsilon_\mathrm{m}})\|_\mathrm{F}.
\end{equation*}
\new{For methods relying on the diagonalization of commuting matrices, it is well known that sensitivity analyses depend on eigenvalue gaps. In~\cite{HeKressner2024}, which studies the simultaneous diagonalization of symmetric matrices, a probabilistic accuracy bound is derived in terms of the smallest eigenvalue gap among the family of matrices. Similarly, the do-one-then-do-the-other algorithm of~\cite[Alg.~2]{Sutton23} requires a refinement of the invariant subspaces that is governed by the eigenvalue gaps. As shown in \cref{sec:accuracy}, the specific setting of real normal matrices, together with \cref{alg:schur_decomposition_floating}, allows for more precise statements, again involving eigenvalue gaps.} 

In the next section, we study how the loss of accuracy occurs and how  \cref{alg:schur_decomposition} should be implemented in order to obtain a desired accuracy $\mu$.
	
\section{Analysis of numerical stability and accuracy}\label{sec:accuracy} In this section, we analyze the effect of floating point arithmetic on the accuracy of \cref{alg:schur_decomposition}. We rely on standard assumptions of floating point arithmetic. In particular, $\varepsilon_\mathrm{m}$ denotes the machine precision and, for every scalar $\alpha\in\mathbb{C}$, $\fl(\alpha)=\alpha(1+\Delta \alpha)$ with $|\Delta\alpha|\leq\varepsilon_\mathrm{m}$ models its float representation. The \emph{absolute gap} between two numbers $\alpha,\widehat{\alpha}\in\mathbb{C}$ is defined as $|\alpha - \widehat{\alpha}|$ while the \emph{relative gap} refers to $\frac{|\alpha - \widehat{\alpha}|}{|\alpha| + |\widehat{\alpha}|}$. The notation $\alpha\lesssim \beta$ means that there is a moderately growing function $\varphi(n)$, e.g., $n$ or $n^2$, independent of the values of $\alpha$ and $\beta$, such that $\alpha\leq \varphi(n)\beta$. 

We will analyze step by step the numerical properties of the algorithm. This will determine how to implement \cref{alg:schur_decomposition} to obtain high accuracy on both the eigenvalues and the Schur vectors/invariant spaces. The floating point arithmetic version of the method is  given in  \cref{alg:schur_decomposition_floating}. \new{The algorithm requires to identify clusters of eigenvalues. A definition for a cluster of eigenvalues, justified by \cref{thm:accuracy,thm:accuracy_repeated,thm:accuracy_real} is given~\cref{def:cluster}.
\begin{definition}\label{def:cluster}
Given $\delta>0$, a matrix $A$, a set of $2m>0$ eigenvalues $\{\pm i\sigma_j\}_{j=1}^{m}$ of $\Omega =\mathrm{skew}(A)$ with $\sigma_1\geq \sigma_2\geq...\geq \sigma_m\geq0$ 
is a \emph{$\delta$-cluster} of $\Omega$ if and only if
	\begin{equation*}
		|\sigma_j - \sigma_{j+1}|\leq \delta \|A\|_\mathrm{F}\quad \text{ and}\quad |\sigma_k - \sigma|> \delta \|A\|_\mathrm{F}\quad \text{for } j=1,...,m-1, \text{ and, } k=1,...,m,
	\end{equation*}
	where $i\sigma$ is any eigenvalue of $\Omega$ not in $\{\pm i\sigma_j\}_{j=1}^m$.
\end{definition}
\cref{def:cluster} is similar to the clustering criterion used in~\cite[Alg.~2]{Sutton23}. Instead of the Frobenius norm~$\|A\|_\mathrm{F}$, the results of \cite{Sutton23} employ the \emph{spread} of the eigenvalues, i.e., the maximum distance between the eigenvalues. This difference in definitions reflects the distinct hypotheses and analytical frameworks underlying the respective results.}
\begin{algorithm}
	\caption{Floating point arithmetic algorithm for the real Schur decomposition of a normal matrix (\texttt{nrmschur})}
	\label{alg:schur_decomposition_floating}
	\begin{algorithmic}
	\STATE \textbf{Input:} A normal matrix $A\in\mathbb{R}^{n\times n}$, $\delta, \delta_r=\sqrt{\varepsilon_\mathrm{m}}$ and $\mu = \varepsilon_\mathrm{m} t$, $t\geq 1$.
	\STATE \textbf{Output:} $\widehat{Q}$ and $\widehat{S}$ where $d_{\mathrm{O}(n)}(\widehat{Q})\lesssim\varepsilon_\mathrm{m}$ and $\|A\widehat{Q} - \widehat{Q} \widehat{S}\|_\mathrm{F}\lesssim \mu \|A\|_\mathrm{F}$.
		\STATE \textbf{step 1:} Compute the skew-symmetric part $\widehat{\Omega}\approx \frac12(A-A^\top)$.
		\STATE \textbf{step 2.1:} Tridiagonalize $\widehat{\Omega} \approx \widetilde{Q} T \widetilde{Q}^\top$.
		\STATE \textbf{step 2.2:} Apply the permutation $\left[\begin{smallmatrix}
		0&-\widetilde{B}^\top \\
		\widetilde{B}&0
		\end{smallmatrix}\right]\coloneq P_{\mathrm{eo}}^\top T P_{\mathrm{eo}}$ where $\widetilde{B}\in\mathbb{R}^{\lfloor\frac{n}{2}\rfloor\times \lceil\frac{n}{2}\rceil}$ is bidiagonal.
		\STATE \textbf{step 3:}  (\wordref{oracle1}{Routine~1}) Compute the (permuted) SVD of  $\widetilde{B} \approx \widehat{U}\left[\begin{smallmatrix} 0&0_{\lfloor\frac{\widehat{r}}{2}\rfloor\times \lceil\frac{\widehat{r}}{2}\rceil}\\ \widehat{\Sigma}  &0 \end{smallmatrix}\right]\widehat{V}^\top $ where $\mathrm{diag}(\widehat{\Sigma})>0$ and $\widehat{r}$ is the size of the $\delta_r$-cluster of eigenvalues of $\widehat{\Omega}$ around 0.
		\STATE \textbf{step 4:} Compute $\widehat{Q}\approx \widetilde{Q} P_{\mathrm{eo}} \left[\begin{smallmatrix}
		V&0\\
		0&U\end{smallmatrix}\right]$, set $\Lambda \sin(\Theta) \coloneq \widehat{\Sigma}$ and $\widehat{p}\coloneq \frac{n-\widehat{r}}{2}$.
		\FOR{each $\delta$-cluster of $2m$ (non-zero) eigenvalues of $\widehat{\Omega}$}
			\STATE \textbf{step 5.1:} Set $V$ as the $2m$ columns of $\widehat{Q}$ associated to the cluster and compute $M\approx V^\top AV$. 
			\STATE \textbf{step 5.2:} (\wordref{oracle3}{Routine~3}) Compute the RSD $M \approx R \widetilde{S}_{M} R^\top $.
			\STATE \textbf{step 5.3:} Replace the $2m$ columns of $\widehat{Q}$ from step 5.1 by $VR$.
		\ENDFOR
		\STATE \textbf{step 6:} Set $\Lambda\cos(\Theta)$ as the diagonal of $\widehat{Q}_{1:\widehat{p}}^\top  \mathrm{sym}(A) \widehat{Q}_{1:\widehat{p}}$.
		\IF{there is a $\delta_r$-cluster of $\widehat{r}$ eigenvalues of $\Omega$ around $0$,}
		\STATE \textbf{step 7.1:} Compute the matrix $H\approx \widehat{Q}_{\widehat{p}+1:\widehat{p}+\widehat{r}}^\top A\widehat{Q}_{\widehat{p}+1:\widehat{p}+\widehat{r}}$. 
		\IF{$\|\mathrm{skew}(H)\|_\mathrm{F}\leq \varepsilon_\mathrm{m} \|H\|_\mathrm{F}$,}
		\STATE \textbf{step 7.2a:} (\wordref{oracle2}{Routine~2}) Compute the EVD of $\mathrm{sym}(H)\approx\breve{R}\breve{\Lambda} \breve{R}^\top $.
		\ELSE
			\STATE \textbf{step 7.2b:} (\wordref{oracle3}{Routine~3}) Compute the RSD of $H\approx \breve{R}\breve{S}_H \breve{R}^\top $.
		\ENDIF
		\STATE \textbf{step 7.3:} Set $\widehat{Q}_{\widehat{p}+1:\widehat{p}+\widehat{r}}\approx \widehat{Q}_{\widehat{p}+1:\widehat{p}+\widehat{r}} \breve{R}$.
		\ENDIF
		\STATE \textbf{step~8:} (Optional, see \cref{subsec:accuracy_cont}) Subspaces correction steps. For each $\frac{1}{t}$-cluster of eigenvalues of $\widehat{\Omega}$, repeat step 5 with $\delta=\frac{1}{t}$ and consider a Jacobi algorithm~\cite{ZhouBrent03} for \wordref{oracle3}{Routine~3}.
		\RETURN $\widehat{Q}$ and $\widehat{S}\coloneq \left[\begin{smallmatrix}
		\Lambda \cos(\Theta)&0&-\Lambda\sin(\Theta)\\
		0&\breve{\Lambda}&0\\
	\Lambda\sin(\Theta)&0 &\Lambda \cos(\Theta)
		\end{smallmatrix}\right]$.
	\end{algorithmic}%
	\end{algorithm}
\subsection{Analysis of steps 1 to 4}
	\textbf{Step 1} already entails some subtleties since it consists of $n^2$ subtractions. Some rare pathological cases exist where this operation might be an initial source of inaccuracy. Indeed, letting $\Omega=\mathrm{skew}(A)$ and $\widehat{\Omega}=\mathrm{skew}(\fl(A))$ such that $\widehat{\omega}_{ij}=\frac{\fl(a_{ij})-\fl(a_{ji})}{2}$, we obtain
	\begin{equation}\label{eq:skew_computation}
		\left|\frac{\omega_{ij}-\widehat{\omega}_{ij}}{\omega_{ij}}\right|=\left|\frac{a_{ij}\Delta a_{ij}-a_{ji}\Delta a_{ji}}{a_{ij}-a_{ji}} \right|\leq \varepsilon_\mathrm{m}\frac{|a_{ij}|+|a_{ji}|}{|a_{ij}-a_{ji}|}.
	\end{equation}
Based on the elementwise relation \eqref{eq:skew_computation} and acknowledging that $\mathrm{Tr}(|A|^2)\leq\|A\|_\mathrm{F}^2$ by Cauchy-Schwarz inequality, we can also deduce the more practical matrix-relation
%\begin{align*}
%\|\Omega-fl(\Omega)\|_\mathrm{F}^2 &= \sum_{i,j}(\omega_{ij}-fl(\omega_{ij}))^2\\&=\sum_{i,j}(a_{ij}\Delta a_{ij}-a_{ji}\Delta a_{ji})^2\\&\leq \varepsilon_\mathrm{m}^2 \sum_{i,j} a_{ij}^2+ a_{ji}^2+2|a_{ij}| |a_{ji}|\\&= \varepsilon_\mathrm{m}^2 (2\|A\|_\mathrm{F}^2 + 2\mathrm{Tr}(|A|^2))\\&\leq 4\varepsilon_\mathrm{m} \|A\|_\mathrm{F}^2 
%\end{align*}
\begin{equation}\label{eq:matrix_accuracy}
 \frac{\|\Omega-\widehat{\Omega}\|_\mathrm{F}}{\|\Omega\|_\mathrm{F}}\leq \varepsilon_\mathrm{m} \frac{2\|A\|_\mathrm{F}}{\|\Omega\|_\mathrm{F}}.
 \end{equation} In view of \eqref{eq:skew_computation}, whenever symmetrically positioned entries of $A$ are \emph{relatively} close, the numerical approximation $\widehat{\omega}_{ij}$ may become inaccurate. % Since only the relative gap intervenes in \eqref{eq:skew_computation}, symmetrically positioned zeros of $A$ are harmless. 
However, relative inaccuracies only occur for $\omega_{ij}$ close to zero. Therefore, if $\|\Omega\|_\mathrm{F}\gg \varepsilon_\mathrm{m}\|A\|_\mathrm{F}$, then $\widehat{\Omega}$ is accurate in the sense of \eqref{eq:matrix_accuracy}, even if some entries are not. An example of a matrix set where such entrywise inaccuracies occur in general with negligible effect is 
	\begin{equation}\label{eq:matrix_set}
	\{U\in\mathrm{O}(2m)\ | \ U = \exp\left[\begin{smallmatrix}
		0&-M^\top \\
		M&0
		\end{smallmatrix}\right],\ M\in\mathbb{R}^{m\times m}\}.
	\end{equation}
	The matrix set \eqref{eq:matrix_set} notably arises in the fundamental Riemannian computation of geodesics on the Grassmann manifold \cite{bendokat2024}. Entrywise inaccuracies can also be neglected if $\|\Omega\|_\mathrm{F}\leq \varepsilon_\mathrm{m} \|A\|_\mathrm{F}$. In this case, the matrix $A$ can be numerically considered as symmetric. In conclusion, computing a skew-symmetric part is dubious when $\varepsilon_\mathrm{m}\|A\|_\mathrm{F}\ll \|\Omega\|_\mathrm{F}\ll\|A\|_\mathrm{F}$. However, it is shown in \cref{sec:best_worst_case} that \cref{alg:schur_decomposition_floating} can handle small skew-symmetric parts.
	
	\textbf{Step 2} is %if $\mathrm{dim}(\ker(\Omega))=r$,\footnote{Due to perturbations, $r$ can not (and need not) to be determined with certitude. We only use it for notation purposes. However, by the Hoffman-Wielandt inequality \cite[Thm.~1]{HoffmanWielandt}, zero eigenvalues are absolutely perturbed by $\|\Delta fl(\Omega)\|_\mathrm{F}$ at most.} 
the tridiagonal reduction using Householder reflectors. It is known to be backward stable \cite[p.~416]{GoluVanl96}, i.e.,%\footnote{\eqref{eq:backward_stability} is in fact a relaxation of the best stability result that Householder transformations can provide \cite[p.~360]{Higham2002}.}
	\begin{equation}\label{eq:backward_stability}
		\widehat{\Omega} + \Delta\widehat{\Omega} = \widetilde{Q}P_{\mathrm{eo}}\left[\begin{smallmatrix}
		0&-\widetilde{B}^\top \\
		\widetilde{B}&0
		\end{smallmatrix}\right] P_{\mathrm{eo}}^\top \widetilde{Q}^\top  \text{ with } \|\Delta\widehat{\Omega}\|_\mathrm{F}\lesssim \varepsilon_\mathrm{m}\|\widehat{\Omega}\|_\mathrm{F},
	\end{equation}
	where $\Delta\widehat{\Omega}$ is exactly skew-symmetric, $\widetilde{Q}\in\mathrm{O}(n)$ is the exact product of the Householder reflectors numerically computed during the reduction and $P_{\mathrm{eo}}$ is an even-odd permutation \eqref{eq:even_odd}. Moreover, Householder reflectors can be assembled very accurately as $\widetilde{Q}_{\mathrm{asm}}$ such that $\|\widetilde{Q}_{\mathrm{asm}}-\widetilde{Q}\|_\mathrm{F} \lesssim \varepsilon_\mathrm{m}$ \cite[Eq.~19.13]{Higham2002}. Defining the distance to orthogonality of a matrix $V$ by 
	\begin{equation}\label{eq:dist_to_orthogonality}
	d_{\mathrm{O}(n)}(V)\coloneq\min_{Q\in\mathrm{O}(n)} \|Q-V\|_\mathrm{F},
	\end{equation} 
we can write $d_{\mathrm{O}(n)}(\widetilde{Q}_\mathrm{asm})\lesssim\varepsilon_\mathrm{m}$.%This result does not hold true for symmetric tridiagonal matrices \cite[Chap.~4.2]{dhillon1997new}, which is why symmetric eigensolver such as \href{https://www.netlib.org/lapack/explore-html/d1/d56/group__heevr_gaa334ac0c11113576db0fc37b7565e8b5.html#gaa334ac0c11113576db0fc37b7565e8b5}{syevr} rely on \emph{relatively robust representations} \cite[Def.~5.2.1]{dhillon1997new}. 

\textbf{Step~3} is the bidiagonal SVD. Algorithms such as \href{https://www.netlib.org/lapack//explore-html/d6/d51/group__bdsqr_gade20fbf9c91aa7de0c3d565b39588dc5.html}{bdsqr} provide a mixed forward-backward stable SVD of $\widetilde{B}\in\mathbb{R}^{\lfloor\frac{n}{2}\rfloor\times \lceil\frac{n}{2}\rceil}$:
	\begin{equation*}
	\widetilde{B} + \Delta\widetilde{B} = \widehat{U}\widehat{\Sigma} \widehat{V}^\top\quad \text{and,}\quad \|\Delta\widetilde{B}\|_\mathrm{F}\lesssim\varepsilon_\mathrm{m} \|\widetilde{B}\|_\mathrm{F},
	\end{equation*}
	such that $\max\{d_{\mathrm{O}(\lfloor\frac{n}{2}\rfloor)}(\widehat{U}),d_{\mathrm{O}(\lceil\frac{n}{2}\rceil)}(\widehat{V})\}\lesssim \varepsilon_\mathrm{m}$ \cite{lapack99, DemmelKahan90}. 
	\begin{remark}
	\emph{The singular values of $\widetilde{B}$ are computed with high \emph{relative} accuracy by algorithms such as \href{https://www.netlib.org/lapack//explore-html/d6/d51/group__bdsqr_gade20fbf9c91aa7de0c3d565b39588dc5.html}{bdsqr} \cite[Thm.~6~and~7]{DemmelKahan90}. This holds because relatively small perturbations of $\widetilde{B}$ yield relatively perturbed singular values \cite[Chap.~4.3]{dhillon1997new}. Moreover, it was conjectured~\cite{DemmelKahan90} and proved \cite[Thm.~6.1]{Deift91} that left and right singular vectors of $\widetilde{B}$ are also very accurately computed. %can be computed with an accuracy inversely proportional to the relative gap between the singular values. %: $\mathrm{relgap}(\sigma_i)\coloneq\min_{j=1,...,\lfloor\frac{n}{2}\rfloor}\frac{|\sigma_i-\sigma_j|}{\sigma_i+\sigma_j}$. If $\angle(\cdot,\cdot)$ denotes the angle between two vectors, then, for $i=1,...,\lfloor\frac{n}{2}\rfloor$, we have\footnote{We recall the following relation for $u,\widehat{u}\in\mathrm{St}(n,1)$: $\frac{\|u-\widehat{u}\|_\mathrm{F}}{2} = \sin\left(\frac{\angle(u,\widehat{u})}{2}\right)$.}\begin{equation*}\max(\sin\angle(u_i,\widehat{u}_i),\sin\angle(v_i,\widehat{v}_i))\lesssim \varepsilon_\mathrm{m}\mathrm{relgap}(\sigma_i)^{-1}.\end{equation*}
However, the Bauer-Fike theorem \cite[Thm.~6.3.2]{horn13} and its corollary \cite[Cor.~6.3.4]{horn13} indicate that, due to the perturbation $\Delta\widehat{\Omega}$ in \eqref{eq:backward_stability}, relative accuracy of the eigenvalues of $\Omega$ cannot be achieved. Consequently, the relative accuracy of the singular values of $\widetilde{B}$ is not relevant.}
	\end{remark}

	%Such values are thus very well relatively separated. A fortunate consequence is that the cluster of numerically zero eigenvalues of $\Omega$, i.e., with modulus less than $\varepsilon_\mathrm{m}$, can be relatively isolated from the other eigenvalues in any practical situation. Therefore, the numerical approximation of the vector space $\ker(\Omega)$ is computed to high accuracy.
	
	\textbf{Step 4} is a product of matrices. Define $\widehat{Q}_{\mathrm{exact}}$ as the \emph{exact} matrix product of $\widetilde{Q}_{\mathrm{asm}}P_{\mathrm{eo}}$ by $\left[\begin{smallmatrix}
		\widehat{V}&0\\
		0&\widehat{U}
		\end{smallmatrix}\right]$. Since the permutation is performed exactly, w.l.o.g., we assume $P_{\mathrm{eo}}=I_n$. We obtain that 
		\begin{align*}
		d_{\mathrm{O}(n)}(\widehat{Q}_{\mathrm{exact}})&\leq \|\widetilde{Q}_{\mathrm{asm}}\left[\begin{smallmatrix}
		\widehat{V}&0\\
		0&\widehat{U}
		\end{smallmatrix}\right] - \widetilde{Q}\left[\begin{smallmatrix}
		V&0\\
		0&U
		\end{smallmatrix}\right]\|_\mathrm{F}\\
		&= \|\widetilde{Q}_{\mathrm{asm}}\left[\begin{smallmatrix}
		\widehat{V}-V+V&0\\
		0&\widehat{U}-U+U
		\end{smallmatrix}\right] - \widetilde{Q}\left[\begin{smallmatrix}
		V&0\\
		0&U
		\end{smallmatrix}\right]\|_\mathrm{F}\\
		&\leq  \|\widetilde{Q}_{\mathrm{asm}} - \widetilde{Q}\|_\mathrm{F} + \|\widetilde{Q}_{\mathrm{asm}} \left[\begin{smallmatrix}
		\widehat{V}-V&0\\
		0&\widehat{U}-U
		\end{smallmatrix}\right]\|_\mathrm{F}\\
		&\lesssim \varepsilon_\mathrm{m}.
		\end{align*}
	Finally, it holds by \cite[Eq.~3.13]{Higham2002} that $\|\widehat{Q} - \widehat{Q}_{\mathrm{exact}}\|_\mathrm{F}\lesssim \varepsilon_\mathrm{m} \|\widetilde{Q}_{\mathrm{asm}}\|_\mathrm{F}\|\left[\begin{smallmatrix}
		V&0\\
		0&U
		\end{smallmatrix}\right]\|_\mathrm{F}\lesssim \varepsilon_\mathrm{m}$. Therefore, by the triangular inequality, we have $d_{\mathrm{O}(n)}(\widehat{Q}) \lesssim \varepsilon_\mathrm{m}$. By summing all errors from steps 1 to 4, we obtain the skew-symmetric error matrix $\Delta_2\widehat{\Omega}$:
		\begin{align*}
		\Delta_2\widehat{\Omega} &\coloneq \Delta\widehat{\Omega} + \widetilde{Q}\left[\begin{smallmatrix}
		0&-\Delta\widetilde{B}^\top \\
		\Delta\widetilde{B}&0
		\end{smallmatrix}\right]\widetilde{Q}^\top \\
		&+ (\widetilde{Q}_\mathrm{asm} -\widetilde{Q}) \left[\begin{smallmatrix}0&-\widetilde{B}^\top \\ \widetilde{B}&0 \end{smallmatrix}\right]\widetilde{Q}^\top - \widetilde{Q}  \left[\begin{smallmatrix} 0&-\widetilde{B}^\top \\ \widetilde{B}&0 \end{smallmatrix}\right]^\top(\widetilde{Q}_\mathrm{asm} - \widetilde{Q})^\top\\
		&+(\widehat{Q} -\widehat{Q}_{\mathrm{exact}}) \left[\begin{smallmatrix}
		0&-\widehat{\Sigma} \\
		\widehat{\Sigma}&0
		\end{smallmatrix}\right]\widehat{Q}^\top-\widehat{Q}\left[\begin{smallmatrix}
		0&-\widehat{\Sigma} \\
		\widehat{\Sigma}&0
		\end{smallmatrix}\right]^\top(\widehat{Q}-\widehat{Q}_{\mathrm{exact}})^\top +\mathcal{O}(\varepsilon_\mathrm{m}^2). 
		\end{align*}
$\mathcal{O}(\varepsilon_\mathrm{m}^2)$ gathers terms where, e.g., $\widetilde{Q}_\mathrm{asm} -\widetilde{Q}$ appears twice. If we consider $(\widehat{Q},\widehat{\Sigma})$ as the output of the skew-symmetric RSD, this yields the mixed forward-backward stable decomposition of the skew-symmetric part~$\widehat{\Omega}$:
	\begin{equation}\label{eq:stable_skew_EVP}
	\widehat{\Omega} + \Delta_2\widehat{\Omega} = \widehat{Q}\begin{bmatrix}
	0&-\widehat{\Sigma}\\
		\widehat{\Sigma}&0
	\end{bmatrix}\widehat{Q}^\top  \quad\text{with}\quad \|\Delta_2\widehat{\Omega}\|_\mathrm{F}\lesssim \varepsilon_\mathrm{m}\|\widehat{\Omega}\|_\mathrm{F}\quad\text{and}\quad d_{\mathrm{O}(n)}(\widehat{Q})\lesssim \varepsilon_\mathrm{m}.
	\end{equation}
\emph{Mixed forward-backward} stability refers to the fact that $\widehat{Q}$ is near orthogonal. In this context, the decomposition would be said \emph{backward} stable if $\widehat{Q}$ was exactly orthogonal. As shown in~\cref{sec:exp_accuracy}, $d_{\mathrm{O}(n)}(\widehat{Q})$ is very close to machine precision in practice.

\subsection{Sensitivity of the invariant subspaces: Set~1}
Let $A$ be a normal matrix with an exact RSD $A=QSQ^\top$. We first assume that the matrix $A$ has no real eigenvalues ($r=0, n=2p$) and that all eigenvalues have distinct imaginary parts. Given a backward stable RSD of $\widehat{\Omega}$, i.e., $\widehat{\Omega} \widehat{Q} -\widehat{Q}\mathrm{skew}(\widehat{S})=\widehat{E}$ such that $ \|\widehat{E}\|_\mathrm{F}\lesssim \varepsilon_\mathrm{m} \|\widehat{\Omega}\|_\mathrm{F}$, we quantify $\|A \widehat{Q} -\widehat{Q}S\|_\mathrm{F}$ in \cref{thm:accuracy}. We first assume that $\widehat{Q}$ is \emph{exactly} orthogonal such that \cref{thm:accuracy} focuses on the important conclusions. We address deviation from orthogonality in \cref{cor:generalized_accuracy}. This quantification is key for two points: first, characterizing how accurately eigenvalues are computed at step 6 by the Rayleigh quotients and second, designing correction steps to achieve a specified accuracy on the Schur vectors, if needed.

\begin{figure}
	\centering
	\hspace{-0.2cm}
	\includegraphics[width = 10cm]{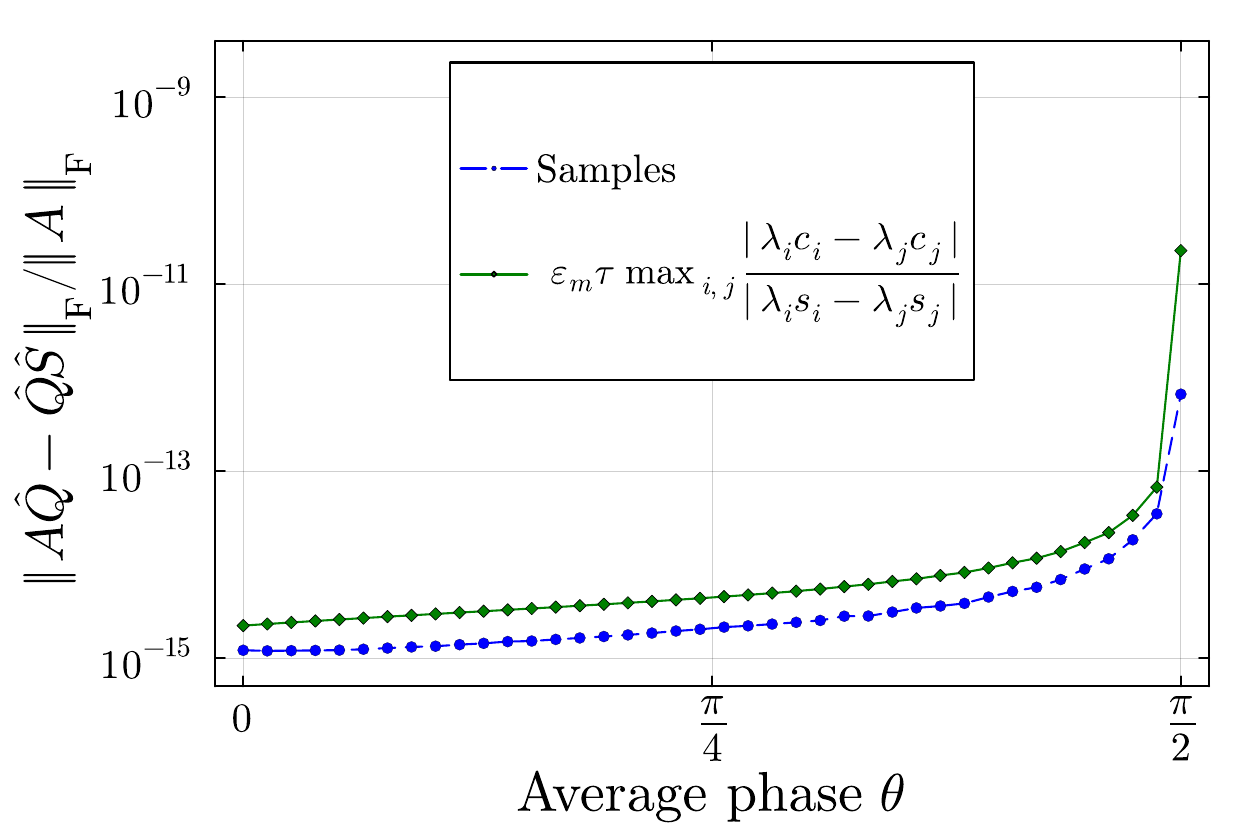}
	\hspace{-0.4cm}
	\includegraphics[width = 10cm]{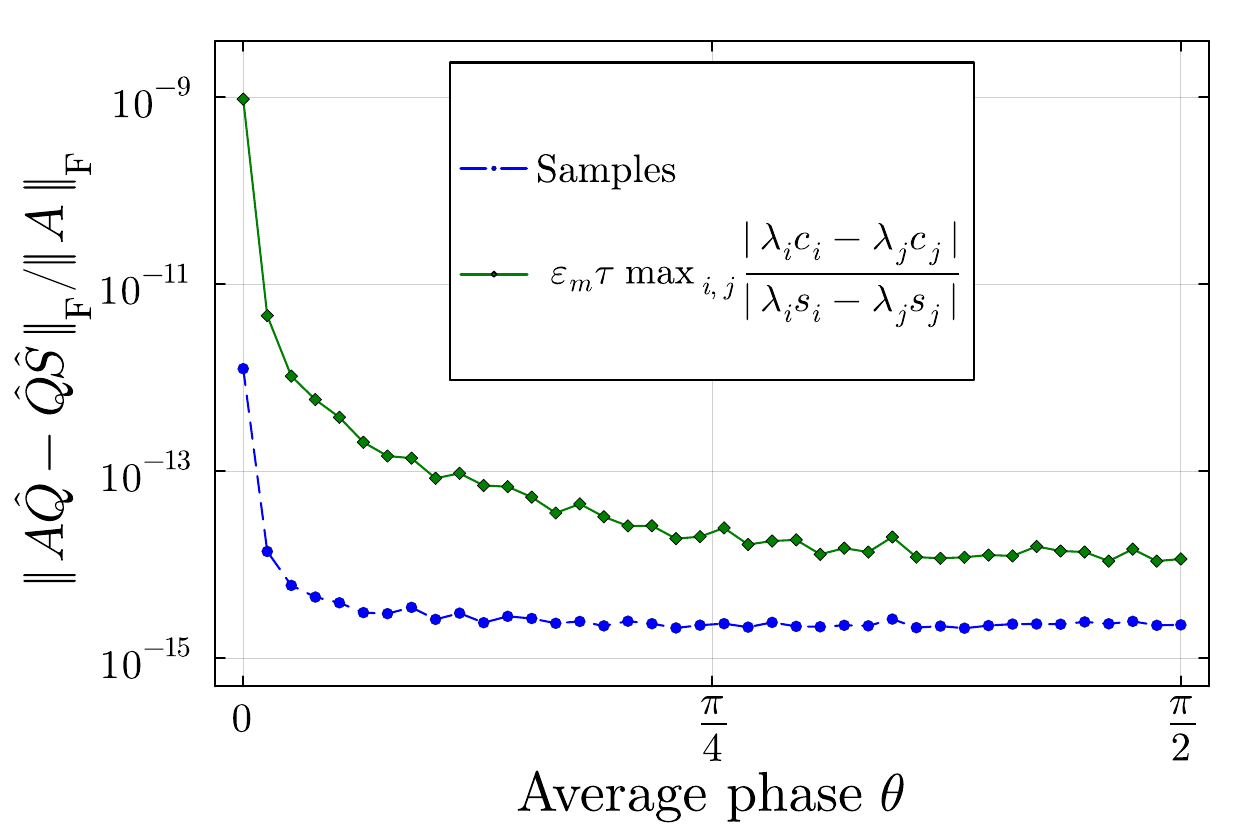}
	\vspace*{-0.3cm}
	\caption{Accuracy in double precision of \cref{alg:schur_decomposition_floating} (without any correction step~8) for $100\times 100$ orthogonal (left plot) and normal (right plot) matrices when the average phase of the eigenvalues varies from $\theta=0$ to $\theta=\frac{\pi}{2}$. The Schur vectors are sampled with Haar distribution on $\mathrm{O}(100)$. The eigenvalues of the matrices are sampled with phase uniformly distributed in $[\theta-\Delta\theta, \theta+\Delta\theta]$, with $\Delta\theta=0.005$. For normal matrices, the moduli of the eigenvalues are uniformly sampled in~$[0, 2]$. Results are averaged over 50 runs. The accuracy follows the behavior predicted by \cref{thm:accuracy}.}
	\label{fig:angles}
	\vspace*{-0.3cm}
\end{figure}

\new{We start by showing in \cref{lem:approx_backward} that it is equivalent to consider the numerical output $\widehat{\Omega} \widehat{Q} -\widehat{Q}\mathrm{skew}(\widehat{S})=\widehat{E}$ with $ \|\widehat{E}\|_\mathrm{F}\leq \widehat{\tau}\varepsilon_\mathrm{m} \|\widehat{\Omega}\|_\mathrm{F}$ from \eqref{eq:stable_skew_EVP} and the more practical relation $\Omega \widehat{Q} -\widehat{Q}\mathrm{skew}(S)=E$ with $\|E\|_\mathrm{F}\leq\tau \|A\|_\mathrm{F}$ and $\tau>0$. In the next lemma, we assume that $S$ and $\widehat{S}$ are organized such that $\mathrm{skew}(S)$ and $\mathrm{skew}(\widehat{S})$ have their eigenvalues sorted, e.g., in decreasing order.
\begin{lemma}\label{lem:approx_backward}
Let $A=QSQ^\top$ be a normal matrix, $\Omega \coloneq\mathrm{skew}(A)$ and $\widehat{\Omega}\coloneq \mathrm{skew}(\mathrm{fl}(A))$.  If the orthogonal matrix $\widehat{Q}\in\mathrm{O}(2p)$ is such that $\widehat{\Omega}\widehat{Q}-\widehat{Q}\mathrm{skew}(\widehat{S})=\widehat{E}$ is a backward stable RSD of $\widehat{\Omega}$, namely, if $\|\widehat{E}\|_\mathrm{F}\leq \widehat{\tau} \varepsilon_\mathrm{m} \|\widehat{\Omega}\|_\mathrm{F}$ for $\widehat{\tau}>0$, then it follows that
\begin{equation*}
	\|\Omega \widehat{Q} - \widehat{Q}\mathrm{skew}(S)\|_\mathrm{F}\leq (2\widehat{\tau}+4)\varepsilon_\mathrm{m} \|A\|_\mathrm{F} + \mathrm{o}(\varepsilon_\mathrm{m}).
\end{equation*}
\end{lemma}
\begin{proof}
Let us relate $\widehat{Q}$ to $\Omega$ instead of $\widehat{\Omega}$.
\begin{align*}
	\widehat{\Omega} \widehat{Q} -\widehat{Q}\mathrm{skew}(\widehat{S})&=\widehat{E}\\
	\iff \Omega \widehat{Q} -\widehat{Q}\mathrm{skew
}(S) &= \widehat{E} - (\widehat{\Omega}-\Omega)\widehat{Q}-\widehat{Q}\mathrm{skew}(S-\widehat{S}).
\end{align*}
Now let $\mathrm{skew}(\widehat{S}_2)$ be the \emph{exact} real Schur form of $\widehat{\Omega}$. Then, we have
\begin{align*}
	\|\Omega \widehat{Q} -\widehat{Q}\mathrm{skew}(S)\|_\mathrm{F} &\leq \|\widehat{E}\|_\mathrm{F} + \|(\widehat{\Omega}-\Omega)\widehat{Q}\|_\mathrm{F}+\|\widehat{Q}\mathrm{skew}(S-\widehat{S})\|_\mathrm{F}\\
	&= \|\widehat{E}\|_\mathrm{F} + \|\widehat{\Omega}-\Omega\|_\mathrm{F}+\|\mathrm{skew}(S-\widehat{S})\|_\mathrm{F}\\
	&\leq \|\widehat{E}\|_\mathrm{F} + \|\widehat{\Omega}-\Omega\|_\mathrm{F}+\|\mathrm{skew}(S-\widehat{S}_2)\|_\mathrm{F} + \|\mathrm{skew}(\widehat{S}_2-\widehat{S})\|_\mathrm{F}.
\end{align*}
By \cite[Thm.~1]{HoffmanWielandt}, we have $\|\mathrm{skew}(S-\widehat{S}_2)\|_\mathrm{F}\leq \|\widehat{\Omega}-\Omega\|_\mathrm{F}$. Moreover, since $\widehat{\Omega} = \widehat{Q}\mathrm{skew}(\widehat{S})\widehat{Q}^\top  + \widehat{E}\widehat{Q}^\top$, again by \cite[Thm.~1]{HoffmanWielandt}, we have $\|\mathrm{skew}(\widehat{S}_2-\widehat{S})\|_\mathrm{F}\leq\| \widehat{E}\widehat{Q}^\top\|_\mathrm{F} = \| \widehat{E}\|_\mathrm{F}$. By \eqref{eq:matrix_accuracy}, we have $\|\widehat{\Omega}-\Omega\|_\mathrm{F}\leq 2\varepsilon_\mathrm{m}\|A\|_\mathrm{F}$ and $\| \widehat{E}\|_\mathrm{F}\leq \widehat{\tau} \varepsilon_\mathrm{m}\|\widehat{\Omega}\|_\mathrm{F} \leq \widehat{\tau} \varepsilon_\mathrm{m}(\|\Omega\|_\mathrm{F} + \varepsilon_\mathrm{m}2\|A\|_\mathrm{F})$. This yields
\begin{align*}
	\|\Omega \widehat{Q} -\widehat{Q}\mathrm{skew}(S)\|_\mathrm{F}&\leq 2\|\widehat{E}\|_\mathrm{F} + 2\|\widehat{\Omega}-\Omega\|_\mathrm{F}\\
	&\leq 2\widehat{\tau} \varepsilon_\mathrm{m} (\|\Omega\|_\mathrm{F}+\varepsilon_\mathrm{m} 2\|A\|_\mathrm{F}) + 2(2\varepsilon_\mathrm{m} \|A\|_\mathrm{F})\\
	&\leq (2\widehat{\tau} +4) \varepsilon_\mathrm{m} \|A\|_\mathrm{F}+4\widehat{\tau} \varepsilon_\mathrm{m}^2  \|A\|_\mathrm{F}.
\end{align*}
\end{proof}}
We can now prove the stability result of \cref{thm:accuracy}.
\begin{theorem}\label{thm:accuracy}
Let $A=QSQ^\top \in\mathbb{R}^{2p\times 2p}$ be a normal matrix with eigenvalues $\lambda_j e^{\pm i\theta_j}$, $\lambda_j>0,\ \theta_j\in(0,\pi)$ for $j=1,...,p$ and $\Omega \coloneq \mathrm{skew}(A)$. Assume moreover that the imaginary parts of the eigenvalues are all distinct. Then, if the matrix $\widehat{Q}\in\mathrm{O}(2p)$ satisfies $\Omega\widehat{Q}-\widehat{Q}\mathrm{skew}(S)= E$ with $\|E\|_\mathrm{F}\leq \tau \varepsilon_\mathrm{m} \|A\|_\mathrm{F}$, we have
\begin{equation}\label{eq:final_result}
	\|A\widehat{Q}-\widehat{Q}S\|_\mathrm{F}\leq \tau \varepsilon_\mathrm{m}\|A\|_\mathrm{F} \left(1+\max_{\substack{i,j=1,..,p \\ i\neq j}}\left|\frac{\lambda_i\cos(\theta_i)-\lambda_j \cos(\theta_j)}{\lambda_i \sin(\theta_i)-\lambda_j \sin(\theta_j)}\right|\right).
\end{equation}
Moreover, if $A$ is orthogonal, then 
\begin{equation}\label{eq:tangent}
\max_{\substack{i,j=1,..,p \\ i\neq j}}\left|\frac{\lambda_i\cos(\theta_i)-\lambda_j \cos(\theta_j)}{\lambda_i \sin(\theta_i)-\lambda_j \sin(\theta_j)}\right|=\max_{\substack{i,j=1,..,p \\ i\neq j}}\left|\tan\left(\frac{\theta_i+\theta_j}{2}\right)\right|.
\end{equation}
\end{theorem}
\begin{proof}

Let us introduce that $\Omega=Q\mathrm{skew}(S)Q^\top$. This leads to
\begin{equation}\label{eq:matrix_sylvester}
	\Omega \widehat{Q} -\widehat{Q}\mathrm{skew}(S)=E\iff \mathrm{skew}(S)X - X\mathrm{skew}(S)=\widetilde{E},
\end{equation}
where we defined $X\coloneq Q^\top \widehat{Q}$ and $\widetilde{E} \coloneq Q^\top E$. Recall that $S=\left[\begin{smallmatrix}\Lambda\cos(\Theta)&-\Lambda\sin(\Theta)\\\Lambda\sin(\Theta)&\Lambda\cos(\Theta)\end{smallmatrix}\right]$ such that~\eqref{eq:matrix_sylvester} is the concatenation of $p^2=\frac{n^2}{4}$ small $2\times 2$ Sylvester equations. Indeed, for every pair $i,j=1,...,p$, using the short-hand notation $s_i\coloneq \sin(\theta_i)$ and $c_i\coloneq \cos(\theta_i)$, we have
\begin{equation}\label{eq:small_sylvester}
	\begin{bmatrix}
		0&-\lambda_i s_i\\
		\lambda_i s_i&0
	\end{bmatrix}X_{ij}-X_{ij}\begin{bmatrix}
		0&-\lambda_j s_j\\
		\lambda_j s_j&0
	\end{bmatrix}=\widetilde{E}_{ij},
\end{equation} where $X_{ij}$ and $\widetilde{E}_{ij}$ are non-overlapping $2\times2$ blocks of $X$ and $\widetilde{E}$ respectively. If $i\neq j$, by \cite[Thm.~VII.2.8]{bhatia97}, we have $\|X_{ij}\|_\mathrm{F}\leq \frac{1}{|\lambda_is_i-\lambda_js_j|}\|\widetilde{E}_{ij}\|_\mathrm{F}$. Therefore, acknowledging that $\|A\widehat{Q}-\widehat{Q}S\|_\mathrm{F}=\|SX-XS\|_\mathrm{F}\leq \|\mathrm{sym}(S)X-X\mathrm{sym}(S)\|_\mathrm{F}+\|\widetilde{E}\|_\mathrm{F}$, we obtain
\begin{align}
\nonumber
	\|A\widehat{Q}-\widehat{Q}S\|_\mathrm{F}^2&=\sum_{i,j=1}^{p} \left\|\lambda_i \begin{bmatrix}c_i&-s_i\\ s_i&c_i\end{bmatrix} X_{ij}-X_{ij}\lambda_j \begin{bmatrix}c_j&-s_j\\ s_j&c_j\end{bmatrix} \right\|_\mathrm{F}^2\\
	\nonumber
	&\leq \sum_{i,j=1}^p\left(|\lambda_i c_i-\lambda_j c_j|\left\|X_{ij}\right\|_\mathrm{F}+\|\widetilde{E}_{ij}\|_\mathrm{F}\right)^2\\
	%\nonumber
	%&\quad +\sum_{i,j=1, i\neq j}^p 2|\lambda_i c_i-\lambda_j c_j|\|X_{ij}\|_\mathrm{F}\|\widetilde{E}_{ij}\|_\mathrm{F}\\
	%\nonumber
	%&\leq \sum_{i,j=1, i\neq j}^p\left|\frac{\lambda_i c_i-\lambda_j c_j}{\lambda_i s_i-\lambda_j s_j}\right|^2\|\widetilde{E}_{ij}\|_\mathrm{F}^2+ \sum_{i,j=1}^p\|\widetilde{E}_{ij}\|_\mathrm{F}^2\\
	\nonumber
	&\leq \left(1+\max_{\substack{i,j=1,..,p \\ i\neq j}}\left|\frac{\lambda_i c_i-\lambda_j c_j}{\lambda_i s_i-\lambda_j s_j}\right| \right)^2 \sum_{i,j=1}^p\|\widetilde{E}_{ij}\|_\mathrm{F}^2.
\end{align}
The claim \eqref{eq:final_result} follows directly by observing that $\sqrt{\sum_{i,j=1}^p\|\widetilde{E}_{ij}\|_\mathrm{F}^2}=\|\widetilde{E}\|_\mathrm{F}=\|E\|_\mathrm{F}\leq \tau \varepsilon_\mathrm{m} \|A\|_\mathrm{F}$. If $A$ is orthogonal, then $\lambda_j=1$ for $j=1,...,p$. \eqref{eq:tangent} follows from simple trigonometric identities.
\end{proof}
Generalizing \cref{thm:accuracy} when the matrix $\widehat{Q}$ is not exactly orthogonal is simple and does not change its conclusion. This is shown in \cref{cor:generalized_accuracy}.
\begin{corollary}\label{cor:generalized_accuracy}
Assume in \cref{thm:accuracy} that $\widehat{Q}\notin \mathrm{O}(2p)$ but $d_{\mathrm{O}(n)}(\widehat{Q})\leq \kappa \varepsilon_\mathrm{m} $ with $\kappa>0$. Then the matrix $\widehat{Q}_*\in \mathrm{arg}\min_{Q\in\mathrm{O}(n)}\|Q-\widehat{Q}\|_\mathrm{F}$ satisfies the hypotheses of \cref{thm:accuracy} and $\|A\widehat{Q} - \widehat{Q}S\|_\mathrm{F}\leq 2\kappa \varepsilon_\mathrm{m} \|A\|_\mathrm{F} + \|A\widehat{Q}_* - \widehat{Q}_*S\|_\mathrm{F}$.
\end{corollary}
	\begin{proof}
	First, we have 
	\begin{align*}
		\|A\widehat{Q} - \widehat{Q}S\|_\mathrm{F}&\leq \|A(\widehat{Q}-\widehat{Q}_*) - (\widehat{Q}-\widehat{Q}_*)S\|_\mathrm{F} +  \|A\widehat{Q}_* - \widehat{Q}_*S\|_\mathrm{F}\\
		&\leq 2\|\widehat{Q}-\widehat{Q}_*\|_\mathrm{F}\|A\|_\mathrm{F} +  \|A\widehat{Q}_* - \widehat{Q}_*S\|_\mathrm{F}\\
		&\leq 2\kappa \varepsilon_\mathrm{m} \|A\|_\mathrm{F} + \|A\widehat{Q}_* - \widehat{Q}_*S\|_\mathrm{F}.
	\end{align*}
	Moreover, $\widehat{Q}_*$ satisfies the hypotheses of \cref{thm:accuracy}:
	\begin{align*}
		\|\Omega\widehat{Q}_*-\widehat{Q}_* \mathrm{skew}(S)\|_\mathrm{F}&\leq \|\widehat{\Omega} (\widehat{Q}_*-\widehat{Q})-(\widehat{Q}_*-\widehat{Q}) \mathrm{skew}(S)\|_\mathrm{F}+\|\widehat{\Omega} \widehat{Q}-\widehat{Q} \mathrm{skew}(S)\|_\mathrm{F}\\
		&\leq\|\widehat{Q}_*-\widehat{Q}\|_\mathrm{F}(\|\Omega\|_\mathrm{F} +  \|\mathrm{skew}(S)\|_\mathrm{F})+\tau \varepsilon_\mathrm{m} \|A\|_\mathrm{F}\\
		&\leq (2\kappa + \tau )\varepsilon_\mathrm{m} \|A\|_\mathrm{F}.
	\end{align*}
	\end{proof}
	\cref{thm:accuracy} is very informative on the accuracy that \cref{alg:schur_decomposition_floating} will achieve when it is allowed to take its fastest route by skipping steps~5 and~7. We recall that this setting occurs with very high probability for Haar-distributed orthogonal matrices, and more generally, for random orthogonal matrices. In particular, if the bound in \eqref{eq:final_result} is not satisfactory to the user, i.e., if $\max_{\substack{i,j=1,..,p \\ i\neq j}}\left|\frac{\lambda_i\cos(\theta_i)-\lambda_j \cos(\theta_j)}{\lambda_i \sin(\theta_i)-\lambda_j \sin(\theta_j)}\right|\geq t$, the sets of eigenvalues responsible for the loss of accuracy can be identified and the associated invariant subspaces can be corrected. We detail this approach in \cref{subsec:accuracy_cont}. In \cref{fig:angles}, we present the results of a numerical experiment corroborating the bounds of \cref{thm:accuracy} and assessing their relevance.
\subsection{Sensitivity of the invariant subspaces: Set~2}
From now on, we use the notation $c,s$ for $\cos(\theta)$ and $\sin(\theta)$. Assume still that $A=QSQ^\top\in\mathbb{R}^{2p\times 2p}$ ($r=0$) is a normal matrix with eigenvalues $\lambda_j e^{\pm i\theta_j}$, $\lambda_j>0,\ \theta_j\in(0,\pi)$ for $j=1,...,p$. We next consider the case of eigenvalues that have repeated imaginary parts (\wordref{set2}{Set~2}). Assume a set of indices $\mathcal{J}\subset \{1,2,..,p\}$, $|\mathcal{J}|\eqcolon m$, such that for all $j\in\mathcal{J}$, we have $ \lambda_{j} s_{j}=\sigma>0$. In other words, we can numerically compute $\widehat{V}$ such that $\widehat{\Omega}\widehat{V}-\widehat{V}\widehat{\sigma} \left[\begin{smallmatrix}0&-I_m\\
	I_m&0\end{smallmatrix}\right]=\widehat{E}$ with $\|\widehat{E}\|_\mathrm{F}\lesssim\varepsilon_\mathrm{m} \|\widehat{\Omega}\|_\mathrm{F}$. In view of \cref{cor:generalized_accuracy}, w.l.o.g., we can consider $\widehat{V}\in\mathrm{St}(n,2m)$. We set $\mathcal{I}\coloneq \{1,2,..,p\}\setminus \mathcal{J}$. The proof of \cref{thm:accuracy} fails when \cite[Thm.~VII.2.8]{bhatia97} is used after \eqref{eq:small_sylvester} because, for $j_1,j_2\in\mathcal{J}$, $\lambda_{j_1} s_{j_1}=\lambda_{j_2}s_{j_2}$. %Indeed, we know from \cref{subsec:case2} that $\widehat{V}$ needs to be post-multiplied by a matrix $R\in\mathrm{O}(2m)$ to approximate the Schur vectors of $A$. 
	
	As expected, the accuracy of the invariant subspace computation depends on the separation between $\pm i\sigma$ and the other eigenvalues values of $\Omega$. In \cref{subsec:accuracy_cont} we use \cref{thm:accuracy_repeated} to design an algorithm that satisfies a desired accuracy.
\begin{theorem}\label{thm:accuracy_repeated}
Let $A=QSQ^\top \in\mathbb{R}^{2p\times 2p}$ be a normal matrix with eigenvalues $\lambda_j e^{\pm i\theta_j}$, $\lambda_j>0,\ \theta_j\in(0,\pi)$ for $j=1,...,p$ and $\Omega \coloneq \mathrm{skew}(A)$. Assume there is $\mathcal{J}\subset\{1,...,p\}$, $|\mathcal{J}|\eqcolon m$, such that for all $j\in\mathcal{J}$, we have $ \lambda_{j} s_{j}=\sigma$. Set $\mathcal{I}\coloneq\{1,...,p\}\setminus \mathcal{J}$. If the matrix $\widehat{V}\in\mathrm{St}(2p,2m)$ satisfies $\Omega\widehat{V}-\widehat{V}\sigma\left[\begin{smallmatrix}0&-I_m\\
	I_m&0\end{smallmatrix}\right]=E$ with $\|E\|_\mathrm{F}\leq \tau \varepsilon_\mathrm{m} \|A\|_\mathrm{F}$, then we have
\begin{align*}
	&\min_{R\in\mathrm{O}(2m)}\|A(\widehat{V}R)-(\widehat{V}R)\left[\begin{smallmatrix}
		D_m&-\sigma I_m\\
		\sigma I_m&D_m
		\end{smallmatrix}\right]\|_\mathrm{F} \\&\leq\tau \varepsilon_\mathrm{m} \|A\|_\mathrm{F} \left(1 + \frac{\max\limits_{i\in\mathcal{I},j\in\mathcal{J}}|\lambda_i c_i-\lambda_j c_j|}{\min\limits_{i\in\mathcal{I}}|\lambda_i s_i-\sigma|}+ \varepsilon_\mathrm{m} \tau\|A\|_\mathrm{F}\frac{ \max\limits_{j_1,j_2\in\mathcal{J}}|\lambda_{j_1} c_{j_1}-\lambda_{j_2} c_{j_2}|}{\min\limits_{i\in\mathcal{I}}|\lambda_i s_i-\sigma|^2}\right).
\end{align*}
\end{theorem}
\begin{proof}	
	Define $X\coloneq Q^\top \widehat{V}$, then it follows that
	\begin{equation*}
	\mathrm{skew}(S)X-X\sigma \begin{bmatrix}0&-I_m\\
	I_m&0\end{bmatrix}=Q^\top E\eqcolon \widetilde{E}.
	\end{equation*}
There exist two permutations $P_\mathrm{l}\in\mathrm{O}(n),P_\mathrm{r}\in\mathrm{O}(2m) $ such that $P_\mathrm{l}^\top X P_\mathrm{r}=\left[\begin{smallmatrix}X_{\mathcal{I},\mathcal{J}}  \\X_{\mathcal{J},\mathcal{J}}\end{smallmatrix}\right]$ where the block $X_{\mathcal{I},\mathcal{J}}$ comprises all $2\times 2$ blocks $X_{ij}$ such that for $i\in\mathcal{I}, j\in\mathcal{J}$. We have
\begin{equation}\label{eq:sylvester_x}
\begin{bmatrix}
		0&-\lambda_i s_i\\
		\lambda_i s_i&0
	\end{bmatrix}X_{ij}-X_{ij}\begin{bmatrix}
		0&-\sigma\\
		\sigma&0
	\end{bmatrix}=(P_\mathrm{l}^\top \widetilde{E} P_\mathrm{r})_{ij}.
\end{equation} 
If it hold that $\widetilde{E}=0$, then we would have $X_{\mathcal{I},\mathcal{J}}=0 $ and $X_{\mathcal{J},\mathcal{J}}\in\mathrm{O}(2m)$. However, since $\widetilde{E}\neq 0$ and $X\in\mathrm{St}(n,2m)$, we have
\begin{align*}
 X_{\mathcal{I},\mathcal{J}}^\top X_{\mathcal{I},\mathcal{J}} +X_{\mathcal{J},\mathcal{J}}^\top X_{\mathcal{J},\mathcal{J}} = I_{2m} \Longrightarrow \| X_{\mathcal{J},\mathcal{J}}^\top X_{\mathcal{J},\mathcal{J}}-I_{2m}\|_\mathrm{F}\leq \| X_{\mathcal{I},\mathcal{J}}\|_\mathrm{F}^2\\
 \text{and by \eqref{eq:sylvester_x},}\quad \| X_{\mathcal{I},\mathcal{J}}\|_\mathrm{F}^2\leq \sum_{i\in\mathcal{I}, j\in\mathcal{J}} \frac{1}{|\lambda_i s_i-\sigma|^2}\|(P_\mathrm{l}^\top \widetilde{E} P_\mathrm{r})_{ij}\|_\mathrm{F}^2\leq \frac{\|E\|_\mathrm{F}^2}{\min_{i\in\mathcal{I}}|\lambda_i s_i-\sigma|^2}.
\end{align*}
The last inequality is obtained using~\cite[Thm.~VII.2.8]{bhatia97}. The deviation of $X_{\mathcal{J},\mathcal{J}}$ from orthogonality can thus be quantified. In addition, by taking the SVD $ X_{\mathcal{J},\mathcal{J}}=U_X\Sigma_X V_X^\top$, we can write $\| X_{\mathcal{J},\mathcal{J}}^\top X_{\mathcal{J},\mathcal{J}}-I_{2m}\|_\mathrm{F}=\|\Sigma_X^2-I_{2m}\|_\mathrm{F}\geq \|\Sigma_X-I_{2m}\|_\mathrm{F}$, where the last inequality stands because $|\alpha^2-1|\geq |\alpha-1|$ for all $\alpha\geq 0$. Define $\left[\begin{smallmatrix}\widetilde{S}_{2(p-m)} &0 \\ 0 &\widetilde{S}_{2m}\end{smallmatrix}\right]\coloneq P_\mathrm{l}^\top S P_\mathrm{l}$, then for all $R\in \mathrm{O}(2m)$, we have 
\begin{align}
\nonumber
	\|A\widehat{V}R-\widehat{V}R\widetilde{S}_{2m}\|_\mathrm{F}&\leq\|\mathrm{sym}(S)XR-XR\mathrm{sym}(\widetilde{S}_{2m})\|_\mathrm{F}+\|E\|_\mathrm{F}\\
	\label{eq:master_equation}
	&\leq \|\mathrm{sym}(\widetilde{S}_{2(p-m)})X_{\mathcal{I},\mathcal{J}}R-X_{\mathcal{I},\mathcal{J}}R\mathrm{sym}(\widetilde{S}_{2m})\|_\mathrm{F}\\
	\nonumber
	&\quad +\|\mathrm{sym}(\widetilde{S}_{2m})X_{\mathcal{J},\mathcal{J}}R-X_{\mathcal{J},\mathcal{J}}R\mathrm{sym}(\widetilde{S}_{2m})\|_\mathrm{F}+\|E\|_\mathrm{F}.
	%\nonumber
	%\|\widetilde{S}_{2m}X_{\mathcal{J},\mathcal{J}}R-X_{\mathcal{J},\mathcal{J}}R\widetilde{S}_{2m}\|_\mathrm{F}^2+\|\widetilde{S}_{2k}X_{\mathcal{I},\mathcal{J}}R-X_{\mathcal{I},\mathcal{J}}R\widetilde{S}_{2m}\|_\mathrm{F}^2
	%\label{eq:master_equation}
	%&+\sum_{i\in\mathcal{I},j\in\mathcal{J}} \left\|\lambda_i G_i(XR)_{ij}-(XR)_{ij}\lambda_j G_j \right\|_\mathrm{F}^2.
\end{align}
The first term of \eqref{eq:master_equation} can be bounded by
\begin{align}
\nonumber
	\|\mathrm{sym}(\widetilde{S}_{2(p-m)})X_{\mathcal{I},\mathcal{J}}R-X_{\mathcal{I},\mathcal{J}}R\mathrm{sym}(\widetilde{S}_{2m})\|_\mathrm{F}&\leq  \max_{i\in\mathcal{I},j\in\mathcal{J}}|\lambda_i c_i-\lambda_j c_j|\|X_{\mathcal{I},\mathcal{J}}R\|_\mathrm{F}\\
	\label{eq:slave1}
	&\leq \frac{\max_{i\in\mathcal{I},j\in\mathcal{J}}|\lambda_i c_i-\lambda_j c_j|}{\min_{i\in\mathcal{I}}|\lambda_i s_i-\sigma|}\|E\|_\mathrm{F}.
\end{align}
Moreover, the second term can also be bounded by
\begin{align}
\nonumber
&\|\mathrm{sym}(\widetilde{S}_{2m})X_{\mathcal{J},\mathcal{J}}R-X_{\mathcal{J},\mathcal{J}}R\mathrm{sym}(\widetilde{S}_{2m})\|_\mathrm{F}\\
%\nonumber
%&=\|\widetilde{S}_{2m}(X_{\mathcal{J},\mathcal{J}}R-I_{2m})-(X_{\mathcal{J},\mathcal{J}}R-I_{2m})\widetilde{S}_{2m}\|_\mathrm{F}\\
\nonumber
&= \|\mathrm{sym}(\widetilde{S}_{2m})(X_{\mathcal{J},\mathcal{J}}R-I_{2m})-(X_{\mathcal{J},\mathcal{J}}R-I_{2m})\mathrm{sym}(\widetilde{S}_{2m})\|_\mathrm{F}\\
\label{eq:slave2}
&\leq  \max_{j_1,j_2\in\mathcal{J}}|\lambda_{j_1} c_{j_1}-\lambda_{j_2} c_{j_2}|\|X_{\mathcal{J},\mathcal{J}}R-I_{2m}\|_\mathrm{F}.
\end{align}
In particular, for $R=V_X U_X^\top$ which minimizes $\|X_{\mathcal{J},\mathcal{J}}R-I_{2m}\|_\mathrm{F}$, we have 
\begin{equation}\label{eq:slave3}
\|X_{\mathcal{J},\mathcal{J}}R-I_{2m}\|_\mathrm{F}=\|\Sigma_X-I_{2m}\|_\mathrm{F}\leq \frac{\|E\|_\mathrm{F}^2}{\min_{i\in\mathcal{I}}|\sigma-\lambda_i s_i|^2}.  
\end{equation}
The theorem follows directly by inserting \eqref{eq:slave1}, \eqref{eq:slave2} and \eqref{eq:slave3} in \eqref{eq:master_equation} and acknowledging that  $\|E\|_\mathrm{F}\leq \tau \varepsilon_\mathrm{m} \|A\|_\mathrm{F}$.
\end{proof}
 
\subsection{Sensitivity of the invariant subspaces: Set~3} The case of real eigenvalues (\wordref{set3}{Set~3}) can be handled in a similarly to the case of repeated imaginary parts. Indeed, the proof of \cref{thm:accuracy_real} is very close to that of \cref{thm:accuracy_repeated}, if $\left[\begin{smallmatrix}
		D_m&-\sigma I_m\\
		\sigma I_m&D_m
		\end{smallmatrix}\right]$ is replaced by $\breve{\Lambda}$. For this reason, the proof is given in \cref{app:theorem_proof}.
\begin{theorem}\label{thm:accuracy_real}
Let $A=QSQ^\top \in\mathbb{R}^{n\times n}$ be a normal matrix with real eigenvalues $\breve{\lambda}_k$ for $k\in\mathcal{K}\coloneq\{1,...,r\}$ and $\Omega \coloneq \mathrm{skew}(A)$. Moreover, assume $\lambda_j e^{\pm i\theta_j}$ is an eigenvalue of $A$ for $j \in\mathcal{I}$. If the matrix $\widehat{V}\in\mathrm{St}(n,r)$ satisfies $\Omega\widehat{V}=\widehat{E}$ with $\|E\|_\mathrm{F}\leq\tau \varepsilon_\mathrm{m} \|A\|_\mathrm{F}$, then we have
\begin{align*}
	&\min_{\breve{R}\in\mathrm{O}(2m)}\|A\widehat{V}\breve{R}-\widehat{V}\breve{R}\breve{\Lambda}\|_\mathrm{F} \\&\leq \tau \varepsilon_\mathrm{m}\|A\|_\mathrm{F} \left(1+ \frac{\max\limits_{i\in\mathcal{I},k\in\mathcal{K}}|\lambda_i c_i-\breve{\lambda}_k|}{\min\limits_{i\in\mathcal{I}}|\lambda_i s_i|}+\varepsilon_\mathrm{m} \tau\|A\|_\mathrm{F}\frac{\max\limits_{k_1,k_2\in\mathcal{K}}|\breve{\lambda}_{k_1}-\breve{\lambda}_{k_2}|}{\min\limits_{i\in\mathcal{I}}|\lambda_i s_i|^2}\right).
\end{align*}
\end{theorem}
\begin{proof}
	See \cref{app:theorem_proof}.
\end{proof}
\subsection{Analysis of steps 5 to 8}\label{subsec:accuracy_cont}
	
	The previous subsections determined the sensitivity of the computation of the Schur vectors of $A$ by the Schur vectors of $\Omega$. We can now provide conditions on steps~5 to~8 of \cref{alg:schur_decomposition_floating} to achieve a desired accuracy.
	
	\textbf{Step 5} requires identifying every approximately invariant subspace of $\Omega$ spanned by $\widehat{V}$, i.e., $\widehat{V}$ such that $\| \Omega \widehat{V} - \widehat{V}\sigma J_{2m}\|_\mathrm{F}\lesssim \varepsilon_\mathrm{m}\|A\|_\mathrm{F}$, $m > 1$. In view of \cref{cor:generalized_accuracy}, we assume $\widehat{V}\in\mathrm{St}(n,2m)$. To recover the Schur vectors of $A$, we need to compute $M\coloneq \widehat{V}^\top A\widehat{V}$. In exact arithmetic, $2m$ is the multiplicity of the singular value $\sigma$ of~$\Omega$. In floating point arithmetic, $2\widehat{m}$ must be the size of the cluster of singular values of~$\widehat{\Omega}$. 
	If we take $\widehat{V}$ as the $2\widehat{m}$ computed Schur vectors associated to a $\delta$-cluster of $\widehat{\Omega}$, in view of \cref{thm:accuracy_repeated}, it holds that $\|A\widehat{V} - \widehat{V} M\|_\mathrm{F}\lesssim \frac{\varepsilon_\mathrm{m}}{\delta}\|A\|_\mathrm{F}$. If $\delta$ is chosen such that the $\delta$-cluster is well-separated from the others and if the RSD $M \approx R \widetilde{S}_{2m} R^\top$ is such that $\|M - R \widetilde{S}_{2m} R^\top\|_\mathrm{F}\lesssim \varepsilon_\mathrm{m} \|A\|_\mathrm{F}$, then after step 5, we have  $\|A\widehat{V}R - \widehat{V}  R \widetilde{S}_{2m}\|_\mathrm{F}\lesssim \frac{\varepsilon_\mathrm{m}}{\delta}\|A\|_\mathrm{F}$ and $d_{\mathrm{St}(n,2m)}(\widehat{V}R)\lesssim \varepsilon_\mathrm{m}$. Notice here that the effective accuracy that is obtained is in fact function of the largest $\delta$ such that the $\delta$-cluster of $\widehat{\Omega}$ remains the same, rather than the user specified value $\delta$. As mentioned in \cref{sec:normal_EVP}, several stable algorithms exist for \wordref{oracle3}{Routine~3} at step~5.2. 
	%\begin{equation*}
	%	\widetilde{A}_m + \Delta \widetilde{A}_m = \widehat{R} \widehat{S}_{2m} \widehat{R}^\top  \text{ with } \|\Delta \widetilde{A}_m\|_2 < c \varepsilon_\mathrm{m} \|\Delta \widetilde{A}_m\|_2,
	%\end{equation*}
	%such that $R^\top R=I_{2m}$ and $\|\widehat{S}_{2m}-S_{2m}\|_2 < c \varepsilon_\mathrm{m} \|S_{2m}\|_2$.
	
	\textbf{Step 6} obtains the real parts of the eigenvalues by computing Rayleigh quotients. It is known that the Rayleigh quotient offers the square of the accuracy of the eigenvector, see, e.g., \cite[Lecture~27]{trefethen97}. Therefore, by \cref{thm:accuracy}, step 6 provides a worst-case $\mathcal{O}(\varepsilon_\mathrm{m}\|A\|_\mathrm{F})$-absolute accuracy on the real parts of the eigenvalues if $\max_{i,j}\left|\frac{\lambda_i c_i - \lambda_j c_j}{\lambda_i s_i - \lambda_j s_j}\right|\lesssim \frac{1}{\sqrt{\varepsilon_\mathrm{m}}}$. It suggests considering $\delta$-clusters of $\Omega$ with $\delta\leq\sqrt{\varepsilon_\mathrm{m}}$ to guarantee high absolute accuracy on the eigenvalues. 
	
	\textbf{Step 7} computes the real eigenvalues and the eigenvectors of $A$. This starts by estimating $r$ by $\widehat{r}$.  By \cref{thm:accuracy_real}, considering a $\delta_r$-cluster of the eigenvalues of $\widehat{\Omega}$ around zero will provide a worst-case $\mathcal{O}(\frac{\varepsilon_\mathrm{m}}{\delta_r}\|A\|_\mathrm{F})$-absolute accuracy on the real eigenspace of~$A$. It is essential not to underestimate $r$, i.e., to have $\widehat{r}\geq r$. Therefore, $\delta_r$ should always be increased until the $\delta_r$-cluster around zero is well separated from the other eigenvalues, if needed at the cost of the symmetry of $H\coloneq \widehat{Q}_{\widehat{p}+1:\widehat{p}+\widehat{r}}^\top A\widehat{Q}_{\widehat{p}+1:\widehat{p}+\widehat{r}}$. If $\|\mathrm{skew}(H)\|_\mathrm{F}\leq \varepsilon_\mathrm{m} \|H\|_\mathrm{F}$, then $H$ can be considered exactly symmetric without loss of accuracy since the skew-symmetric part can be considered as a negligible perturbation. If $H$ is symmetric, the decomposition at step~7.2a using, e.g.,  \href{https://www.netlib.org/lapack/explore-html/d1/d56/group__heevr_gaa334ac0c11113576db0fc37b7565e8b5.html#gaa334ac0c11113576db0fc37b7565e8b5}{syevr}, yields
	\begin{equation*}
		 H + \Delta H = \breve{R}\breve{\Lambda}\breve{R}^\top\quad\text{with}\quad \|\Delta H\|_\mathrm{F}\lesssim \varepsilon_\mathrm{m}\|H\|_\mathrm{F},
	\end{equation*}
	  with $d_{\mathrm{O}(\widehat{r})}(\breve{R})\lesssim \varepsilon_\mathrm{m}$ \cite{lapack99}. This yields $\|A\widehat{Q}_{\widehat{p}+1:\widehat{p}+\widehat{r}} \breve{R} - \widehat{Q}_{\widehat{p}+1:\widehat{p}+\widehat{r}} \breve{R}\breve{\Lambda}\|_\mathrm{F}\lesssim \frac{\varepsilon_\mathrm{m}}{\delta_r} \|A\|_\mathrm{F}$ with $d_{\mathrm{St}(n,\widehat{r})}(\widehat{Q}_{p+1:p+\widehat{r}}\breve{R})\lesssim \varepsilon_\mathrm{m} $. Only in certain circumstances will symmetric eigensolvers achieve relative accuracy \cite{DemmelVeseli92,dhillon1997new,EisenstatIpsen98,VESELIC199381}. %, notably for tridiagonal matrices \cite{dhillon1997new} and positive definite matrices \cite{}. 
If $H$ is not numerically symmetric, it means that $\widehat{r}$ overestimated $r$, \wordref{oracle3}{Routine~3} should then be preferred to \wordref{oracle2}{Routine~2} to maintain accuracy.
	
	\new{\textbf{Step~8} can finally be considered to improve the accuracy of the RSD to a desired tolerance $\mu =\varepsilon_\mathrm{m} t$, $t\geq 1$ if this accuracy has not yet been attained. By combining \cref{thm:accuracy,thm:accuracy_repeated,thm:accuracy_real}, we deduce that $\|A\widehat{Q} -\widehat{Q}\widehat{S}\|_\mathrm{F}\lesssim \left(1+ \max_{i\neq j, k}\{\frac{|\lambda_i c_i - \lambda_j c_j|}{|\lambda_i s_i - \lambda_j s_j|},\frac{|\breve{\lambda}_k - \lambda_j c_j|}{| \lambda_j s_j|}\}\right)\varepsilon_\mathrm{m} \|A\|_\mathrm{F}$. If $\max_{i\neq j, k}\{\frac{|\lambda_i c_i - \lambda_j c_j|}{|\lambda_i s_i - \lambda_j s_j|},\frac{|\breve{\lambda}_k - \lambda_j c_j|}{| \lambda_j s_j|}\}\geq t$, the groups of eigenvalues responsible for making the constant larger than $t$ can be identified. Then, to restore the desired accuracy, one should repeat step~5 with $\delta=\frac{1}{t}$ on the invariant subspaces associated with $\frac{1}{t}$-clusters of $\widehat{\Omega}$. Because the matrix $\widehat{Q}$ is already close to the Schur vectors of $A$, the matrix $M$ at step~5 is a small perturbation of the Schur form. For such matrices, a Jacobi algorithm such as \cite{ZhouBrent03} is efficient since it converges quadratically. One sweep is enough in practice to obtain the desired accuracy. However, for specific matrices and $t$ too small, the $\frac{1}{t}$-clusters may include many eigenvalues. The user should be aware that step~8 may then require a significant additional computational work. In double precision, the authors recommend using $t\geq 100$ based on empirical experience. If the output of \wordref{oracle3}{Routine~3} satisfies $\|M - R \widetilde{S}_{2m} R^\top\|_\mathrm{F}\lesssim \varepsilon_\mathrm{m} \|A\|_\mathrm{F}$, then finally, we obtain the decomposition $\|A\widehat{Q}-\widehat{Q}\widehat{S}\|_\mathrm{F}\lesssim\mu \|A\|_\mathrm{F}$ with $d_{\mathrm{O}(n)}(\widehat{Q})\lesssim \varepsilon_\mathrm{m}$.} %The second approach applies on the complete decomposition and is suited if $t$ is close to $1$. It refines the RSD of the full matrix $\widehat{Q}^\top A \widehat{Q}$. In both cases, the factorized matrices are small perturbations of the real Schur form. Therefore, the Jacobi method for normal matrices should converge in only a few sweeps~\cite{ZhouBrent03}, as shown in the experiments of \Cref{tab:accuracy}. Since the method is locally quadratically convergent, one sweep is often enough in practice.
	\section{Numerical experiments on the accuracy}\label{sec:exp_accuracy} 
\begin{figure}
	\centering
	\includegraphics[width = 10cm]{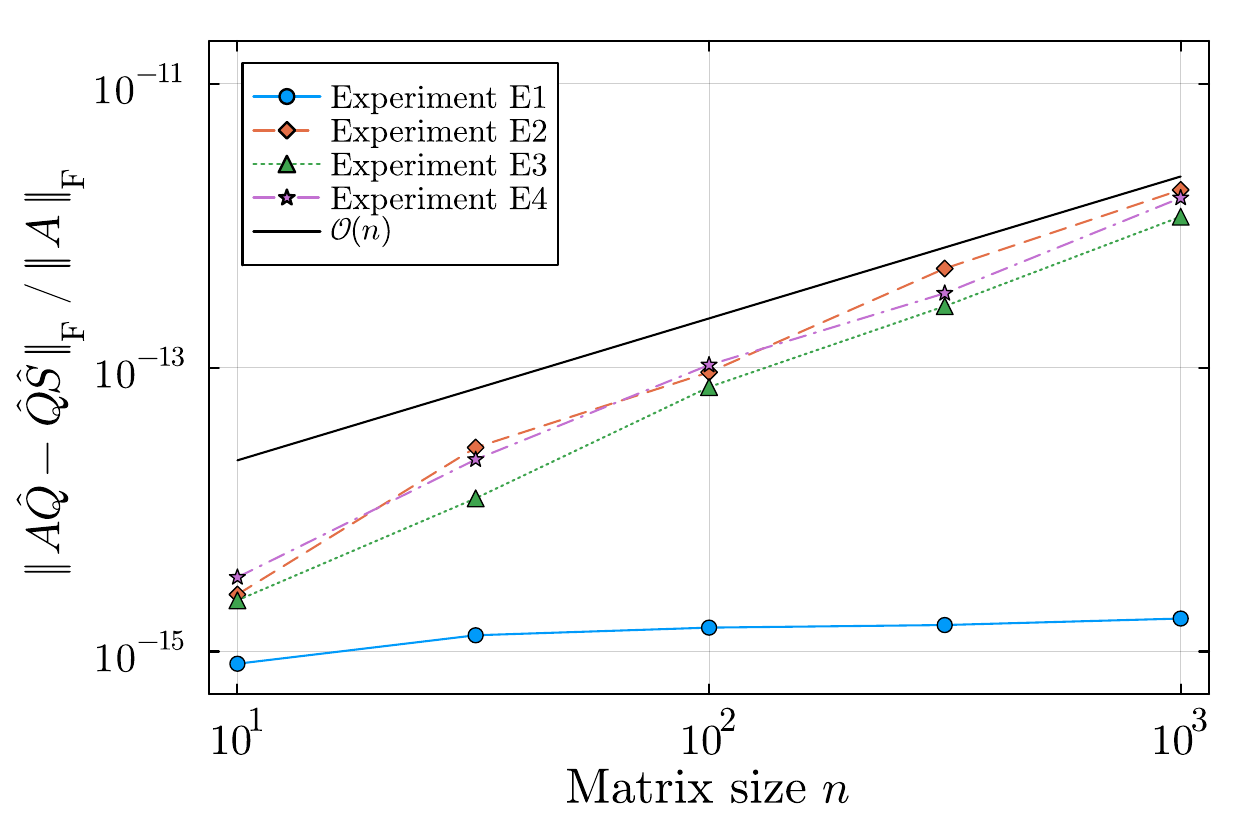}
	\vspace*{-0.3cm}
	\caption{\new{Evolution of the accuracy in double precision of \texttt{nrmschur} (\cref{alg:schur_decomposition_floating}) without correction step~8. The matrix size $n$ varies from $10$ to $1000$. The plots shows results for experiments E1 to E4 from \cref{sec:exp_accuracy}. The results are averaged over 100 runs.}}
	\label{fig:evolution_accuracy}
	\vspace*{-0.3cm}
\end{figure}

	\begin{table}[ht]
	 
	 \centering
	 %\small
	 
	 \begin{tabular}{||c||c|c||c|c||c|c||}
	 \hline
	 $n$&\multicolumn{2}{l||}{\large$\quad\quad\quad\frac{\|A\widehat{Q} - \widehat{Q}\widehat{S}\|_\mathrm{F}}{\|A\|_\mathrm{F}}$}&\multicolumn{2}{l||}{\large $\quad\quad\frac{\|\widehat{Q}^\top \widehat{Q}- I\|_\mathrm{F}}{\sqrt{n}}\quad\quad$}&\multicolumn{2}{l||}{\large $\quad \ \frac{\|\mathrm{diag}(S-\widehat{S})\|_\mathrm{F}}{1+\|\mathrm{diag}(S)\|_\mathrm{F}}$}\\
	 \hline
	 &\texttt{nrmschur}&\texttt{schur}&\texttt{nrmschur}&\texttt{schur}&\texttt{nrmschur}&\texttt{schur}\\
	 \hline
	 \multicolumn{7}{||l||}{E1: $\lambda = 1,\ \theta\sim \mathcal{U}(0, \frac{\pi}{4})$: Best accuracy scenario on $\mathrm{SO}(n)$.}\\
 
	 \hline
	10 & 8.0e-16 & 1.1e-15 & 6.7e-16 & 1.2e-15 & 3.8e-16 & 4.6e-16 \\
32 & 1.3e-15 & 1.9e-15 & 1.2e-15 & 2.0e-15 & 6.6e-16 & 8.1e-16 \\
100 & 1.5e-15 & 3.5e-15 & 1.6e-15 & 3.7e-15 & 6.6e-16 & 1.6e-15 \\
316 & 1.5e-15 & 7.2e-15 & 2.0e-15 & 7.8e-15 & 5.8e-16 & 4.3e-15 \\
1000 & 1.7e-15 & 2.3e-14 & 2.8e-15 & 2.3e-14 & 6.4e-16 & 1.9e-14 \\
	 \hline
	 \multicolumn{7}{||l||}{E2: $\lambda \sim \mathcal{U}(0,2),\ \theta\sim \mathcal{U}(0,\pi)$: Random normal scenario.}\\
	 \hline
	 10 & 2.7e-15 & 1.4e-15 & 6.2e-16 & 1.3e-15 & 3.4e-16 & 5.1e-16 \\
32 & 1.6e-14 & 2.6e-15 & 1.1e-15 & 2.2e-15 & 6.2e-16 & 1.0e-15 \\
100 & 1.2e-13 & 4.2e-15 & 1.5e-15 & 3.8e-15 & 7.2e-16 & 1.9e-15 \\
316 & 4.0e-13 & 7.4e-15 & 1.9e-15 & 7.0e-15 & 6.0e-16 & 3.9e-15 \\
1000 & 1.6e-12 & 1.3e-14 & 2.6e-15 & 1.2e-14 & 6.2e-16 & 8.4e-15 \\
	 \hline
	 \multicolumn{7}{||l||}{E3: $\lambda,\breve{\lambda} \sim \mathcal{U}(0,2),\ \theta\sim \mathcal{U}(0,\pi)$: Random normal scenario, $20\%$ real eigenvalues.}\\
	 \hline
	10 & 2.3e-15 & 1.4e-15 & 5.7e-16 & 1.3e-15 & 2.9e-16 & 4.9e-16 \\
32 & 1.3e-14 & 2.6e-15 & 1.2e-15 & 2.2e-15 & 5.8e-16 & 1.0e-15 \\
100 & 6.7e-14 & 4.1e-15 & 3.8e-15 & 3.7e-15 & 6.5e-16 & 1.9e-15 \\
316 & 2.7e-13 & 6.9e-15 & 1.2e-14 & 6.5e-15 & 5.8e-16 & 3.7e-15 \\
1000 & 8.7e-13 & 1.1e-14 & 2.9e-14 & 1.1e-14 & 5.7e-16 & 7.3e-15 \\
	 \hline
	 \multicolumn{7}{||l||}{E4: $\lambda\sim \mathcal{U}(0,2),\ \theta\sim \mathcal{U}(0,\pi)$: Random normal scen., $20\%$ repeated imaginary parts.}\\
	 \hline
	 10 & 3.5e-15 & 1.3e-15 & 6.6e-16 & 1.3e-15 & 3.4e-16 & 4.6e-16 \\
32 & 2.4e-14 & 2.0e-15 & 1.1e-15 & 1.9e-15 & 5.8e-16 & 7.5e-16 \\
100 & 1.2e-13 & 3.6e-15 & 1.5e-15 & 3.5e-15 & 6.9e-16 & 1.4e-15 \\
316 & 4.1e-13 & 7.4e-15 & 2.1e-15 & 7.9e-15 & 8.4e-16 & 4.2e-15 \\
1000 & 1.5e-12 & 2.3e-14 & 3.1e-15 & 2.4e-14 & 1.2e-15 & 1.9e-14 \\
 \hline
 \multicolumn{7}{||l||}{E5: $\lambda\sim \mathcal{U}(0,2),\ \theta\sim \pi\sqrt{\varepsilon_\mathrm{m}}\mathcal{N}(1,1)$: Worst-case normal scenario.}\\
	 \hline
	 %$10$&9.0e-10$\overset{\text{step 8}}{ \rightarrow \ }$\textcolor{green!70!blue!80!black!70}{ 1.2e-15}& 1.01e-15 & 4.86e-16\\
10 & 5.1e-10$\overset{\text{step 8}}{ \rightarrow \ }$\textcolor{Green}{1.2e-15} & 1.1e-15 & 9.9e-16 & 1.1e-15 & 5.3e-16 & 5.3e-16 \\
32 & 2.1e-9 $\rightarrow \ $\textcolor{Green}{2.9e-15}& 1.8e-15 & 2.2e-15 & 1.9e-15 & 1.2e-15 & 6.1e-16 \\
100 & 5.6e-9 $\rightarrow \ $\textcolor{Green}{5.8e-15} & 3.7e-15 & 4.1e-15 & 3.6e-15 & 2.9e-15 & 1.7e-15 \\
316 & 9.7e-9 $\rightarrow \ $\textcolor{Green}{1.4e-14}& 6.1e-15 & 6.8e-15 & 6.5e-15 & 3.5e-15 & 2.5e-15 \\
1000 & 1.1e-8 $\rightarrow \ $\textcolor{Green}{8.4e-14}& 9.7e-15 & 1.1e-14 & 1.0e-14 & 5.8e-15 & 3.8e-15 \\
	 \hline
	 \end{tabular}
	 \caption{\justifying Numerical experiment on the accuracy in double precision of \texttt{nrmschur} (\cref{alg:schur_decomposition_floating}) and \texttt{schur} (\texttt{LAPACK}'s \texttt{gees}) to compute $A\approx\widehat{Q}\widehat{S}\widehat{Q}^\top $. Step~8 is considered only for the last case where exactly one Jacobi sweep is performed. We chose $\delta,\delta_r,\frac{1}{t}$ equal to $\sqrt{\varepsilon_\mathrm{m}}$ to define clusters. $S$ is the ground truth real Schur form. Each reported accuracy is averaged over 100 runs. The normal matrices $A$ are randomly sampled with eigenvalues of different distributions specified in the table. The Schur vectors are always sampled with a Haar distribution on $\mathrm{O}(n)$. For all $a,b\in\mathbb{R}$, $a$\emph{e}-$b$  denotes $a\cdot10^{-b}$.}
	 \label{tab:accuracy}
	 \vspace*{-0.8cm}
	 \end{table}
In this section, we present the results of several numerical experiments on the accuracy of \cref{alg:schur_decomposition_floating} for a diverse set of $n\times n$ problems. The results are given in \Cref{tab:accuracy} and \cref{fig:evolution_accuracy}. For five different distributions of the eigenvalues, we report the average accuracy of  \cref{alg:schur_decomposition_floating} over $100$ runs. Codes reproducing the numerical experiments can be found at \url{https://github.com/smataigne/NormalEVP.jl}.
\begin{itemize}

\item Experiment E1 verifies that in the setting of \cref{thm:accuracy}, if the bound of \eqref{eq:final_result} remains small, then \cref{alg:schur_decomposition_floating} achieves high accuracy on the RSD. 

\item Experiment E2 tests a random normal scenario representing an ``average normal case''. In this case, a moderate loss of accuracy of the Schur decomposition occurs as $n$ grows because the separation between the imaginary parts of the eigenvalues decreases with $n$. However, the absolute accuracy on the eigenvalues remains close to machine precision thanks to the Rayleigh quotients.

\item Experiments E3 and E4 are similar to E2, but enforce nonempty sets \wordref{set2}{Set~2} and \wordref{set3}{Set~3}, respectively. The accuracy of the RSD remains equivalent to E2, as well as the absolute accuracy on the eigenvalues.

\item Finally, experiment E5 tests the worst case for \cref{alg:schur_decomposition_floating} in terms of accuracy. The skew-symmetric part $\Omega$ is small such that the eigenvalues of $\Omega$ are relatively close and induce sensitivity of the invariant subspaces. Moreover, $\Omega$ is large enough such that eigenvalues of $\Omega$ are not considered as $\sqrt{\varepsilon_\mathrm{m}}$-clusters and the accuracy is not restored by step~5. However, the high accuracy on the eigenvalues is maintained by the Rayleigh quotients. We show that one sweep of the Jacobi method~\cite{ZhouBrent03} allows to restore high accuracy of the RSD.
\end{itemize}

In \cref{tab:accuracy}, the deviation from orthogonality is measured by the simple criterion $\|\widehat{Q}^\top\widehat{Q}-I\|_\mathrm{F}$. We recall in \cref{lem:ortho_measures} that this measure is equivalent to $d_{\mathrm{O}(n)}(\widehat{Q})$ when it is small. \new{The results of \cref{fig:evolution_accuracy} show that for the experiments E2 to E4, the accuracy of the RSD grows proportionally to the matrix size $n$. Moreover, for the experiment E1, the accuracy is almost insensitive to $n$.}

	\section{Complexity analysis and running time performance}\label{sec:best_worst_case}
	We analyze the best and worst case computational complexities of~\cref{alg:schur_decomposition_floating}. We exclude the cases where the matrix $A$ is either purely skew-symmetric or symmetric since~\cref{alg:schur_decomposition_floating} reduces to obtaining, respectively, the skew-symmetric RSD with \wordref{oracle1}{Routine~1} or the symmetric EVD with \wordref{oracle2}{Routine~2}. These ideal cases are not relevant to our discussion. By analogy with algorithms such as quicksort, we will show that worst-case complexity analysis is not the most relevant here, as worst cases are very unlikely to occur.
	
	Since the cost of the method is $\mathcal{O}(n^3)$, the $\mathcal{O}(n^2)$ costs are neglected. Following the notation from \cite[Lec.~8]{trefethen97}, we say that a method requires $\sim g(n,p,r)$ flops if \begin{equation}
	\lim_{n,p,r\rightarrow \infty} \frac{\text{number of flops}}{g(n,p,r)}=1.
\end{equation} The numerical experiments presented in this section are performed on a processor \emph{Intel(R) Core(TM) i7-8750H CPU @ 2.20GHz} in single-threaded framework. Codes reproducing the numerical experiments can be found at \url{https://github.com/smataigne/NormalEVP.jl}.
%For completeness, we recall the following computational costs \cite[Sec.~1.2.4]{GoluVanl96}:
%\begin{itemize}
%\item The product $Ab$ with $A\in\mathbb{R}^{n\times p},\ b\in\mathbb{R}^p$ requires $\sim 2np$ flops.
%\item  The product $AB$ with $A\in\mathbb{R}^{n\times p},\ B\in\mathbb{R}^{p\times k}$ requires $\sim 2npk$ flops. 				\simcom{I am not convinced that recalling this is necessary}
%\end{itemize}
	\subsection{The intuitive (but not) best-case complexity}\label{subsec:best_complexity}
	Intuitively, the best case for \cref{alg:schur_decomposition_floating} is when all eigenvalues are complex and none of the imaginary parts are clustered (relatively to $\varepsilon_\mathrm{m}$). The assumptions of~\cref{thm:distinct_singular_values} are then matched with $k=\lfloor\frac{n}{2}\rfloor=p$ ($r=1$ for $n$ odd). We recall that for Haar-distributed special orthogonal matrices~\cite{Anderson87, Stewart80}, this case occurs with probability $1$~\cite{Fasi21}. %This is intuitive since, for every $A\in\mathrm{SO}(n)$, an arbitrarily small perturbation of the eigenvalues is enough to have $r=1$ and all imaginary parts distinct from each other.
	Step 1 is negligible since it requires $\mathcal{O}(n^2)$ flops.
	
	The cost of step 2 is $\sim	\frac{4}{3}n^3$ for the skew-symmetric tridiagonal reduction \cite[Lec.~26]{trefethen97}, plus an additional cost of $\sim \frac{4}{3}n^3$ to assemble the Householder reflectors \cite[Sec.~5.1.6]{GoluVanl96}. Indeed, given $M\in\mathbb{R}^{k\times k}$ and $v\in\mathbb{R}^k$, performing $(I-vv^\top )M$ requires $\sim 4k^2$ flops. Therefore, assembling all reflectors of size $k=1,...,n-1$ has a cost $\sim\sum_{k=1}^{n-1}4k^2\sim\frac{4}{3}n^3$ flops. 
	
	At step 3, the call to \wordref{oracle1}{Routine~1} can be neglected since an $\mathcal{O}(p^2)$ implementation exists \cite{GROER200345}. 
	
	Step 4 requires $\sim 4np^2$ flops, step 5 is skipped and step 6 costs $\sim 2n^2p$ flops. Step~7 is skipped (or costs $\mathcal{O}(n^2)$ if $r=1$). The total cost is thus
	\begin{equation}\label{eq:best_cost}
	\sim\frac{8}{3}n^3 + 4np^2  + 2n^2p \text{ flops}.
\end{equation}
Given $p=\lfloor\frac{n}{2}\rfloor$, the total cost is $\sim \frac{14}{3}n^3$, among which $\sim\frac{6}{3}n^3$ are level 3 \texttt{BLAS} operations \cite{blackford2002updated}. Moreover, the $\sim\frac{8}{3}n^3$ flops from the tridiagonalization also admit blocked implementations, maximizing the use of level~3 \texttt{BLAS} (see, e.g., \href{https://github.com/JuliaLinearAlgebra/SkewLinearAlgebra.jl}{SkewLinearAlgebra.jl}).
	
	In the best case, the total cost of~\cref{alg:schur_decomposition_floating} is thus equivalent to a Hessenberg factorization, which is known to cost $\frac{10}{3}n^3$ flops to obtain the $H$ factor \cite[Lec.~26]{trefethen97} and $\sim \frac{4}{3}n^3$ to assemble the Householder reflectors, hence $\sim \frac{14}{3}n^3$ flops in total. This result is \cref{alg:schur_decomposition_floating}'s best advantage over methods that include an initial reduction of the matrix $A$ to Hessenberg form \cite{AmmarReichel86,graggreichel90,wanggragg02,wanggragg01}, while the latter methods usually offer better accuracy on the eigenvectors in general.
	
	% In particular, random normal matrices belong to this category with probability $1$. %in the context of orthogonal matrices and statistics on manifolds, the matrices that deviate from the conditions of~\cref{thm:distinct_singular_values} are the exception.
	
\begin{figure}
		\center
		
		\includegraphics[width = 10cm]{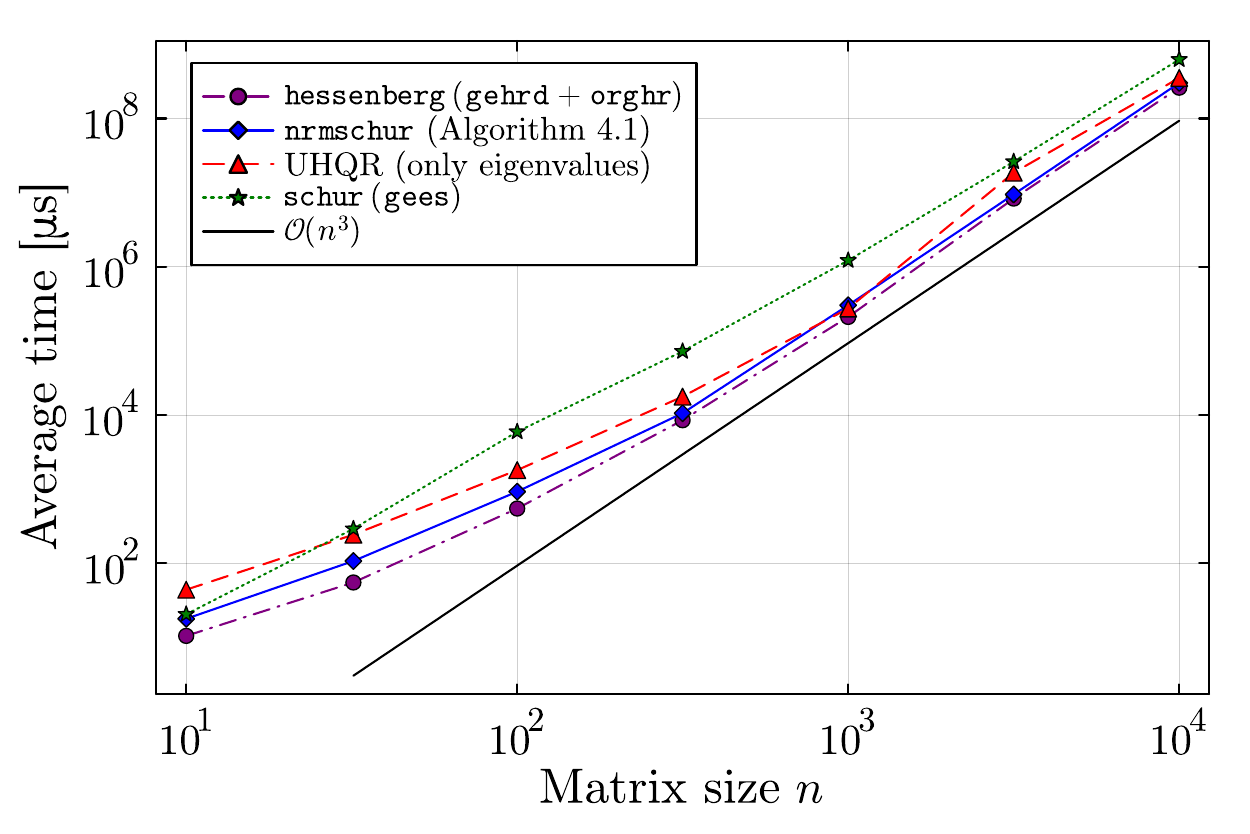}
		\includegraphics[width = 10cm]{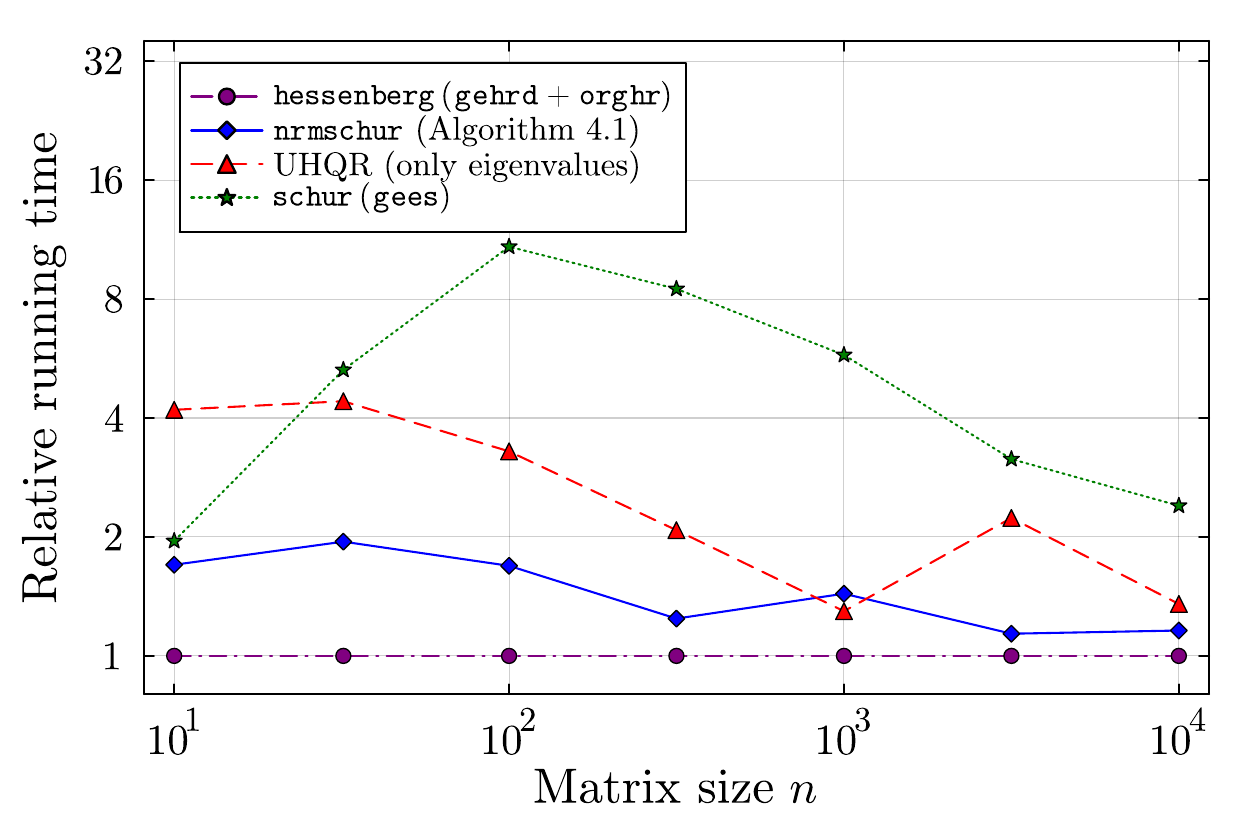}
		\vspace*{-0.4cm}
		\caption{Top: Benchmark of \texttt{nrmschur} (\cref{alg:schur_decomposition_floating}) on Haar-distibuted matrices of $\mathrm{SO}(n)$ w.r.t.~the \texttt{Julia} routines \href{https://github.com/JuliaLang/julia/blob/master/stdlib/LinearAlgebra/src/schur.jl}{\texttt{schur}} (\texttt{LAPACK}'s \href{https://www.netlib.org/lapack/explore-html/d5/d38/group__gees_gab48df0b5c60d7961190d868087f485bc.html}{\texttt{gees}}) and \href{https://github.com/JuliaLang/julia/blob/master/stdlib/LinearAlgebra/src/hessenberg.jl}{\texttt{hessenberg}} (\texttt{LAPACK}'s \href{https://www.netlib.org/lapack/explore-html/d2/d28/group__gehrd_ga74cea8f05a014cca243674999f71c238.html}{\texttt{gehrd}} and \href{https://www.netlib.org/lapack/explore-html/de/d26/group__unghr_gae1e32db855babbb134fd2a273aceed36.html}{\texttt{orghr}}). The performance of the UHQR algorithm for eigenvalues~\cite{GRAGG19861} as well as a representative $\mathcal{O}(n^3)$ are also drawn. The plot shows experiments for $n\in\{10,32,100,316,1000,3162,10000\}$. Bottom: the same results are displayed but divided by the computational time of \texttt{hessenberg}.}
		\label{fig:benchmark}
		\vspace*{-0.4cm}
	\end{figure}

	\Cref{fig:benchmark} illustrates the performance of \href{https://github.com/smataigne/NormalEVP.jl/blob/main/src/normal_schur.jl}{\texttt{nrmschur}} (\cref{alg:schur_decomposition_floating}) compared to \texttt{Julia} calls to \href{https://www.netlib.org/lapack/explore-html/index.html}{\texttt{LAPACK}}'s Hessenberg factorization (\href{https://www.netlib.org/lapack/explore-html/d2/d28/group__gehrd_ga74cea8f05a014cca243674999f71c238.html}{\texttt{gehrd}}+\href{https://www.netlib.org/lapack/explore-html/de/d26/group__unghr_gae1e32db855babbb134fd2a273aceed36.html}{\texttt{orghr}}) and \texttt{schur} routine (\href{https://www.netlib.org/lapack/explore-html/d5/d38/group__gees_gab48df0b5c60d7961190d868087f485bc.html}{\texttt{gees}}). We also include the performance of our \texttt{Julia} implementation of the UHQR algorithm \cite{GRAGG19861}. UHQR's computational complexity for the eigenvalues is $\sim \frac{10}{3}n^3$ for the Hessenberg reduction plus $\mathcal{O}(n^2)$ for QR iterations, thus faster than \texttt{LAPACK's} eigenvalue routine (\href{https://www.netlib.org/lapack/explore-html/d9/dc6/group__hseqr_ga62c3f96d2f67f96d6dc10334e118e451.html#ga62c3f96d2f67f96d6dc10334e118e451}{hseqr}). However, if the Schur/eigen-vectors are desired, the complexity of UHQR is $\sim \frac{14}{3}n^3$ for the Hessenberg reduction plus $\mathcal{O}(n^3)$  because Givens rotations of the QR are accumulated.
	
	%For Riemannian logarithms and Karcher means applications, we also draw the performance of \texttt{Julia}'s \href{https://github.com/JuliaLang/julia/blob/0af2f7407b39aee61385149475e849777c3aa85a/stdlib/LinearAlgebra/src/dense.jl}{\texttt{log}} implementation of the scaling and squaring algorithm for the matrix logarithm \cite{Higham05}. 
With due caution in interpreting complexity from timing experiments, the numerical results in \cref{fig:benchmark} corroborate that the cost of \cref{alg:schur_decomposition_floating} is comparable to that of a Hessenberg factorization. Moreover, it outperforms \texttt{gees} by a factor $10$ for $n\approx 100$ and a factor $2$ for $n=10000$. We can observe that UHQR (only for eigenvalues) lies above \cref{alg:schur_decomposition_floating} however, a \texttt{FORTRAN} implementation of UHQR would be faster. Nonetheless, this highlights that because \cref{alg:schur_decomposition_floating} relies exclusively on high performance software, namely \texttt{LAPACK}, it can be efficiently implemented in any computer language where an interface to this software is available.
	
	\subsection{The true best-case complexity}\label{subsec:oddity}
\begin{figure}[ht]
	\label{fig:alpha_figure}
	\center
	\includegraphics[width = 10cm]{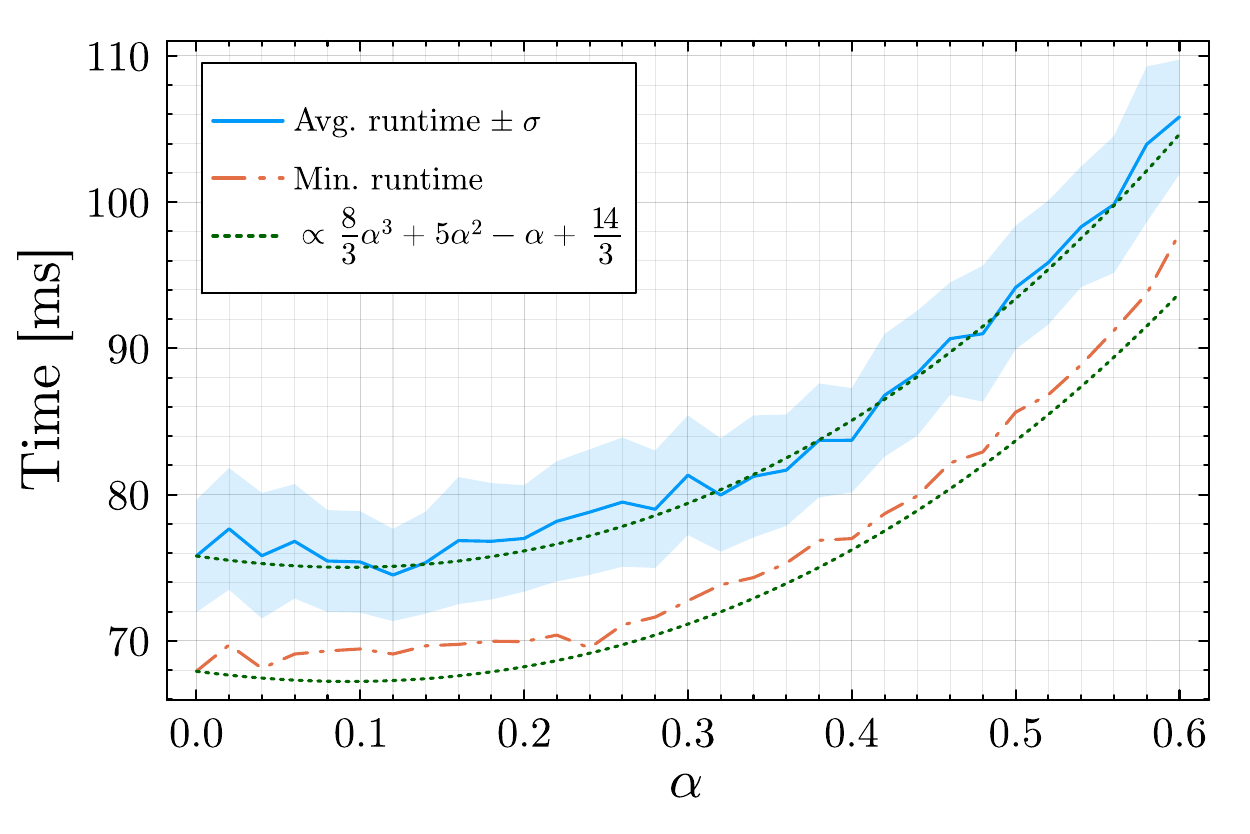}
	\vspace*{-0.4cm}
	\caption{Numerical experiments over $600\times 600$ random normal matrices with proportion $\alpha$ of real eigenvalues. For each value $\alpha$ ($\Delta \alpha=0.02$), 200 matrices are sampled. The Schur vectors are Haar-distributed on $\mathrm{O}(n)$, the real eigenvalues are normally distributed with unit variance and pairs of conjugate eigenvalues are sampled with uniformly distributed phases in $[0,\pi]$ and radii in $[0,2]$. The theoretical complexity from \eqref{eq:alpha_polynomial} is represented by the dotted green line. The polynomial $\frac{8}{3}\alpha^3+5\alpha^2-\alpha +\frac{14}{3}$ is scaled by a constant factor to match the intercept ($\alpha=0$) of the numerical experiments.}
	\vspace*{-0.4cm}
	\end{figure}

	If we consider $r>1$ real eigenvalues, then there is an additional cost of $\sim 2n^2r+2nr^2$ flops at step~7.1. The call to \wordref{oracle2}{Routine~2} at step~7.2a requires $\sim(\frac{4}{3}+\frac{4}{3})r^3$ flops for tridiagonalization plus $\mathcal{O}(r^2)$ flops to solve the tridiagonal EVP~\cite{dhillon1997new}. Finally, step~7.3 demands $\sim 2nr^2$ flops. However, notice that $p:=\frac{n-r}{2}$ decreases when $r$ rises. Hence, all costs from~\eqref{eq:best_cost} depending on $p$ decrease accordingly. For $r$ sufficently small, we show that the decrease of $p$ has a larger effect on the operation count than the increase of $r$. For $r>0$, a pessimistic (but realistic, see \cref{rem:low_rank}) total cost of~\cref{alg:schur_decomposition_floating} is 
	\begin{equation}\label{eq:cost_alpha}
	\sim \frac{8}{3}(n^3+r^3)+n(n-r)^2+n^2(n-r)+2n^2r+4nr^2 \text{ flops.}
	\end{equation}
	By setting $r=\alpha n$ with $\alpha\in[0,1]$, we obtain the complexity
	\begin{equation}\label{eq:alpha_polynomial}
		\sim \left(\frac{8}{3}\alpha^3+5\alpha^2-\alpha +\frac{14}{3}\right)n^3\text{ flops.}
	\end{equation}
	 The polynomial in $\alpha$ achieves a local minimum at $\alpha=\frac{\sqrt{5^2+8}-5}{8}\approx 0.093$ for an optimal cost of $4.61n^3$ flops, instead of $\sim\frac{14}{3}n^3\approx4.66n^3$ when $\alpha=0$. Although, the variation of the constant is too small to yield significant gain in performance, the take-home conclusion is that ``a few'' real eigenvalues (up to $20\%$) do not impact the performance negatively. The numerical experiment of \cref{fig:alpha_figure} confirms the relevance of this operation count. The average runtime fits the expected behavior when $\alpha$ varies.
	\begin{remark}\label{rem:low_rank}
\emph{When the number $r$ of real eigenvalues grows, exact arithmetic could use only $n-r$ Householder reflectors to tridiagonalize $\mathrm{skew}(A)$ instead of $n-1$. Indeed, since $\mathrm{skew}(A)$ has rank $n-r$ and since each column of a tridiagonal matrix with a non-zero subdiagonal element is independent from all previous columns, the (exact) tridiagonal form can not have more than $n-r$ non-zero subdiagonal elements. In short, tridiagonalization is a rank-revealing method. The most expensive reflectors being the $n-k$ ones of largest size, the worst case complexity in exact arithmetic of the tridiagonalization is  $\sim\frac{8}{3}(n^3-r^3)$ instead of $\sim\frac{8}{3}n^3$. However, due to numerical roundoff, logical zeros are perturbed. Given random matrices $U,V\in \mathbb{R}^{n\times l}$, $l\ll n$, it can be verified that \texttt{LAPACK}'s methods \href{https://www.netlib.org/lapack/explore-html/df/d34/group__hetrd_gac6c89000449da7196e0273911fe492b2.html}{\texttt{sytrd}} on $UU^\top $ and \href{https://www.netlib.org/lapack/explore-html/d2/d28/group__gehrd_ga74cea8f05a014cca243674999f71c238.html}{\texttt{gehrd}} on $UV^\top $ produce $n-1$ non-trivial reflectors. This validates in practice the cost $\sim\frac{8}{3}n^3$ considered in~\eqref{eq:cost_alpha}.}
	\end{remark}
\subsection{The worst-case complexity}\label{subsec:worst_complexity}
	The worst case for \cref{alg:schur_decomposition_floating} is when all eigenvalues are complex (all but one if $n$ is odd) and their imaginary parts form two conjugate  clusters. This case maximizes the additional work of \wordref{oracle3}{Routine~3} at step~5. For simplicity, assume $n$ even such that $r=0$ and at step~5, $\Sigma =\sigma I_p$ such that $V=\widehat{Q}$. Essentially, all steps but step 5.2 are a waste a computational resources.
	
Assuming we ignore the (potentially inaccurate) symmetric skew-Hamiltonian structure, step 5.1 requires $\sim 4n^3$ flops. Step 5.2 requires $\sim \frac{14}{3}n^3$ and $\mathcal{O}(n^3)$ flops for, respectively, the Hessenberg reduction and \wordref{oracle3}{Routine~3}. Step~5.3 demands $\sim 2n^3$ flops. \cref{alg:schur_decomposition_floating} returns after step~5 since $\Lambda\cos(\Theta)$ was computed at step~5.2. Taking into account the $\sim \frac{8}{3}n^3+\mathcal{O}(n^2)+n^3$ flops of steps 2 to 4, the total cost is $\sim \frac{43}{3}n^3+\mathcal{O}(n^3)$, thus a large increase of the computational cost.

\begin{figure}[ht]
	\center
	\includegraphics[width = 10cm]{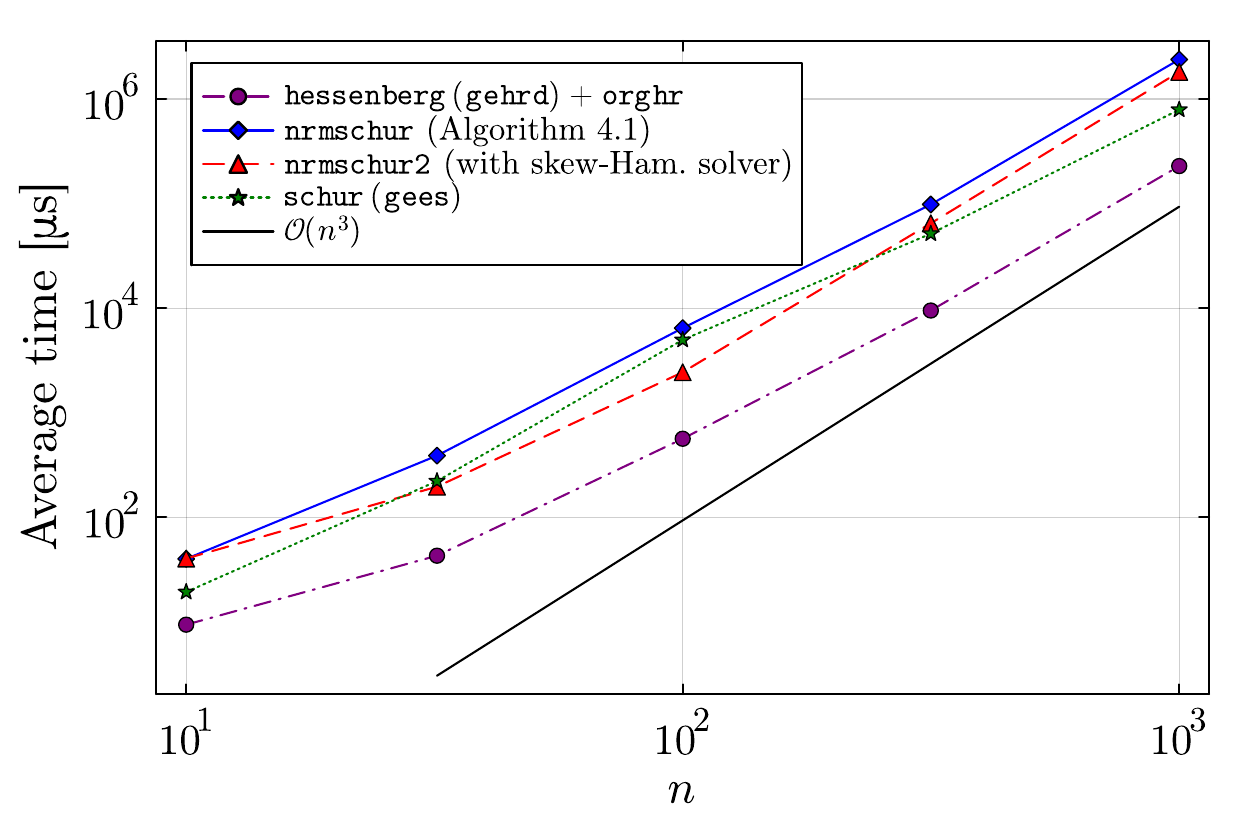}
	\vspace*{-0.4cm}
	\caption{Worst-case benchmark, as in \cref{fig:benchmark}, with $n\times n$ random normal matrices such that $\Sigma= \sigma I_p$, $\sigma\in\mathbb{R}$ and $D_p\sim\mathcal{N}(0,1)$. The Schur vectors are sampled with Haar-distribution on $\mathrm{O}(n)$~\cite{Stewart80}. The plot shows $n\in\{10,32,100,316,1000\}$.}
	\label{fig:worst_benchmark}
	% Failure of \texttt{LAPACK}'s bidiagonal SVD \texttt{dbdsdc} at $n=3334$ prevented numerical experiments at this size (LAPACKException(1)). Similar failures were also observed with \texttt{dstemr}.}
	%\vspace*{-0.4cm}
	\end{figure}

	\Cref{fig:worst_benchmark} reports an experiment that is similar to the best case (\cref{fig:benchmark}) but for the worst case. Without surprise, the average runtime of \cref{alg:schur_decomposition_floating} exceeds that of \texttt{gees} since \texttt{gees} becomes a subroutine of~\cref{alg:schur_decomposition_floating}. The practical advantage of~\cref{alg:schur_decomposition_floating} stems from the unlikelihood of having $\Sigma=\sigma I_p$. Most likely, the multiplicity of the entries of $\Sigma$, if any, remains small compared to~$n$. Then, the cost of step~5 remains negligible and the performance does not significantly differ from the $\sim\frac{14}{3}n^3$ flops of the best case.

	 \section{Application: Riemannian barycenter on manifolds} \label{sec:karcher_mean} The Riemannian barycenter is a classical problem from statistics on manifolds, in particular on $\mathrm{SO}(n)$ and $\mathrm{St}(n,p)$, see, e.g., \cite{BINI2013,krakowski07,LIM2012,YING2016,ZimmermannHuper22}. For both manifolds, intensive computations of matrix logarithms of orthogonal matrices are required, as detailed next for $\mathrm{SO}(n)$. Given a collection of $N$ matrices $\{X_i\}_{i=1}^N$ where $X_i\in\mathrm{SO}(n)$ for $i=1,...,N$, one is interested in finding a Riemannian center of mass $X_\mathrm{c}\in\mathrm{SO}(n)$ of the collection, i.e., a matrix satisfying
	 \begin{equation}\label{eq:karcher_mean}
	 	X_\mathrm{c}\in\mathrm{arg}\hspace{-4pt}\min_{\hspace{-12pt} X\in \mathrm{SO}(n)}\frac{1}{2}\sum_{i=1}^N d(X_i, X)^2,
	 \end{equation}
	 where $d(X_\mathrm{a},X_\mathrm{b})\coloneq\|\log(X_\mathrm{a}^\top X_\mathrm{b})\|_\mathrm{F}$ is the Riemannian distance function on $\mathrm{SO}(n)$. Well-posedness of the problem, existence and uniqueness of $X_\mathrm{c}$ are discussed in \cite{Afsari11}. The definition \eqref{eq:karcher_mean} extends the usual notion of center of mass in the Euclidean space. To solve \eqref{eq:karcher_mean}, we consider the Riemannian gradient descent \cite{Afsari13, krakowski07}, i.e., for $k=1,2,...$ we compute
	 \begin{equation}\label{eq:gradient_step}
	 X_\mathrm{c}^{k+1} = X_\mathrm{c}^{k}\exp\left( - \frac{1}{N} \sum_{i=1}^N \log(X_i^\top X_\mathrm{c}^k)\right).
	 \end{equation}
	 To test the performance of \cref{alg:schur_decomposition_floating}, we perform $100$ iterations of Riemannian gradient descent using~\eqref{eq:gradient_step}. We perform a fixed number of iterations instead of taking a convergence criterion to evaluate more clearly the impact of the RSD routine. The Riemannian gradient step is computed using the relation $\log(A)=Q\log(S)Q^\top $ (see \eqref{eq:log_normal}). The RSD is computed using either \texttt{schur} (\texttt{LAPACK}'s \texttt{gees}) or \texttt{nrmschur} (\cref{alg:schur_decomposition_floating}). We vary the size $n$ of $\mathrm{SO}(n)$, and the number $N$ of samples in the collection. The results are displayed in \Cref{tab:karcher}. The larger  $n$ and $N$, the larger the improvement. Indeed, the factor of improvement in running time varies from $1.72$ in the smallest configuration to $4.38$ for $(N=64, n=200)$.

	 \begin{table}[ht]
	 \centering
	 %\scriptsize
	 %\scalefont{0.8}
	 %\fontsize{8.5}{10}\selectfont
	 \begin{tabular}{||c||c c||c c||c c||c c||}
	 %\hline
	 %& \multicolumn{8}{c||}{$n=$}\\
	 \hline
	  \diagbox[dir=SE]{$N$}{$n$}&\multicolumn{2}{c||}{$25$}&\multicolumn{2}{c||}{$50$}&\multicolumn{2}{c||}{$100$}&\multicolumn{2}{c||}{$200$}\\
	  \hline
	  & \texttt{schur}&\texttt{nrmschur}& \texttt{schur}&\texttt{nrmschur}& \texttt{schur}&\texttt{nrmschur}&\texttt{schur}&\texttt{nrmschur}\\
	 \hline
	  $16$& $0.307$&$ 0.178$ &  $1.31$& $0.663$& $10.1 $&$ 2.56$& $42.9$&$11.6$ \\
	  $32$&$0.570$&$0.327$&$  2.90$&$1.35$& $ 21.0$&$4.97$&$  77.4$&$17.8 $ \\
	 $64$&$1.02$&$0.527$ & $ 5.87$&$ 2.23$&$ 43.2$&$9.56$&$ 185$&$42.2$ \\
	 \hline
	 \end{tabular}
	 \caption{\justifying Time in seconds to perform $100$ steps of Riemannian Gradient Descent (RGD) for the Riemannian barycenter problem. The sizes $n$ of $\mathrm{SO}(n)$ varies from $25$ to $200$ and the number of samples from $16$ to~$64$. Orthogonal EVDs are respectively computed with \texttt{schur} and \texttt{nrmschur}.}
	 \label{tab:karcher}
	 \vspace*{-0.7cm}
	 \end{table}
	%\rotatebox{90}{$N=$}
	 
	 \section{Conclusion} In this paper, we proposed a new method for computing the real Schur decomposition of a dense real normal matrix, or equivalently, its eigenvalue decomposition. We showed that the skew-symmetric part of the matrix is a valuable tool for obtaining this decomposition, as its two-dimensional invariant subspaces coincide with those of the original matrix.

Our approach was the following. First, we developed the matrix theory results that motivate the algorithm. This naturally led to an exact arithmetic version of the method, which we detail in \cref{alg:schur_decomposition}. Next, we examined the method’s robustness under floating point arithmetic, resulting in the floating point variant \cref{alg:schur_decomposition_floating}.

By analyzing the sensitivity of the invariant subspaces of the skew-symmetric part, we theoretically quantified in \cref{sec:accuracy} and numerically corroborated in \cref{sec:exp_accuracy} how an unfavorable distribution of eigenvalues can reduce the numerical accuracy of the real Schur decomposition if no additional computations are performed. While this loss of accuracy may often be negligible in applications such as Riemannian computations, we demonstrated how it can be mitigated by applying correction steps to selected inaccurate invariant subspaces, using the bounds provided in \cref{thm:accuracy,thm:accuracy_repeated,thm:accuracy_real}. The numerical accuracy of \cref{alg:schur_decomposition_floating} was verified in \cref{sec:exp_accuracy}.

We then analyzed the computational complexity of the method. In particular, we showed that in scenarios such as Haar-distributed orthogonal matrices, the cost of \cref{alg:schur_decomposition_floating} is comparable to that of performing a Hessenberg reduction followed by assembling the Householder reflectors, i.e., $\sim \frac{14}{3}n^3$.

Finally, we highlighted the relevance of this method in the context of statistics on manifolds, and computations on the Stiefel manifold and the orthogonal group, where the Riemannian logarithm plays a central role and requires an efficient eigenvalue decomposition.

\appendix

\section{Special invariance groups}\label{sec:similarities} In this section, we provide important lemmas about invariance groups and we characterize two of them in \cref{lem:Invariance_lambda,lem:invariance_Im}.
\begin{lemma}\label{lem:permute}

	For every $A\in\mathbb{R}^{n\times n}$, all $Q\in\mathrm{ig}(A)$ and all $V\in\mathrm{O}(n)$, we have $VQV^\top \in\mathrm{ig}(VAV^\top )$. In particular, $V$ can be a permutation matrix.
	\end{lemma}
	\begin{proof}
	Since $QAQ^\top =A$ and $V^\top V=I_n$, we have $(VQV^\top )(VAV^\top )(VQV^\top )=(VAV^\top )$.
	\end{proof}
	\begin{lemma}
	\label{lem:Invariance_lambda}
		Let $\Sigma\in\mathbb{R}^{p\times p}$ be a diagonal matrix with distinct non-zero diagonal values. An orthogonal matrix $Q\in\mathrm{O}(2p)$ is such that
		\begin{equation}\label{eq:invariance_sigma}
			Q\begin{bmatrix}
			0&-\Sigma\\
			\Sigma& 0
			\end{bmatrix}Q^\top  = \begin{bmatrix}
			0&-\Sigma\\
			\Sigma& 0
			\end{bmatrix},
		\end{equation}
		if and only if $Q = \begin{bmatrix}
		\cos(\Phi)&-\sin(\Phi)\\
		\sin(\Phi)&\cos(\Phi)
		\end{bmatrix}$ for some diagonal matrix $\Phi\in\mathbb{R}^{p\times p}$.
	\end{lemma}
	\begin{proof}
	We proceed by induction. Assume first $p=1$, i.e., $\Sigma=\sigma \in \mathbb{R}$. Then, for every $Q\in\mathrm{O}(2)$, there is $\phi\in\mathbb{R}$ such that either  $Q=\left[\begin{smallmatrix}
		\cos(\phi)&-\sin(\phi)\\
		\sin(\phi)&\cos(\phi)
		\end{smallmatrix}\right]$ ($Q$ is a rotation) or $Q=\left[\begin{smallmatrix}
		\cos(\phi)&\sin(\phi)\\
		\sin(\phi)&-\cos(\phi)
		\end{smallmatrix}\right]$ ($Q$ is a reflection). It is easily verified that only rotations belong to $\mathrm{ig}\left(\left[\begin{smallmatrix}
		0&-\sigma\\
		\sigma&0
		\end{smallmatrix}\right]\right)$. Thus, \cref{lem:Invariance_lambda} holds for $p=1$. 
		
		Assume \cref{lem:Invariance_lambda} holds for $p\geq 1$ and let $\Sigma\in\mathbb{R}^{p\times p}$, $S_1 = \left[\begin{smallmatrix}
			0&-\Sigma\\
			\Sigma& 0
			\end{smallmatrix}\right]$ and $S_2 = \left[\begin{smallmatrix}
			0&-\sigma\\
			\sigma& 0
			\end{smallmatrix}\right]$ with $\sigma\in\mathbb{R}$ be non-zero and distinct from all entries of $\Sigma$. Let $Q=\left[\begin{smallmatrix}
			Q_1&Q_2\\
			Q_3&Q_4
			\end{smallmatrix}\right]\in\mathrm{O}(2p+2)$ with $Q_1\in\mathbb{R}^{2p\times 2p}$ and $Q_4\in\mathbb{R}^{2\times 2}$. Then, if one imposes\begin{equation*}
			Q\begin{bmatrix}
			S_1&0\\
			0&S_2
			\end{bmatrix}Q^\top =\begin{bmatrix}
			S_1&0\\
			0&S_2
			\end{bmatrix} \text{ and } Q\in \mathrm{O}(2p+2),
			\end{equation*}
			it must hold in particular that
			\begin{equation*}
			\begin{cases}
			Q_3 S_1 Q_3^\top  + Q_4 S_2 Q_4^\top  = S_2,\\
			Q_3 S_1 Q_1^\top  + Q_4 S_2 Q_2^\top  = 0,\\
			\end{cases} \text{ and } \begin{cases}
			Q_1^\top Q_1 + Q_3^\top Q_3 = I_{2p},\\
			Q_1^\top Q_2 + Q_3^\top Q_4 = 0.
			\end{cases}
			\end{equation*}
			By multiplying the second equation by $Q_1$ on the right, we obtain $Q_3 S_1 Q_1^\top Q_1 + Q_4 S_2 Q_2^\top Q_1 = 0$. Leveraging both orthogonality constraints, we obtain $Q_3S_1(I_{2p}-Q_3^\top Q_3)-Q_4S_2Q_4^\top Q_3=0$. Introducing the first equation, we obtain $Q_3S_1=S_2Q_3$. Let $Q_3:=[Q_{31}\ Q_{32}]$ such that \begin{align*}
		Q_3S_1=S_2Q_3 \iff Q_{32}\Sigma = S_2 Q_{31} &\text{ and } \ -Q_{31}\Sigma= S_2Q_{32}.\\
		\iff Q_{31}\Sigma^2 = \sigma^2 Q_{31} &\text{ and } \ Q_{32}\Sigma^2 = \sigma^2 Q_{32}.
\end{align*}			 
 Since $\sigma$ is distinct from all entries of $\Sigma$, $Q_3=0$. Hence, $Q_4\in\mathrm{O}(2)$ and $Q_2=0$. Finally, $Q\in \mathrm{ig}\left(\left[\begin{smallmatrix}
			S_1&0\\
			0&S_2
			\end{smallmatrix}\right]\right)$ if and only if it is block-diagonal, $Q_1\in\mathrm{ig}(S_1)$ and $Q_4\in\mathrm{ig}(S_2)=\mathrm{SO}(2)$.
			
			To conclude, we can define $\widetilde{\Sigma}:=\mathrm{diag}(\Sigma,\sigma)$ and take a permutation matrix $P$ such that $P\left[\begin{smallmatrix}
			S_1&0\\
			0&S_2
			\end{smallmatrix}\right]P^\top =\left[\begin{smallmatrix}
			0&-\widetilde{\Sigma}\\
			\widetilde{\Sigma}&0
			\end{smallmatrix}\right]$. Applying \cref{lem:permute} is enough to conclude that \cref{lem:Invariance_lambda} holds for $p+1$. Therefore, it holds for all $p\geq 1$.
	\end{proof}
	\begin{comment}
	\begin{proof}
	\simcom{Proof to be finished}
		Let us divide $Q$ by blocks, $Q:=\begin{bmatrix}
		A&C\\
		B&D
		\end{bmatrix}$. The orthogonality of $Q$ yields, among others, the following conditions
		\begin{equation}\label{eq:cond_orthogonality}
			\begin{cases}
				C^\top C+D^\top D = I_p,\\
				A^\top C+B^\top D=0.
			\end{cases}
		\end{equation}
		Moreover, the condition \eqref{eq:invariance_sigma} yields
		\begin{equation}
			\begin{cases}
				C\Sigma A^\top  - A\Sigma C^\top  = 0,\\
				D\Sigma B^\top  - B\Sigma D^\top =0,\\
				D\Sigma A^\top  - B\Sigma C^\top  = \Sigma.
			\end{cases}
		\end{equation}
	\end{proof}
	
	Multiplying the second equation of \eqref{eq:cond_orthogonality} by $C\Sigma$ on the left, we get
	\begin{align*}
		C\Sigma A^\top C + C\Sigma B^\top D &= 0\\
		\iff A\Sigma C^\top C +(A\Sigma D^\top -\Sigma)D &=0\\
		\iff A\Sigma C^\top C +A\Sigma(I_p-C^\top C)-\Sigma D &=0\\
		\iff A\Sigma  &= \Sigma D.
	\end{align*}
	Now, multiplying the second equation of \eqref{eq:cond_orthogonality} by $D\Sigma$ on the left, we get
	\begin{align*}
		D\Sigma A^\top C + D\Sigma B^\top D &= 0\\
		\iff (\Sigma + B\Sigma C^\top )C +B\Sigma D^\top  D &=0\\
		\iff \Sigma C + B\Sigma (I-D^\top D)+B\Sigma D^\top  D&=0\\
		\iff \Sigma C  &= -B\Sigma.
	\end{align*}
	\end{comment}
	\begin{lemma}
		\label{lem:invariance_Im}
		An orthogonal matrix $Q = \begin{bmatrix}
		A&C\\
		B&D
		\end{bmatrix}\in\mathrm{O}(2m)$, $A\in\mathbb{R}^{m\times m}$, is such that 
		\begin{equation}\label{eq:invariance}
		Q\begin{bmatrix}
		0&-I_m\\
		I_m&0
		\end{bmatrix}Q^\top =\begin{bmatrix}
		0&-I_m\\
		I_m&0
		\end{bmatrix},
		\end{equation}
		if and only if $A=D$ and $B = -C$. The matrix $Q$ is called \emph{ortho-symplectic}~\cite{dopico09}.
	\end{lemma}
	\begin{proof}
	It is enough to notice that \eqref{eq:invariance} is equivalent to
	\begin{equation}
		Q\begin{bmatrix}
		0&-I_m\\
		I_m&0
		\end{bmatrix}=\begin{bmatrix}
		0&-I_m\\
		I_m&0
		\end{bmatrix}Q\iff \begin{bmatrix}
		C&-A\\
		D&-B
		\end{bmatrix} = \begin{bmatrix}
		-B&-D\\
		A&C
		\end{bmatrix}.
		\end{equation}
		Therefore $Q\in\mathrm{O}(2m)$, $A=D$ and $B=-C$ are necessary and sufficient conditions.
	\end{proof}
\new{Finally, we briefly recall the relation between the distance to the orthogonal group $d_{\mathrm{O}(n)}(\widehat{Q})$ from \eqref{eq:dist_to_orthogonality} and $\|Q^\top Q-I\|_\mathrm{F}$. These measures are equivalent for $d_{\mathrm{O}(n)}(\widehat{Q})$ small.
	\begin{lemma}\label{lem:ortho_measures}
For all matrices $\widehat{Q}\in\mathbb{R}^{n\times n}$, it holds that
\begin{equation*}
d_{\mathrm{O}(n)}(\widehat{Q})\leq \|\widehat{Q}^\top\widehat{Q}-I\|_\mathrm{F}\leq 2d_{\mathrm{O}(n)}(\widehat{Q}) + d_{\mathrm{O}(n)}(\widehat{Q})^2.
\end{equation*}
\end{lemma}	
\begin{proof}
Let $\Sigma$ be the diagonal matrix of the singular values of $\widehat{Q}$. Then, it is well-known that $ d_{\mathrm{O}(n)}(\widehat{Q})=\|\Sigma-I\|_\mathrm{F}$. Moreover, $\|\widehat{Q}^\top\widehat{Q}-I\|_\mathrm{F} = \|\Sigma^2-I\|_\mathrm{F}$.  It follows that $ \|\Sigma^2-I\|_\mathrm{F} = \|(\Sigma-I)(\Sigma+I)\|_\mathrm{F}\geq \|\Sigma-I\|_\mathrm{F}$. Moreover, 
\begin{equation*}
	\|\Sigma^2 - I\|_\mathrm{F}= \|(\Sigma - I +I)^\top(\Sigma - I +I) - I\|_\mathrm{F} \leq 2\|\Sigma - I\|_\mathrm{F} + \|\Sigma - I\|_\mathrm{F}^2.
\end{equation*}
This concludes the proof.
\end{proof}}
	\begin{comment}
	\begin{proof}
	We provide here one of the possible proofs.
		Orthogonality of $Q$ and \eqref{eq:invariance} yield the following conditions
		\begin{equation}
		\begin{cases}
			C^\top C + D^\top D =I_m\\
			DA^\top  - BC^\top  = I_m\\
			A^\top C + B^\top D =0\\
			CA^\top -AC^\top  = 0\\
			DB^\top  - BD^\top  = 0
		\end{cases}
		\end{equation}
		Multiplying the third equation by $C$ on the left, we get
		\begin{align*}
			CA^\top C + CB^\top D &=0\\
			\iff AC^\top C +(AD^\top  - I_m)D &= 0\\
			\iff A(I-D^\top D)+(AD^\top  - I_m)D &= 0\\
			\iff A &= D
		\end{align*}
		Moreover, multiplying the third equation by $D$ on the left, we get
		\begin{align*}
			DA^\top C + DB^\top D &=0\\
			\iff(I_m+BC^\top )C +BD^\top D &= 0\\
			\iff (I_m+BC^\top )C+B(I_m - C^\top C) &= 0\\
			\iff C &= -B
		\end{align*}
		The conditions $A=D$ and $B=-C$ are thus necessary for \eqref{eq:invariance} to hold. One can also easily check that it is also sufficient, so that the conditions $A=D$ and $B=-C$ are necessary and sufficient conditions.
	\end{proof}
	\end{comment}
	\section{Symmetric skew-Hamiltonian matrices}\label{app:sym_skew_hamiltonian} In this section, we highlight two important properties of the symmetric skew-Hamiltonian matrix $\left[\begin{smallmatrix}
	\widetilde{H}&-\widetilde{\Omega}\\
	\widetilde{\Omega}&\widetilde{H}
	\end{smallmatrix}\right]$ from \cref{thm:equal_singular_values}. This matrix appears when a normal matrix has complex eigenvalues with repeated imaginary parts.
	\begin{lemma}\label{lem:symplectic_vectors}
	A matrix $A\in\mathbb{R}^{2m\times 2m}$ is symmetric skew-Hamitonian, i.e., $A=\left[\begin{smallmatrix}
	\widetilde{H}&-\widetilde{\Omega}\\
	\widetilde{\Omega}&\widetilde{H}
	\end{smallmatrix}\right]$ with $\widetilde{H}$ symmetric and $\widetilde{\Omega}$ skew-symmetric, if and only if there is an eigenvalue decomposition $A=Q\left[\begin{smallmatrix}
	\Lambda&0\\
	0&\Lambda
	\end{smallmatrix}\right]Q^\top $ where $Q$ is ortho-symplectic and $\Lambda$ is diagonal.
	\end{lemma}
	\begin{proof}
	($\Longrightarrow$) Assume $A=\left[\begin{smallmatrix}
	\widetilde{H}&-\widetilde{\Omega}\\
	\widetilde{\Omega}&\widetilde{H}
	\end{smallmatrix}\right]$. Since $A$ is symmetric, all eigenvalues are real. Moreover, $\left[\begin{smallmatrix}
	v_1\\
	v_2
	\end{smallmatrix}\right]$ is a normed eigenvector of $A$ if and only if $\left[\begin{smallmatrix}
	-v_2\\
	v_1
	\end{smallmatrix}\right]$ is an eigenvector of~$A$. Indeed, 
	\begin{equation}\label{eq:even_multiplicity}\footnotesize
	\begin{bmatrix}
	\widetilde{H}&-\widetilde{\Omega}\\
	\widetilde{\Omega}&\widetilde{H}
	\end{bmatrix}\begin{bmatrix}
	v_1\\
	v_2
	\end{bmatrix}=\lambda\begin{bmatrix}
	v_1\\
	v_2
	\end{bmatrix}\iff \begin{cases}
	(\widetilde{H}-\lambda I)v_1 = \widetilde{\Omega}v_2,\\
	(\widetilde{H}-\lambda I)v_2 = -\widetilde{\Omega}v_1.
	\end{cases} \hspace{-0.3cm}\iff \begin{bmatrix}
	\widetilde{H}&-\widetilde{\Omega}\\
	\widetilde{\Omega}&\widetilde{H}
	\end{bmatrix}\begin{bmatrix}
	-v_2\\
	v_1
	\end{bmatrix}=\lambda\begin{bmatrix}
	-v_2\\
	v_1
	\end{bmatrix}.
	\end{equation}
	Therefore, the multiplicity of every eigenvalue is even and greater or equal to two. Since \eqref{eq:even_multiplicity} holds for every eigenvector, the eigenvector's matrix always admits a form $Q=\left[\begin{smallmatrix}E&-F\\
	F&E\end{smallmatrix}\right]$ with $A=Q\left[\begin{smallmatrix}
	\Lambda&0\\
	0&\Lambda
	\end{smallmatrix}\right]Q^\top $. By \cref{lem:invariance_Im}, $Q$ is ortho-symplectic.
	
	($\Longleftarrow$) Since $Q$ is ortho-symplectic, $Q=\left[\begin{smallmatrix}E&-F\\
	F&E\end{smallmatrix}\right]$ by \cref{lem:invariance_Im}. Multiplying $Q\left[\begin{smallmatrix}
	\Lambda&0\\
	0&\Lambda
	\end{smallmatrix}\right]Q^\top $ such as in \eqref{eq:symmetric_subproblem} concludes the proof.
	\end{proof}
\begin{theorem}\label{thm:symplectic_tridiag}
For every symmetric skew-Hamiltonian matrix $A\in\mathbb{R}^{2m\times 2m}$, i.e., $A=\left[\begin{smallmatrix}
	\widetilde{H}&-\widetilde{\Omega}\\
	\widetilde{\Omega}&\widetilde{H}
	\end{smallmatrix}\right]$ with $\widetilde{H}$ symmetric and $\widetilde{\Omega}$ skew-symmetric, if $\mathrm{rank}(A)=2m$ and every eigenvalue of $A$ has multiplicity $2$, then the $m$ first steps of Lanczos tridiagonalization are enough in exact arithmetic to build an ortho-symplectic matrix $Y$ where $A = Y\left[\begin{smallmatrix}
	T_m&0\\
	0&T_m
	\end{smallmatrix}\right]Y^\top $ and $T_m$ is an $m\times m$ tridiagonal symmetric matrix.
	\end{theorem}
	\begin{proof}
	We denote $K=[b\ Ab\ ...A^{2m-1}b]$ a Krylov basis of $A$ generated by $b$. Choose any $v\in\mathbb{R}^n$ such that $b\coloneq Av\neq 0$ and $K_{1:m}$ has column rank $m$.\footnote{Notice that $\mathrm{rank}(K)\leq m$ since every eigenvalue has multiplicity $2$ by \cref{lem:symplectic_vectors}.} By \cref{lem:symplectic_vectors}, we know that there is an ortho-symplectic matrix $Q$ such that $A=Q\left[\begin{smallmatrix}
	\Lambda&0\\
	0&\Lambda
	\end{smallmatrix}\right]Q^\top $ and $\Lambda\in\mathbb{R}^{m\times m}$ is diagonal. Let $k_j$ denote the $j$th column of $K$, it follows that \begin{equation*}
		k_j = Q\begin{bmatrix}
	\Lambda^j&0\\
	0&\Lambda^j
	\end{bmatrix}Q^\top v.
	\end{equation*}
	Moreover, for any integer pair $i,j$, with $1\leq i,j\leq n$, we have
	\begin{align}
	\nonumber
		k_i^\top \begin{bmatrix}
		0&-I_m\\
		I_m&0
		\end{bmatrix}k_j &=v^\top Q\begin{bmatrix}
	\Lambda^i&0\\
	0&\Lambda^i
	\end{bmatrix}Q^\top \begin{bmatrix}
		0&-I_m\\
		I_m&0
		\end{bmatrix}Q\begin{bmatrix}
	\Lambda^j&0\\
	0&\Lambda^j
	\end{bmatrix}Q^\top v\\
	\nonumber
	&=v^\top Q\begin{bmatrix}
	\Lambda^i&0\\
	0&\Lambda^i
	\end{bmatrix}\begin{bmatrix}
		0&-I_m\\
		I_m&0
		\end{bmatrix}\begin{bmatrix}
	\Lambda^j&0\\
	0&\Lambda^j
	\end{bmatrix}Q^\top v\\
	\nonumber
	&=v^\top Q\begin{bmatrix}
		0&-\Lambda^{i+j}\\
		\Lambda^{i+j}&0
		\end{bmatrix}Q^\top v\\
\label{eq:symplectic_orthogonality}
		&=0.
	\end{align}
	The last equality stands because $v^\top \Omega v=0$ for every skew-symmetric matrices $\Omega$. Letting $J_{2m}\coloneq\left[\begin{smallmatrix}
		0&-I_m\\
		I_m&0
		\end{smallmatrix}\right]$, it follows from \eqref{eq:symplectic_orthogonality} that $K_{1:m}^\top J_{2m}K_{1:m}=0$. Therefore, the matrix $[K_m\ \ J_{2m}K_{1:m}]$ is full rank. Consequently, the Gram-Schmidt QR factorization of $K_{1:m}$, namely, $K_{1:m} = Y_{1:m}R$, generates $Y\coloneq[Y_{1:m}\ \ J_{2m} Y_{1:m}]$ where $Y$ is ortho-symplectic. Moreover, we also have 
	\begin{equation}\label{eq:biorthogonality}
	k_i^\top  A J_{2m} k_j = k_{i+1}^\top  J_{2m}k_j = 0, 
\end{equation}
and thus $K_{1:m}^\top AJ_{2m}K_{1:m} = 0$. Finally, it is the praised property of Lanczos algorithm that $Y_{1:m}^\top AY_{1:m} = T_m$ where $T_m$ is tridiagonal and symmetric \cite[Thm.~9.1.1]{GoluVanl96}. Combining \eqref{eq:symplectic_orthogonality} and \eqref{eq:biorthogonality}, we obtain
		\begin{align*}
		Y^\top AY &= \begin{bmatrix}
		Y_{1:m}^\top \\
		Y_{1:m}^\top J_{2m}^\top 
		\end{bmatrix}A [Y_{1:m} \ J_{2m}Y_{1:m}]\\
		&=\begin{bmatrix}
		Y_{1:m}^\top AY_{1:m}&Y_{1:m}^\top AJ_{2m}Y_{1:m}\\
		Y_{1:m}^\top J_{2m}^\top AY_{1:m}&Y_{1:m}^\top J_{2m}^\top AJ_{2m}Y_{1:m}
		\end{bmatrix}\\
		& = \begin{bmatrix}
		T_m&0\\
		0&T_m
		\end{bmatrix}.
		\end{align*}
		We used that $J_{2m}^\top AJ_{2m}=A$. This concludes the proof.
	\end{proof}
\cref{alg:symplectic_lanczos} is a simple modification of the usual modified Gram-Schmidt to build a matrix that is not only orthogonal, but also symplectic. In floating point arithmetic, the matrix $Y$ is approximately orthogonal but \emph{exactly} symplectic by construction. Instead of orthogonalizing any matrix, \cref{alg:symplectic_lanczos} orthogonalizes the Krylov basis generated by $A$ symmetric skew-Hamiltonian, we directly obtain Lanczos's procedure with full re-orthogonalization. 
	
\begin{figure}
	\centering
	\includegraphics[width = 8.1cm]{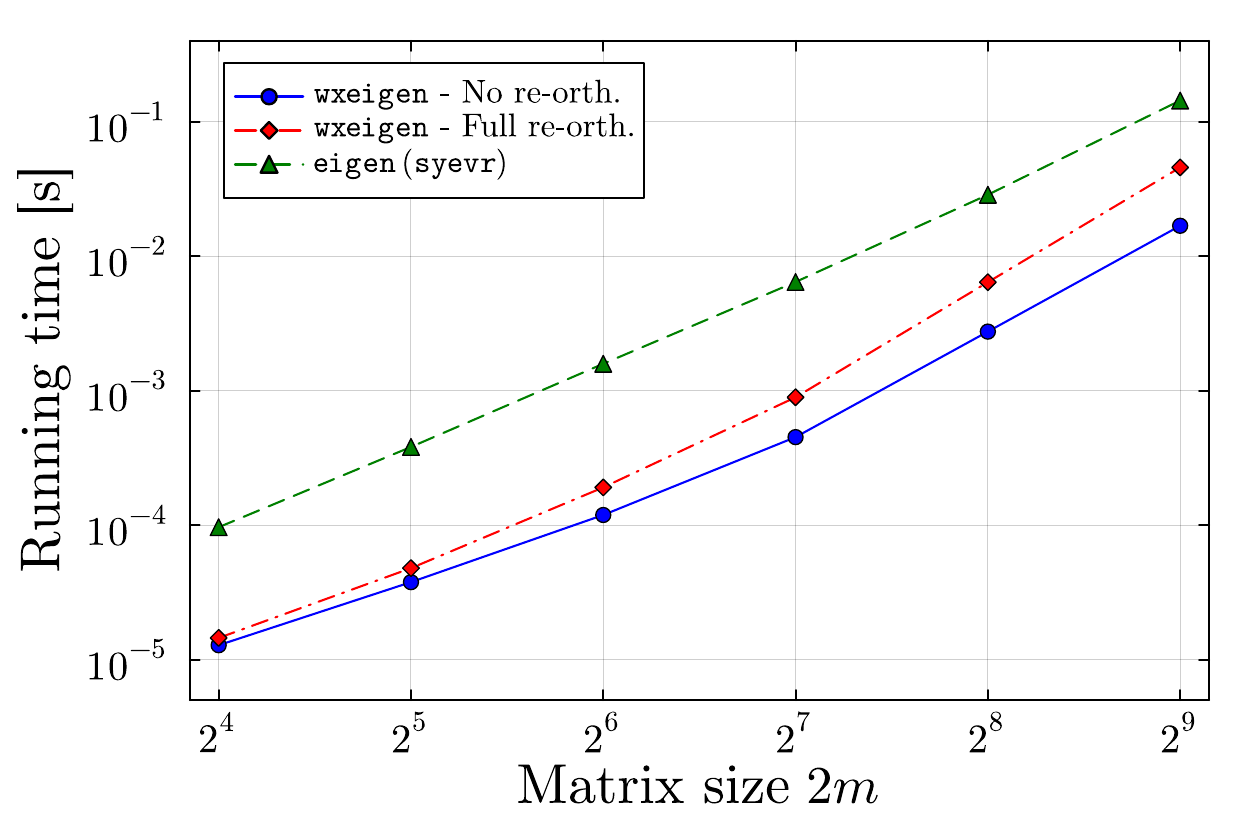}
	\includegraphics[width = 8.1cm]{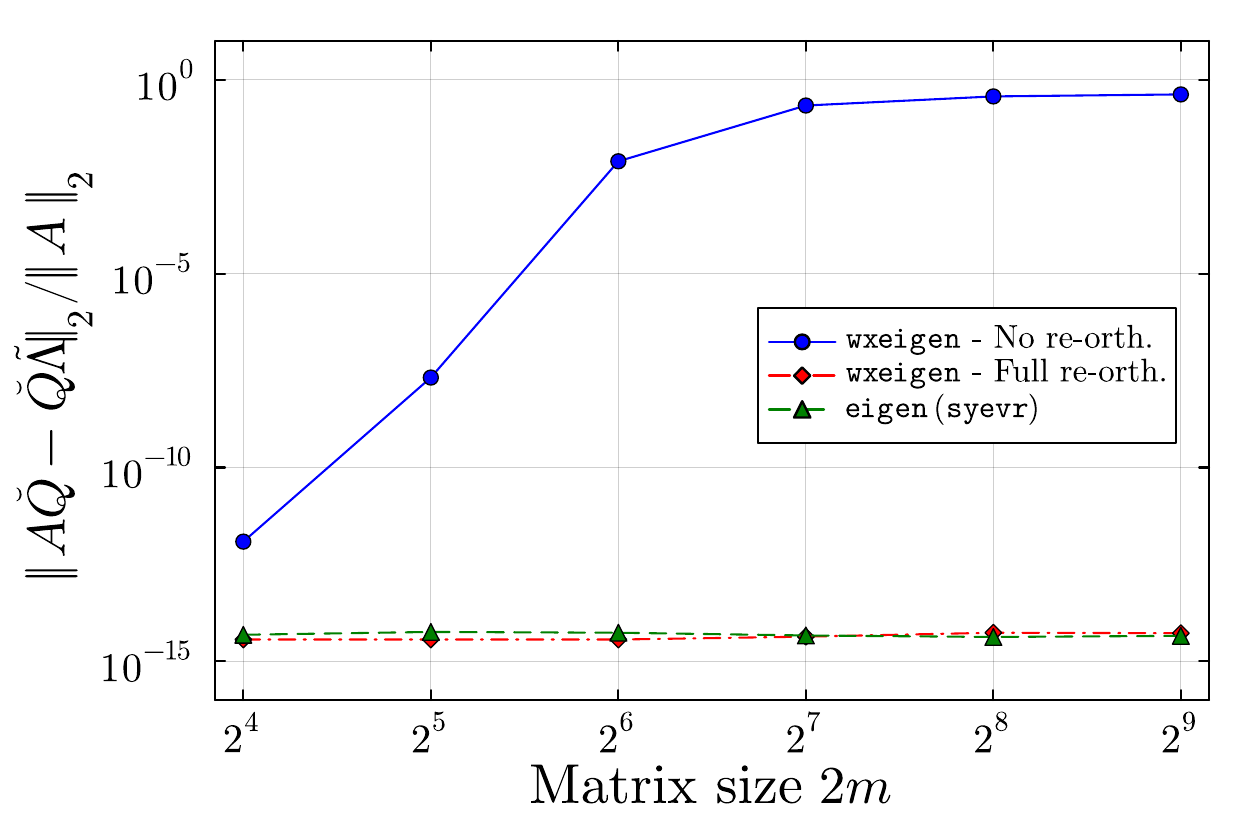}
	\caption{Numerical performance evaluation of \texttt{wxeigen}, without and with full re-orthogonalization of the Lanczos vectors, compared to \texttt{eigen} (\texttt{LAPACK}'s \texttt{syevr}). The plot on the left shows the running time performance. The plot on the right shows the accuracy performance in spectral norm. Taking advantage of the symmetric skew-Hamiltonian structure allows to design faster methods than the generic symmetric eigensolver.}
	\label{fig:WXbenchmark}
	\end{figure}

In \cref{fig:WXbenchmark}, we evaluate how much \cref{thm:symplectic_tridiag} can improve the performance compared to \texttt{LAPACK}'s symmetric eigenvalue routine for small EVP's. The method using only $m$ steps of Lanczos tridiagonalization is called \texttt{wxeigen} and is available at \href{https://github.com/smataigne/NormalEVP.jl/blob/main/src/wxeigen.jl}{wxeigen.jl} in the repository. We observe that \texttt{wxeigen} with full re-orthogonalization achieves good running time performance and accuracy.
\begin{algorithm}
		\caption{Re-orthogonalized Lanczos for symmetric skew-Hamiltonian matrix.}
		\begin{algorithmic}\label{alg:symplectic_lanczos}
		\STATE \textbf{INPUT}: $v\in\mathbb{R}^{2m}$, $A\in\mathbb{R}^{2m\times 2m}$.
		\STATE $\beta_0=0$, $y_1 = Av/\|Av\|_2$ and $y_{m+1} = J_{2m} y_1$.
		\FOR{$i=2,...,m$}
			\STATE $y_i=Ay_{i-1}$ and $\alpha_{i-1} = y_{i-1}^\top y_i$.
			\STATE $y_i = y_i - \alpha_{i-1}y_{i-1} - \beta_{i-2}y_{i-2}$. \hfill Lanczos step.
			\FOR{$j=1,..,i-1$}
				\STATE $y_i = y_i - (y_i^\top y_j)y_{j} - (y_i^\top y_{m+j})y_{m+j}$. \hfill Full re-orthogonalization.
			\ENDFOR
			\STATE $\beta_i=\|y_i\|_2$.
			\STATE $y_i = y_i/\beta_i$ and $y_{m+i} = J_{2m} y_i$.\hfill Build symplectic basis.
		\ENDFOR
		\STATE $\alpha_m = y_mAy_m$.
		\RETURN $Y\in\mathrm{O}(2m)\cap \mathrm{OSp}(2m)$, $\{\alpha_i\}_{i=1}^m$ and $\{\beta_i\}_{i=1}^{m-1}$.
		\end{algorithmic}
	\end{algorithm}

\section{Proof of \cref{thm:accuracy_real}}\label{app:theorem_proof}

\begin{proof}
 Let $A=QSQ^\top \in\mathbb{R}^{n\times n}$ be a normal matrix with real eigenvalue $\breve{\lambda}_k$ for $k\in\mathcal{K}\coloneq\{1,...,r\}$. Moreover, assume $\lambda_j e^{\pm i\theta_j}$ is an eigenvalue of $A$ for $j \in\mathcal{I}$. Finally, assume we have computed $\widehat{V}\in\mathrm{St}(n,r)$ such that $\Omega\widehat{V}=E$ with $\|E\|_\mathrm{F}\leq\tau \varepsilon_\mathrm{m} \|A\|_\mathrm{F}$. We have
\begin{equation}\label{eq:matrix_equation}
	\Omega\widehat{V}=E\iff \mathrm{skew}(S)Q^\top \widehat{V} =Q^\top E\eqcolon\widetilde{E}.
\end{equation}
Defining $X\coloneq Q^\top\widehat{V}$, there are two permutation $P_\mathrm{l}\in\mathrm{O}(n),P_\mathrm{r}\in\mathrm{O}(r) $ such that $P_\mathrm{l}^\top X P_\mathrm{r}=\left[\begin{smallmatrix}X_{\mathcal{I},\mathcal{K}}  \\X_{\mathcal{K},\mathcal{K}}\end{smallmatrix}\right]$ where the block $X_{\mathcal{I},\mathcal{K}}$ comprises all $2\times 1$ blocks $X_{ik}$ such that for $i\in\mathcal{I}, k\in\mathcal{K}$, by \eqref{eq:matrix_equation}, we have
\begin{equation*}
\begin{bmatrix}
		0&-\lambda_i s_i\\
		\lambda_i s_i&0
	\end{bmatrix}X_{ik}=(P_\mathrm{l}^\top \widetilde{E} P_\mathrm{r})_{ik}.
\end{equation*} 
Since $X\in\mathrm{St}(n,r)$, we have
\begin{align*}
 X_{\mathcal{I},\mathcal{K}}^\top X_{\mathcal{I},\mathcal{K}} +X_{\mathcal{K},\mathcal{K}}^\top X_{\mathcal{K},\mathcal{K}} = I_{r} \Longrightarrow \| X_{\mathcal{K},\mathcal{K}}^\top X_{\mathcal{K},\mathcal{K}}-I_{r}\|_\mathrm{F}\leq \| X_{\mathcal{I},\mathcal{K}}\|_\mathrm{F}^2\\
 \text{and,}\quad \| X_{\mathcal{I},\mathcal{K}}\|_\mathrm{F}^2= \sum_{i\in\mathcal{I}, k\in\mathcal{K}} \frac{1}{|\lambda_i s_i|^2}\|(P_\mathrm{l}^\top \widetilde{E}  P_\mathrm{r})_{ik}\|_\mathrm{F}^2\leq \frac{\|E\|_\mathrm{F}^2}{\min_{i\in\mathcal{I}}|\lambda_i s_i|^2},
\end{align*}
where \cite[Thm.~VII.2.8]{bhatia97} is used. The deviation of $X_{\mathcal{K},\mathcal{K}}$ from orthogonality is quantified by that of $X_{\mathcal{I},\mathcal{K}}$ from zero. In addition, by taking the SVD $ X_{\mathcal{K},\mathcal{K}}=U_X\Sigma_X V_X^\top$, we can write $\| X_{\mathcal{K},\mathcal{K}}^\top X_{\mathcal{K},\mathcal{K}}-I_{r}\|_\mathrm{F}=\|\Sigma_X^2-I_{r}\|_\mathrm{F}\geq \|\Sigma_X-I_{r}\|_\mathrm{F}$ by \cref{lem:ortho_measures}. Define $\left[\begin{smallmatrix}\widetilde{S}_{n-r} &0 \\ 0 &\breve{\Lambda}\end{smallmatrix}\right]\coloneq P_\mathrm{l}^\top S P_\mathrm{l}$. For all $\breve{R}\in \mathrm{O}(r)$, we have 
\begin{align}
\nonumber
	\|A\widehat{V}\breve{R}-\widehat{V}\breve{R}\breve{\Lambda}\|_\mathrm{F}&\leq\|\mathrm{sym}(S)XR-XR\breve{\Lambda}\|_\mathrm{F}+\|E\|_\mathrm{F}\\
	\label{eq:newmaster_equation}
	&\leq \|\mathrm{sym}(\widetilde{S}_{n-r})X_{\mathcal{I},\mathcal{K}}\breve{R}-X_{\mathcal{I},\mathcal{K}}\breve{R}\breve{\Lambda}\|_\mathrm{F}\\
	\nonumber
	&\quad +\|\breve{\Lambda}X_{\mathcal{K},\mathcal{K}}\breve{R}-X_{\mathcal{K},\mathcal{K}}\breve{R}\breve{\Lambda}\|_\mathrm{F}+\|E\|_\mathrm{F}.
	%\nonumber
	%\|\widetilde{S}_{2m}X_{\mathcal{J},\mathcal{J}}\breve{R}-X_{\mathcal{J},\mathcal{J}}\breve{R}\widetilde{S}_{2m}\|_\mathrm{F}^2+\|\widetilde{S}_{2k}X_{\mathcal{I},\mathcal{J}}\breve{R}-X_{\mathcal{I},\mathcal{J}}\breve{R}\widetilde{S}_{2m}\|_\mathrm{F}^2
	%\label{eq:master_equation}
	%&+\sum_{i\in\mathcal{I},j\in\mathcal{J}} \left\|\lambda_i G_i(XR)_{ij}-(XR)_{ij}\lambda_j G_j \right\|_\mathrm{F}^2.
\end{align}
The first term of \eqref{eq:master_equation} can be bounded by
\begin{align}
\nonumber
	\|\mathrm{sym}(\widetilde{S}_{n-r})X_{\mathcal{I},\mathcal{K}}\breve{R}-X_{\mathcal{I},\mathcal{K}}\breve{R}\breve{\Lambda}\|_\mathrm{F}&\leq  \max_{i\in\mathcal{I},k\in\mathcal{K}}|\lambda_i c_i-\breve{\lambda}_k|\|X_{\mathcal{I},\mathcal{K}}\breve{R}\|_\mathrm{F}\\
	\label{eq:newslave1}
	&\leq \frac{\max_{i\in\mathcal{I},k\in\mathcal{K}}|\lambda_i c_i-\breve{\lambda}_k|}{\min_{i\in\mathcal{I}}|\lambda_i s_i|}\|E\|_\mathrm{F}.
\end{align}
Moreover, the second term can also be bounded by
\begin{align}
\nonumber
\|\breve{\Lambda}X_{\mathcal{K},\mathcal{K}}\breve{R}-X_{\mathcal{K},\mathcal{K}}\breve{R}\breve{\Lambda}\|_\mathrm{F}&= \|\breve{\Lambda}(X_{\mathcal{K},\mathcal{K}}\breve{R}-I_{r})-(X_{\mathcal{K},\mathcal{K}}\breve{R}-I_{r})\breve{\Lambda}\|_\mathrm{F}\\
\nonumber%\label{eq:newslave2}
&\leq  \max_{k_1,k_2\in\mathcal{K}}|\breve{\lambda}_{k_1}-\breve{\lambda}_{k_2}|\|X_{\mathcal{K},\mathcal{K}}\breve{R}-I_{r}\|_\mathrm{F}.
\end{align}
In particular, for $\breve{R}=V_X U_X^\top$, which minimizes $\|X_{\mathcal{K},\mathcal{K}}\breve{R}-I_{r}\|_\mathrm{F}$, we have 
\begin{equation}\label{eq:newslave3}
\|X_{\mathcal{K},\mathcal{K}}\breve{R}-I_{r}\|_\mathrm{F}=\|\Sigma_X-I_{r}\|_\mathrm{F}\leq \frac{\|E\|_\mathrm{F}^2}{\min_{i\in\mathcal{I}}|\lambda_i s_i|^2}.  
\end{equation}
Inserting \eqref{eq:newslave1} and \eqref{eq:newslave3} in \eqref{eq:newmaster_equation} results in \cref{thm:accuracy_real}.
\end{proof}
%%%%%%%%%%%%%%%%%%%%%%%%%%%%%%%%%%%%%%%%%%%%%%%%%%%%%%%%%%%%%

\bigskip
{\bf Acknowledgment.} The authors would like to thank Paul Van Dooren and Pierre-Antoine Absil for their invaluable advice during the preparation of this manuscript.

\end{document}